\documentclass[a4paper]{amsart}

\usepackage{a4wide}
\usepackage{hyperref}

\newcommand{\dd}{\;{\rm d}}
\newcommand{\de}{{\rm d}}
\newcommand{\Lp}{\mathcal{L}}
\newcommand{\Sbb}{\mathbb{S}}
\newcommand{\pshi}{\varphi}
\renewcommand{\epsilon}{\varepsilon}
\DeclareMathOperator{\Leb}{Leb}

\newcommand{\Z}{\mathbb{Z}}
\newcommand{\N}{\mathbb{N}}
\newcommand{\norm}[1]{\left\| #1 \right\|}
\newcommand{\E}{\mathcal{E}}
\newcommand{\R}{\mathbb{R}}
\newcommand{\C}{\mathbb{C}}
\newcommand{\CC}{\mathcal{C}}
\newcommand{\B}{\mathcal{B}}
\newcommand{\D}{\mathbb{D}}
\DeclareMathOperator{\Lip}{Lip} \DeclareMathOperator{\Id}{Id}
\DeclareMathOperator{\type}{type}
\newcommand{\moins}{-}
\newcommand{\AAA}{\mathcal{A}}
\newcommand{\HH}{\mathcal{H}}
\newcommand{\Mp}{\mathcal{M}}
\newcommand{\TT}{\mathcal{T}}
\newcommand{\UU}{\mathcal{U}}
\newcommand{\KK}{\mathcal{J}}
\newcommand{\BB}{\mathcal{B}}
\newcommand{\st}{\;|\;}
\newcommand{\given}{\ |\ }
\renewcommand{\Id}{I}
\newcommand{\K}{\mathcal{K}}
\DeclareMathOperator{\Ide}{Id}
\newcommand{\kk}{\underline{k}}
\DeclareMathOperator{\End}{End}
\newcommand{\boN}{\mathcal{N}}
\newcommand{\dist}{d}
\newcommand{\PP}{\mathcal{F}}
\newcommand{\SL}{\textrm{SL}(2,\R)}
\newcommand{\Pbb}{\mathbb{P}}
\newcommand{\Cone}{\mathbf{C}}
\newcommand{\eps}{\tau_0}
\DeclareMathOperator{\diam}{diam}
\newcommand{\tC}{\bar C}
\newcommand{\tD}{\bar D}
\newcommand{\TY}{\overline{T}}
\DeclareMathOperator{\Card}{Card}

\newcommand{\Ca}{C_2}
\newcommand{\Cb}{C_3}
\newcommand{\Cc}{C_1}
\newcommand{\Cd}{C_4}

\newtheorem{thm}{Theorem}[section]
\newtheorem{lem}[thm]{Lemma}
\newtheorem{prop}[thm]{Proposition}
\newtheorem{cor}[thm]{Corollary}
\newtheorem{definition}[thm]{Definition}
\newtheorem{rmq}[thm]{Remark}

\numberwithin{equation}{section}

\title[Local limit theorem and Farey sequences]
{Local limit theorem for nonuniformly partially hyperbolic
skew-products, and Farey sequences}

\date{March 22, 2007}

\author{S\'{e}bastien Gou\"{e}zel}
\address{
IRMAR, Universit\'{e} de Rennes 1, Campus de Beaulieu, B\^{a}timent 22,
35042 Rennes Cedex, France. }

\email{sebastien.gouezel@univ-rennes1.fr}

\begin{document}

\begin{abstract}
We study skew-products of the form $(x,\omega)\mapsto (Tx,
\omega+\phi(x))$ where $T$ is a nonuniformly expanding map on a
space $X$, preserving a (possibly singular) probability measure
$\tilde\mu$, and $\phi:X\to \Sbb^1$ is a $C^1$ function. Under
mild assumptions on $\tilde\mu$ and $\phi$, we prove that such
a map is exponentially mixing, and satisfies the central and
local limit theorems. These results apply to a random walk
related to the Farey sequence, thereby answering a question of
Guivarc'h and Raugi.
\end{abstract}

\maketitle

\section{Results}

Let $\TT$ be a transformation on a compact manifold. If $\TT$
is uniformly expanding or hyperbolic, the transfer operator
associated to $\TT$ admits a spectral gap on a well chosen
Banach space, which makes it possible to prove virtually any
limit theorem (for example the local limit theorem) by using
Nagaev's method (see e.g.\ \cite{guivarch-hardy,
hennion_herve}). This article is devoted to the proof of the
local limit theorem for transformations of the form
$\TT:(x,\omega) \mapsto (Tx, \omega+\phi(x))$ where $T$ is a
nonuniformly expanding transformation on a compact manifold
$X$, and $\phi: X \to\Sbb^1$ is a $C^1$ function. This
transformation $\TT$ is an isometry in the fibers $\Sbb^1$,
which prevents us from obtaining a spectral gap.

Limit theorems have been obtained (in the more general setting
of partially hyperbolic transformations) by Dolgopyat in
\cite{dolgopyat:limit} (when $T$ is uniformly hyperbolic, and
for a measure which is absolutely continuous with respect to
Lebesgue measure in the unstable direction). However, he uses
elementary arguments (moment methods) which can not be used to
get the local limit theorem. To the best of our knowledge, the
only partially hyperbolic transformations for which a local
limit theorem is proved in the literature are the Anosov flows,
in  \cite{waddington} (the specific algebraic structure of
flows makes it possible to reduce the problem to the study of
Axiom A maps, which are uniformly hyperbolic). With the
techniques of \cite{tsujii:semiflow}, it is probably possible
to obtain it also for skew-products over uniformly expanding
maps, for an absolutely continuous measure. Unfortunately, the
main motivating example of our study, described in the next
paragraph, is nonuniformly hyperbolic, and its invariant
measure is singular. Hence, we will need to introduce a new
technique, essentially based on renewal theory.

The qualitative theory of skew-products as above has been
studied by Brin. We will need more quantitative results, and
will obtain them by using tools which are mainly due to
Dolgopyat \cite{dolgopyat:decay, dolgopyat:compact_extensions}.
These techniques of Dolgopyat have already proved very powerful
in a variety of contexts (see \cite{pollicott_sharp:Dolgo,
anantharaman:Dolgo, stoyanov:eclipse, naud:Dolgo,
baladi_vallee:ContinuedFractions, baladi_vallee, AGY:teich}),
the present paper is yet another illustration of their
usefulness.

\subsection{Farey sequences}
\label{subsec:Farey}

Before we give the precise definition of the systems to which
our results apply, let us describe an interesting example,
which is in fact the main motivation for this article. The
following discussion is essentially taken from
\cite{conze_guivarch}.

If $p/q$ and $p'/q'$ are two irreducible rational numbers in
$[0,1]$, they are \emph{adjacent} if $|pq'-p'q|=1$. We can then
construct their median $p''/q''=(p+p')/(q+q')$, which lies
between $p/q$ and $p'/q'$, and is adjacent to any of them. Let
$\PP_0=\{ 0/1, 1/1\}$, and define inductively $\PP_n$ by
enumerating the elements of $\PP_{n-1}$ in increasing order,
which gives a sequence of adjacent rational numbers, and by
inserting the successive medians. For example, $\PP_1=\{0/1,
1/2, 1/1\}$ and $\PP_2=\{0/1, 1/3, 1/2, 2/3, 1/1\}$. The set
$\PP_n$ has cardinality $2^n+1$. Let also $\PP_n^*=\PP_n \moins
\{0\}$, it has cardinality $2^n$. Any rational number of
$(0,1]$ belongs to $\PP_n^*$ for any large enough $n$. Let
$\mu_n= \frac{1}{2^n}\sum_{x\in \PP_n^*} \delta_x$, this
sequence of measures converges exponentially fast to a measure
$\mu$, in the following sense: for any $\alpha>0$, there exist
$C>0$ and $\theta<1$ such that, for any function $f:[0,1]\to
\C$ which is H\"{o}lder continuous of exponent $\alpha$,
  \begin{equation}
  \label{eq:ConvergenceExpMarcheAleatoire}
  \left| \int f \dd\mu_n -\int f\dd\mu \right| \leq C \theta^n
  \norm{f}_{C^\alpha}.
  \end{equation}
The measure $\mu$ is \emph{Minkowski's measure}, it has full support
in $[0,1]$ and is totally singular with respect to Lebesgue measure.
It is the Stieltjes measure associated to Minkowski's $?$ function.

To prove the exponential convergence
\eqref{eq:ConvergenceExpMarcheAleatoire}, it is more convenient
to reformulate everything in terms of a random walk on a
homogeneous space for the group $\SL$. Consider the two
matrices $A=\left(\begin{smallmatrix} 1&0\\1&1
\end{smallmatrix}\right)$ and $B=\left(\begin{smallmatrix}
0&1\\-1&2 \end{smallmatrix}\right)$ in $\SL$. Their linear
action on $\R^2$ leaves invariant the cone $\Cone=\{ (x,y) \st
0\leq x \leq y\}$, and its projectivization $\Pbb(\Cone)$ is
the unique closed subset of $\Pbb(\R^2)$ which is invariant and
minimal for the action of the semigroup $\Sigma$ generated by
$A$ and $B$. Let us identity $\Pbb(\Cone)$ with the interval
$[0,1]$ by intersecting $\Cone$ with the line $y=1$, we obtain
an action of $\Sigma$ on $[0,1]$. The actions of the matrices
$A$ and $B$ are given by the transformations
  \begin{equation}
  h_A(x)=\frac{x}{1+x}, \quad h_B(x)=\frac{1}{2-x}.
  \end{equation}
It can easily be checked inductively that
  \begin{equation}
  \PP_n^*=\{ M_n\cdots M_1\cdot 1 \st M_i\in \{A,B\}\text{ for
  }i=1,\dots,n\}.
  \end{equation}
In particular, setting $\nu=(\delta_A+\delta_B)/2$, we have $\mu_n=
\nu^n \star\delta_1$. The measure $\mu$ is the unique stationary
measure for the random walk given by $\nu$, i.e., such that
$\nu\star\mu=\mu$. Finally, the exponential convergence
\eqref{eq:ConvergenceExpMarcheAleatoire} is proved by showing that
the Markov operator associated to the random walk has a spectral gap
when it acts on the space of H\"{o}lder continuous functions.

In \cite{conze_guivarch} (see also \cite{guivarch_raugi}),
Conze and Guivarc'h have considered the same random walk, but
on homogeneous spaces which are larger than $\Pbb(\R^2)$. More
precisely, let us fix $r>1$, and consider the quotient of
$\R^2\moins\{0\}$ by the subgroup $H_r$ of homotheties of ratio
$\pm r^n,\ n\in \Z$. This is a compact space, endowed with an
action of $\SL$. In particular, the semigroup $\Sigma$ acts on
$\bar \Cone=\Cone/H_r$, which is a compact extension (with
fiber $\Sbb^1$) of $\Pbb(\Cone)$. Let us identify $\bar \Cone$
with $[0,1] \times \R/ (\log r)\Z$ by $(x,y) \mapsto (x/y, \log
y + (\log r)\Z)$. The random walk given by $\nu$ on $\bar
\Cone$ jumps from $(x,\omega)$ to $\bar h_A(x,\omega):=(h_A(x),
\omega+ \log (1+x))$ or $\bar h_B(x,\omega):=(h_B(x),
\omega+\log(2-x))$ with probability $1/2$. Let $\bar
\PP_n^*=\{( p/q, \log q) \st p/q \in \PP_n^*\} \subset [0,1]
\times \R/ (\log r)\Z$, the measure $\bar\mu_n:=\nu^n \star
\delta_{(1,0)}$ is the average of the Dirac masses at the
points of $\bar\PP_n^*$. Hence, the random walk given by $\nu$
and starting from the point $(1,0)$ describes the rational
numbers obtained by the Farey process, as well as the logarithm
of their denominators, modulo $\log r$. By general results on
random walks on compact extensions, Conze, Guivarc'h and Raugi
proved in \cite{conze_guivarch, guivarch_raugi} that
$\bar\mu_n$ converges weakly to $\mu \otimes \Leb$, where
$\Leb$ denotes the normalized Lebesgue measure on $\R/(\log
r)\Z$. This is an equirepartition result of the denominators
modulo $\log r$.

In this article, we are interested in more precise results for this
random walk. First of all, we prove that the previous convergence is
exponentially fast:
\begin{thm}
\label{thm:MelangeExponentielFarey} For any $\alpha>0$, there exist
$C>0$ and $\theta<1$ such that, for any function $f:\bar \Cone\to
\C$ which is H\"{o}lder-continuous of exponent $\alpha$,
  \begin{equation}
  \left| \int f \dd\bar\mu_n -\int f\dd(\mu\otimes \Leb) \right| \leq C \theta^n
  \norm{f}_{C^\alpha}.
  \end{equation}
\end{thm}

We also obtain limit theorems for this random walk. In
particular, we prove that it satisfies the local limit theorem.
This answers a question raised by Guivarc'h and Raugi in
\cite{guivarch_raugi}.

\begin{thm}
\label{thm:ThmLimiteFarey} Let $\psi:\bar \Cone \to \R$ be a $C^6$
function. Assume that there does not exist a continuous function
$f:\bar \Cone\to \R$ such that $\psi \circ \bar h_M=f\circ \bar
h_M-f$ for $M=A$ and $B$. Then the Markov chain $X_n$ on
$\bar\Cone$, starting from $(1,0)$ and whose transition probability
is given by $\nu$, satisfies a nondegenerate central limit theorem
for the function $\psi$, i.e., there exists $\sigma^2>0$ such that,
for any $a\in \R$,
  \begin{equation}
  P \left( \frac{1}{\sqrt{n}} \sum_{k=1}^{n} \psi(X_k) < a
  \right) \to \frac{1}{\sigma \sqrt{2\pi}} \int_{-\infty}^a
  e^{-\frac{t^2}{2\sigma^2}} \dd t.
  \end{equation}
Assume additionally that there do not exist constants $a>0$,
$\lambda>0$ and a continuous function $f:\bar \Cone \to \R/ \lambda
\Z$ such that $\psi\circ \bar h_M= f\circ \bar h_M -f +a \mod
\lambda \Z$ for $M=A$ and $B$. Then $\psi$ satisfies the local limit
theorem:  for any compact subinterval $I$ of $\R$ and any real
sequence $k_n$ such that $k_n/\sqrt{n} \to \kappa\in \R$, then
  \begin{equation}
  \sqrt{n}\; P \left( \sum_{k=1}^{n} \psi(X_k) \in I+k_n
  \right) \to \Leb(I)
  \frac{e^{-\frac{\kappa^2}{2\sigma^2}}}{\sigma \sqrt{2\pi}}.
  \end{equation}
\end{thm}
This result as well as Theorem
\ref{thm:MelangeExponentielFarey} in fact hold for any starting
point of the random walk, there is nothing specific about
$(1,0)$. Note that aperiodicity conditions on $\psi$ are
clearly necessary to get the theorem. For $\kappa=0$, the local
limit theorem can be reformulated as follows. Consider a random
walk on $\bar \Cone \times \R$ whose transition probability is
$Q((x,\omega,z) \to (x',\omega',z'))= P((x,\omega)\to
(x',\omega'))1_{z'=z+\psi(x,\omega)}$. The local limit theorem
simply means that the measure $\sqrt{n}Q^n\delta_{(1,0,0)}$
converges weakly to an explicit multiple of the measure $\mu
\otimes \Leb_{\R/(\log r)\Z} \otimes \Leb_\R$.

Let $T$ be the transformation on the interval $[0,1]$ given by
  \begin{equation}
  \label{def:TFarey}
  T(x)=\frac{x}{1-x}\text{ if }x<1/2,\quad
  T(x)=2-\frac{1}{x}\text{ if }x\geq 1/2.
  \end{equation}
Then $h_A$ and $h_B$ are the inverse branches of the
transformation $T$. The Markov operator corresponding to the
random walk on $[0,1]$ is therefore the adjoint (for the
measure $\mu$) of the composition by $T$, i.e., the transfer
operator associated to $T$. The transformation $T$ is
topologically conjugate to the transformation $x\mapsto 2x$ on
$[0,1]$, and $\mu$ is simply the maximal entropy measure of
$T$, i.e., the pullback of Lebesgue measure under this
conjugacy. Note that $T$ is not uniformly expanding, since it
has neutral fixed points at $0$ and $1$. We can then define a
transformation $\TT$ on $[0,1] \times \R/(\log r)\Z$ whose
inverse branches are $\bar h_A$ and $\bar h_B$, by
  \begin{equation}
  \label{def:TTFarey}
  \TT(x,\omega)=(Tx, \omega+\phi(x)),
  \end{equation}
where $\phi(x)=\log(1-x)$ if $x<1/2$, and $\phi(x)=\log(x)$ if
$x\geq 1/2$. By construction, the Markov operator corresponding to
the random walk on $\bar \Cone$ is the transfer operator associated
to $\TT$ (for the measure $\mu\otimes \Leb$).

With the preceding discussion, we can reformulate the previous
theorems in the general setting of this article: we are going
to study transformations of the form $(x,\omega)\mapsto (Tx,
\omega+\phi(x))$ where $T$ is a nonuniformly expanding
transformation of a manifold $X$, and $\phi$ is a $C^1$
function from $X$ to the circle $\Sbb^1$. Hence, to integrate
the study of Farey sequences in our general setting, it will be
important not to demand uniform expansion, and to be able to
deal with measures which are singular with respect to Lebesgue
measure. These two constraints will justify the forthcoming
definitions, but they will bring along a certain number of
technical difficulties.

\subsection{Definition of nonuniformly partially hyperbolic
skew-products}

\begin{definition}
\label{def:WeakFedererProperty} Let $Z$ be a riemannian
manifold, endowed with a finite measure $\nu$. An open subset
$O$ of $Z$ is said to have the \emph{weak Federer property}
(for the measure $\nu$) if it satisfies the following property.
We work on $O$, with the induced metric, and the geodesic
distance it defines. For any $C>1$, there exist $D=D(O,C)>1$
and $\eta_0=\eta_0(O,C)>0$ such that, for any $\eta<\eta_0$,
there exist disjoint balls $B(x_1,C \eta),\dots, B(x_k,C \eta)$
which are compactly included in $O$, and sets $A_1,\dots, A_k$
contained respectively in $B(x_1,DC\eta),\dots, B(x_k,
DC\eta)$, whose union covers a full measure subset of $O$, and
such that, for any $x'_i \in B(x_i, (C-1)\eta)$, we have
$\nu(B(x'_i,\eta)) \geq \nu(A_i)/D$.

A family of open subsets $(O_n)_{n\in \N}$ is said to uniformly
have the weak Federer property (for the measure $\nu$) if each
set $O_n$ has the weak Federer property and, furthermore, for
any $C>1$, $\sup_{n\in\N} D(O_n,C)<\infty$
\end{definition}
This is a technical covering condition. It is a kind of
weakening of the classical doubling condition, having the
following advantages. On the one hand, it will be satisfied in
many examples (and in particular for Farey sequences, where the
doubling condition does not hold). On the other hand, it is
sufficient to carry out the forthcoming proofs (essentially, it
is the technical condition which is required for Dolgopyat type
arguments to work). The main point of the definition is that
$D$ can be chosen independently of $\eta$: in some sense, the
weak Federer property is a covering lemma with built-in
uniformity.

The following definition describes the class of applications $T$ to
which the results of this article apply. It is large enough to
contain the map \eqref{def:TFarey}, as we will see later on.

\begin{definition}
\label{def:NonUnifDilatant} Let $T$ be a nonsingular
transformation on a riemannian compact manifold $X$ (possibly
with boundary), endowed with a Borel measure $\mu$. Let $Y$ be
a connected open subset of $X$, with finite measure and finite
diameter for the induced metric. We will say that $T$ is a
\emph{nonuniformly expanding transformation of base $Y$, with
exponential tails and the uniform weak Federer property}, if
the following properties are satisfied:
\begin{enumerate}
\item There exist a finite or countable partition (modulo $0$)
$(W_l)_{l\in \Lambda}$ of $Y$, and times $(r_l)_{l\in\Lambda}$
such that, for all $l\in \Lambda$, the restriction of $T^{r_l}$
to $W_l$ is a diffeomorphism between $W_l$ and $Y$, satisfying
$\kappa \norm{v} \leq \norm{DT^{r_l}(x)v} \leq C_l \norm{v}$
for any $x\in W_l$ and $v$ a tangent vector at $x$, for some
constants $\kappa>1$ (independent of $l$) and $C_l$. We will
denote by $T_Y :Y\to Y$ the map which is equal to $T^{r_l}$ on
each set $W_l$.
\item
Let $\HH=\HH_1$ denote the set of inverse branches of $T_Y$ and,
more generally, let $\HH_n$ denote the set of inverse branches of
$T_Y^n$. Let $J(x)$ be the inverse of the jacobian of $T_Y$ at $x$,
with respect to $\mu$. We assume that there exists a constant $C>0$
such that, for any inverse branch $h\in \HH$, $\norm{D((\log J)\circ
h)}\leq C$.
\item
There exists a constant $C$ such that, for any $l$, if
$h_l:Y\to W_l$ denotes the corresponding inverse branch of
$T_Y$, for any $k\leq r_l$, $\norm{ T^k \circ h_l}_{C^1(Y)}
\leq C$.
\item
Let $r:Y\to \N$ be the function which is equal to $r_l$ on
$W_l$. Then there exists $\sigma_0>0$ such that $\int_Y
e^{\sigma_0 r} \dd\mu<\infty$.
\item
Let $\mu_Y$ denote the probability measure induced by $\mu$ on $Y$.
Then the sets $h(Y)$, for $h\in\bigcup_{n\in \N} \HH_n$, uniformly
have the weak Federer property (with respect to $\mu_Y$).
\end{enumerate}
\end{definition}

In this article, we will only consider transformations $T$ of that
type. Hence, we will simply say that $T$ is nonuniformly expanding
with base $Y$.

The first four conditions roughly mean that $T$ is nonuniformly
expanding, and that an induced map $T_Y$ (which is not
necessarily a first return map) is uniformly expanding and
Markov, with exponential tails. This kind of assumptions is
described in \cite{lsyoung:annals, lsyoung:recurrence}, and is
often called a \emph{Young tower structure} in the literature.
The fifth condition is a covering condition. It is probably not
very natural to require it uniformly over the inverse branches
of the iterates of $T_Y$, but it will be satisfied in all the
examples we are going to consider.

Under the first two assumptions, it is a folklore result that
$T_Y$ preserves a probability measure which is equivalent to
$\mu_Y$, whose density is $C^1$ and bounded away from $0$ and
$\infty$. Without loss of generality, we may replace $\mu_Y$ by
this measure (which does not change the assumptions), and we
will therefore always assume that $\mu_Y$ is invariant under
$T_Y$ (and has mass $1$). Inducing from $\mu_Y$ (and using the
fourth assumption), and then renormalizing, we obtain a
probability measure $\tilde\mu$ on $X$ which is invariant under
$T$ and ergodic. However, the restriction of $\tilde\mu$ to $Y$
is in general not proportional to $\mu_Y$, when the return
times $r_l$ are not first return times.

The measure $\tilde\mu$ is always ergodic for $T$, but
sometimes not for its iterates: in general, there exists a
divisor $d$ of $\gcd\{r_l \st l\in \Lambda\}$ and open sets
$(O_i)_{i\in \Z/d\Z}$ such that $T$ maps $O_i$ to $O_{i+1}$,
and the restriction of $T^d$ to each $O_i$ is mixing. For the
sake of simplicity, we will only consider in what follows
transformations $T$ which are mixing, i.e., for which $d=1$.
However, the results we will give have their counterpart in the
general case, since they can be applied to $T^d$ on each set
$O_i$. Note that the mixing of $T$ is equivalent to the
ergodicity of all the iterates $T^n$, and is implied by the
equality $\gcd\{r_l\}=1$.

\begin{rmq}
Under the first four assumptions of Definition
\ref{def:NonUnifDilatant}, and if $T$ is mixing for the
probability measure $\tilde \mu$, then it is exponentially
mixing (for H\"{o}lder continuous functions). This has been proved
by Young in \cite{lsyoung:annals} (in a slightly different
setting) using a spectral gap argument, and again in
\cite{lsyoung:recurrence} using coupling. We will not use these
results of Young. Indeed, our arguments will yield yet another
proof of this exponential mixing, through operator renewal
theory (see in particular Corollary \ref{cor:ExprimeTn0}). This
proof is not new, it is already implicit in \cite{sarig:decay}
and explicit in \cite{gouezel:these}.
\end{rmq}

In a similar setting (the study of expanding semiflows), Ruelle
shows in \cite{ruelle:flot_lent} that a suspension over an
expanding map cannot be exponentially mixing if the roof
function is locally constant. Therefore, it is not surprising
that this case should be excluded from our study, since we will
(among other results) prove exponential mixing.

\begin{definition}
Let $T$ be a nonuniformly expanding transformation of base $Y$,
on a manifold $X$. Let $\phi:X\to \R$ be a $C^1$ function.
Denote by $\phi_Y$ the induced function on $Y$, given by
$\phi_Y(x)=\sum_{i=0}^{r(x)-1}\phi(T^i x)$. We say that $\phi$
is \emph{cohomologous to a locally constant function} if there
exists a $C^1$ function $f:Y\to\R$ such that the function
$\phi_Y-f+f\circ T_Y$ is constant on each set $W_l,\ l\in
\Lambda$.
\end{definition}

If $\phi$ is not cohomologous to a locally constant function, we
define a map $\TT:X\times \Sbb^1 \to X\times \Sbb^1$ by
$\TT(x,\omega)=(Tx, \omega+\phi(x))$. It preserves the probability
measure $\tilde \mu\otimes \Leb$ (in this article, the Lebesgue
measure on the circle $\Sbb^1=\R/2\pi\Z$, denoted by $\Leb$ or
$\de\omega$, will always be normalized of mass $1$). The
transformation $\TT$ is ``nonuniformly partially hyperbolic'', in
the following sense: in each fiber $\Sbb^1$, $\TT$ is an isometry,
while it is expanding in the direction of $X$. Hence, we would like
to talk of partial hyperbolicity. However, since the expansion of
$T$ is not uniform, $T$ can have neutral fixed points or even
critical points. Hence, there may exist points where the
``expansion'' in the $X$ direction does not dominate what is
happening in the fiber. Therefore, the partial hyperbolicity is
rather asymptotic than instantaneous.

\subsection{Limit theorems for nonuniformly partially hyperbolic skew-products}
\label{subsec:ResultatsLimites}

Let $T$ be a nonuniformly expanding map with base $Y$,
preserving the probability measure $\tilde\mu$, and mixing.
Assume that $\mu_Y$ has full support in $Y$. Let $\phi:X\to\R$
be a $C^1$ function which is not cohomologous to a locally
constant function. We consider the skew-product
$\TT(x,\omega)=(Tx,\omega+\phi(x))$.

\begin{thm}
\label{thm:MelangeExponentiel} For any $\alpha>0$, there exist $\bar
\theta<1$ and $C>0$ such that, for all functions $f,g$ from
$X\times\Sbb^1$ to $\C$ respectively bounded and H\"{o}lder continuous
with exponent $\alpha$, and for all $n\in\N$,
  \begin{equation}
  \left| \int f\circ \TT^n \cdot g \dd(\tilde\mu\otimes \Leb)-
  \left(\int f\dd(\tilde\mu\otimes \Leb)\right)
  \left(\int g\dd(\tilde\mu\otimes \Leb)\right)\right| \leq C \bar\theta^n
  \norm{f}_{L^\infty}\norm{g}_{C^\alpha}.
  \end{equation}
\end{thm}

We will then be interested in limit theorems for the transformation
$\TT$. Let $\psi: X\times \Sbb^1 \to \R$ be a H\"{o}lder continuous
function, such that $\int\psi\dd(\tilde\mu\otimes\Leb)=0$. Let
  \begin{equation}
  \label{eq:DefinitSigma2}
  \sigma^2= \int \psi^2 \dd(\tilde\mu\otimes \Leb) + 2\sum_{k=1}^\infty \int
  \psi \cdot \psi \circ \TT^k \dd(\tilde \mu \otimes \Leb).
  \end{equation}
This quantity is well defined, by Theorem
\ref{thm:MelangeExponentiel}.

\begin{prop}
\label{prop:CaracteriseSigma2} We have $\sigma^2\geq 0$. Moreover,
$\sigma^2=0$ if and only if there exists a measurable function
$f:X\times \Sbb^1 \to \R$ such that $\psi = f-f \circ \TT$ almost
everywhere. In this case, the function $f$ has a version which is
continuous on $Y\times\Sbb^1$, and it belongs to $L^p(X\times
\Sbb^1)$ for all $p<\infty$.
\end{prop}

Let us denote by $S_n\psi$ the Birkhoff sums
$\sum_{i=0}^{n-1}\psi\circ \TT^i$. When $\sigma^2$ is nonzero, i.e.,
$\psi$ is not a coboundary, then $\psi$ satisfies the central limit
theorem:
\begin{thm}
\label{thm:LimiteTCL} Let $\psi$ be a H\"{o}lder continuous function on
$X\times \Sbb^1$ with zero average, such that $\sigma^2>0$. Then
$S_n \psi/\sqrt{n}$ satisfies the central limit theorem, i.e.,
$S_n\psi/\sqrt{n}$ converges in distribution (for the probability
measure $\tilde \mu\otimes \Leb$) towards the gaussian distribution
$\boN(0,\sigma^2)$.
\end{thm}

Let us say that $\psi$ is \emph{aperiodic} if there does not exist
$a>0$, $\lambda>0$ and $f:X\times \Sbb^1 \to \R/\lambda\Z$
measurable, such that $\psi= f-f\circ \TT + a \mod \lambda$ almost
everywhere. This implies in particular that $\psi$ is not a
coboundary, hence $\sigma^2>0$.

\begin{prop}
\label{prop:CaracteriseAperiodique} If $\psi$ is a periodic $C^6$
function, there exist $a>0$, $\lambda>0$ and $f:X\times \Sbb^1 \to
\R/\lambda\Z$ measurable such that $\psi= f-f\circ \TT + a \mod
\lambda$ almost everywhere, and $f$ is continuous on $Y\times
\Sbb^1$.
\end{prop}

The notion of periodicity is interesting, since it gives the only
obstruction to the local limit theorem:

\begin{thm}
\label{thm:LimiteTLL} Let $\psi$ be a $C^6$ function on $X\times
\Sbb^1$, with vanishing average, aperiodic (which implies
$\sigma^2>0$). Then the Birkhoff sums $S_n\psi$ satisfy the local
limit theorem, in the following sense: for any compact interval $I$,
any real sequence $k_n$ such that $k_n/\sqrt{n} \to \kappa\in \R$,
we have when $n\to \infty$
  \begin{equation}
  \sqrt{n}\; (\tilde\mu\otimes \Leb) \{(x,\omega)\in X\times \Sbb^1 \st S_n
  \psi(x,\omega) \in I +k_n \} \to \Leb(I)
  \frac{e^{-\frac{\kappa^2}{2\sigma^2}}}{\sigma \sqrt{2\pi}}.
  \end{equation}
\end{thm}

We also obtain numerous other limit theorems (such as the
Berry-Esseen theorem on the speed of $1/\sqrt{n}$ in the central
limit theorem, the renewal theorem, and so on). Instead of giving
precise statements, we will rather give the key estimate which
implies all of them, by showing that the Birkhoff sums $S_n\psi$
essentially behave like a sum of independent identically distributed
random variables:
\begin{thm}
\label{thm:EstimeesValeurPropreFinal} Let $\psi$ be a $C^6$ function
with zero average, such that $\sigma^2>0$. There exist $\eps>0$,
$C>0$, $c>0$ and $\bar\theta<1$ such that, for all functions $f,g$
from $X\times\Sbb^1$ to $\C$ respectively bounded and $C^6$, for any
$n\in\N$, for any $t\in [-\eps,\eps]$,
  \begin{multline}
  \label{eq:GoodEstimeFinal}
  \left| \int e^{itS_n \psi}\cdot f\circ \TT^n \cdot g \dd(\tilde\mu\otimes
  \Leb)- \left( 1-\frac{\sigma^2 t^2}{2}\right)^n
  \left(\int f\dd(\tilde\mu\otimes \Leb)\right)
  \left(\int g\dd(\tilde\mu\otimes \Leb)\right)\right|
  \\
  \leq C (\bar\theta^n + |t| (1-ct^2)^n)
  \norm{f}_{L^\infty}\norm{g}_{C^6}.
  \end{multline}
Moreover, if $\psi$ is aperiodic, for all $t_0>\eps$, there exist
$C>0$ and $\bar\theta<1$ such that, for all $|t|\in [\eps,t_0]$,
  \begin{equation}
  \label{eq:ControleHorsVois0}
  \left| \int e^{itS_n \psi}\cdot f\circ \TT^n\cdot g \dd(\tilde\mu\otimes
  \Leb)\right|
  \leq C \bar\theta^n  \norm{f}_{L^\infty}\norm{g}_{C^6}.
  \end{equation}
\end{thm}
Taking $f=g=1$, we obtain that the characteristic function of
$e^{itS_n\psi}$ essentially behaves like $(1-\sigma^2 t^2/2)^n$,
which makes it possible to prove Theorem \ref{thm:LimiteTCL} for
$C^6$ functions, Theorem \ref{thm:LimiteTLL}, as well as numerous
limit theorems, by mimicking the classical methods in probability
theory for sums of independent identically distributed random
variables. It should just be checked that the additional error term
$\bar\theta^n + |t| (1-ct^2)^n$ does not spoil the arguments. This
has already been done in \cite{gouezel:local}. We will not give
further details on these classical arguments in the following.

Note that, taking $t=0$, Theorem \ref{thm:EstimeesValeurPropreFinal}
implies Theorem \ref{thm:MelangeExponentiel} (for $\alpha=6$, but
this easily implies the general case by a regularization argument).
However, the proof of Theorem \ref{thm:MelangeExponentiel} is
considerably easier than the proof of Theorem
\ref{thm:EstimeesValeurPropreFinal}. Hence, we will give its proof
with full details -- it will also be the occasion to introduce, in a
simple setting, some tools which will be used later on in more
sophisticated versions.

\begin{rmq}
Propositions \ref{prop:CaracteriseSigma2} and
\ref{prop:CaracteriseAperiodique} give automatic regularity for
solutions of the cohomological equation, with a loss of regularity
(arbitrarily small in Proposition \ref{prop:CaracteriseSigma2}, of
$6$ derivatives in Proposition \ref{prop:CaracteriseAperiodique}).
The loss of $6$ derivatives is probably not optimal but, with the
method of proof we use, some loss seems to be unavoidable.

The continuity of $f$ on $Y\times \Sbb^1$ can in general not be
extended to a continuity on the whole space (think for example of a
map $T$ with discontinuities). Nevertheless, using the specificities
of $T$, it is often possible to obtain the continuity of $f$ on
larger sets.
\end{rmq}

\begin{rmq}
Theorem \ref{thm:LimiteTCL} will first be proved for $C^6$ functions
by using Theorem \ref{thm:EstimeesValeurPropreFinal}, and then
extended to H\"{o}lder continuous functions by an approximation
argument. This argument does not apply for the local limit theorem,
which explains our stronger regularity assumption in Theorem
\ref{thm:LimiteTLL}.
\end{rmq}

\begin{rmq}
We require that $\mu_Y$ has full support in $Y$. For some
interesting maps (e.g.\ maps on Cantor sets, see
\cite{naud:Dolgo}), this condition is not satisfied. The full
support condition is used only to get Dolgopyat-like
contraction, in the proof of Lemma \ref{lem:ExisteBranches},
and can be dispensed with, under a stronger condition on
$\phi$. Indeed, if there exist two sequences $h_1,h_2,\dots$
and $h'_1,h'_2,\dots$ of elements of $\HH$, and a point $x$ in
the support of $\mu_Y$, such that the series $\sum_{n=1}^\infty
D(\phi_Y \circ h_n \cdots h_1)(x)$ and $\sum_{n=1}^\infty
D(\phi_Y \circ h'_n \cdots h'_1)(x)$ converge and are not
equal, then the proof of this lemma goes through (note that
this condition is very similar to (NLI) in \cite{naud:Dolgo}).
When $\mu_Y$ has full support, this condition is equivalent to
$\phi$ not being cohomologous to a locally constant function,
as shown in the proof of Lemma \ref{lem:ExisteBranches}.
\end{rmq}

\subsection{Examples}

In the examples, if $T$ and $\phi$ are given, and one wants to apply
the previous results, one should first check that $T$ is
nonuniformly expanding of base $Y$, for some $Y$, and then prove
that $\phi$ is not cohomologous to a locally constant function. The
first issue depends strongly on the map $T$ (see the following list
of examples), but the second one is in general easy to check as
follows, by using periodic orbits.

Assume -- this will be the case in all our examples -- that every
inverse branch $h\in \HH$ of $T_Y$ has a unique fixed point $x_h$.
Let $f$ be a $C^1$ function on $Y$. If $\phi_Y-f+f\circ T_Y$ is
constant on each set $h(Y)$, it has to be equal to $\phi_Y(x_h)$
there. Consequently, the function $g$, equal to $\phi_Y
-\phi_Y(x_h)$ on each set $h(Y)$, is cohomologous to $0$. In
particular, if one can find a periodic orbit of $T_Y$ along which
the Birkhoff sum of $g$ is nonzero, then this is a contradiction,
and $\phi$ can not be cohomologous to a locally constant function.
This can easily be checked in practice: for example, we will use
this argument in the specific case of Farey sequences.

If $1\leq k\leq\infty$, the previous argument moreover shows that,
in the space of $C^k$ functions on $X$, the set of functions $\phi$
which are cohomologous to a locally constant function is contained
in a closed vector subspace of infinite codimension. Hence, the
theorems of Paragraph \ref{subsec:ResultatsLimites} can be applied
for most (in a very strong sense) functions $\phi$.

Let us now describe different classes of maps $T$ which satisfy
Definition \ref{def:NonUnifDilatant}.

\subsubsection*{Nonuniformly expanding maps, and Lebesgue measure}

Let $T$ be a $C^2$ map on a compact riemannian manifold $X$
(possibly with boundary). We assume that $T$ is nonuniformly
expanding, in the following sense (see \cite{alves_bonatti_viana,
alves_luzzatto_pinheiro, gouezel:viana}). Let $S$ be a closed subset
of $X$ with zero Lebesgue measure (corresponding to the
singularities of $T$), possibly empty, and containing the boundary
of $X$. We assume that $T$ is a local diffeomorphism on $X\moins S$,
nondegenerate close to $S$: there exist $B>1$ and $\beta>0$ such
that, for any $x\in X\moins S$ and any nonzero tangent vector $v$ at
$x$,
  \begin{equation}
  \frac{1}{B}\dist(x,S)^\beta \leq \frac{\norm{DT(x)v}}{\norm{v}}
  \leq B \dist(x,S)^{-\beta}.
  \end{equation}
Assume also that, for any $x,y\in X$ with $\dist(x,y)<\dist(x,S)/2$,
  \begin{equation}
  \Bigl|\log \norm{DT(x)^{-1}} -\log \norm{DT(y)^{-1}} \Bigr|
  \leq B \frac{\dist(x,y)}{\dist(x,S)^\beta}
  \end{equation}
and
  \begin{equation}
  \bigl|\log |\det DT(x)^{-1}| - \log |\det DT(y)^{-1}| \bigr|
  \leq B \frac{\dist(x,y)}{\dist(x,S)^\beta}.
  \end{equation}
For $\delta>0$, let $\dist_{\delta}(x,S)=\dist(x,S)$ if
$\dist(x,S)<\delta$, and $\dist_{\delta}(x,S)=1$ otherwise. Let
$\delta:(0,\epsilon_0)\to \R_+$ be a positive function, and let
$\kappa>0$. Assume that, for any $\epsilon<\epsilon_0$, there exist
$C>0$ and $\theta<1$ such that, for any $N\in \N$,
  \begin{equation*}
  \Leb\Bigl\{ x\in X \st
  \exists n \geq N,
  \frac{1}{n}\sum_{k=0}^{n-1} \log
  \norm{DT(T^k x)^{-1}}^{-1} < \kappa \text{ or }
  \frac{1}{n}\sum_{k=0}^{n-1}
  - \log \dist_{\delta(\epsilon)}(T^k x,S) > \epsilon\Bigr\}
  \leq C \theta^N.
  \end{equation*}
This assumption means that the points that do not see the expansion
or are too close to the singularities, after time $N$, have an
exponentially small measure.

As examples of such applications, let us first mention uniformly
expanding maps, of course, but also multimodal maps with infinitely
many branches \cite{araujo_pacifico} (which have thereby infinitely
many critical points), as well as small perturbations of uniformly
expanding maps (such perturbations can have saddle fixed points),
see \cite[section 6]{alves:stability}.

\begin{prop}
Under these assumptions, there exists a subset $Y$ of $X$ such that
$T$ is nonuniformly expanding of base $Y$, for Lebesgue measure.
\end{prop}
\begin{proof}
This theorem is essentially proved in \cite[Theorem
4.1]{gouezel:viana}. More precisely, this theorem constructs a
subset $Y$ of $X$ and a partition of $Y$ such that the first four
properties of Definition \ref{def:NonUnifDilatant} are satisfied.
The set $Y$ is an open set with piecewise $C^1$ boundary, and each
inverse branch $h$ can be extended to a neighborhood of $Y$.

If the boundary of $Y$ were $C^1$ (and not merely piecewise $C^1$),
each set $h(Y)$ would also be an open set with $C^1$ boundary, and
the uniform weak Federer property would directly result from the
good doubling properties of Lebesgue measure. However, if the
boundary of $Y$ is only piecewise $C^1$, the images of the boundary
components by an inverse branch $h$ could make smaller and smaller
angles, which could prevent the uniform weak Federer property from
holding.

Therefore, we have to modify slightly the construction in
\cite{gouezel:viana} to obtain a set $Y$ with $C^1$ boundary. In
that article, one starts from a partition $U_i$ of $X$ (into sets
with piecewise $C^1$ boundary), and one subdivides each set $U_i$
into subsets $V_j$ which are sent by some iterate of $T$ on one of
the sets $U_k$. The set $Y$ is then one of the $U_i$'s, and the
desired partition of $Y$ is obtained by inducing from the $V_j$'s
(see \cite[section 4]{gouezel:viana} for details).

To obtain a smooth $Y$, we also start from a partition $U_i$, but we
decompose $U_i$ as $U_i^1 \cup U_i^2$ where $U_i^1$ is a ball inside
$U_i$ and $U_i^2$ is its complement. Applying the construction of
\cite{gouezel:viana} separately to each set $U_i^1$ and $U_i^2$, we
subdivide them into sets $V_j$ which are sent by some iterate of $T$
to some $U_k$. We finish the construction by taking for $Y$ one of
the sets $U_i^1$, and inducing on it.
\end{proof}

To apply the results of Paragraph \ref{subsec:ResultatsLimites}, one
needs an additional mixing assumption, which is satisfied as soon as
all the iterates of $T$ are topologically transitive on the
attractor $\bigcap_{n\geq 0} T^n(X)$ (see \cite{gouezel:viana}).

\subsubsection*{Multimodal maps of Collet--Eckmann type}

Let $T$ be a multimodal map on a compact interval $I$. If the
derivative of $T^n$ along the postcritical orbits grow exponentially
fast, and $T$ is not renormalizable (which prevents periodicity
problems), \cite{bruin_luzzatto_strien} shows that there exists a
unique absolutely continuous invariant probability measure
$\tilde\mu$, and that $T$ is exponentially mixing for this measure.

To prove this result, the authors show that there exist an interval
$Y$ and a subpartition $W_l$ of $Y$ satisfying the first four
properties of Definition \ref{def:NonUnifDilatant}, for Lebesgue
measure. Since the sets $h(Y)$ (for $h\in \bigcup_{n\in \N} \HH_n$)
are all intervals, the uniform weak Federer property is also
trivially satisfied by Lebesgue measure.

\subsubsection*{Gibbs measures in dimension $1$}

If $T$ is a $C^2$ uniformly expanding map on a compact connected
manifold $X$, and $u:X\to \R$ is a $C^1$ function, there exists a
unique invariant probability measure $\mu$ which maximizes the
quantity $h_\nu(T)+\int u \dd\nu$ over all invariant probability
measures $\nu$. This is the so-called \emph{Gibbs measure}
associated to the potential $u$.

In general, it is unlikely that such a Gibbs measure satisfies
the weak Federer property (unless $\mu$ is equivalent to
Lebesgue measure, which corresponds to potentials $u$ which are
cohomologous to $-\log \det(DT)$). Indeed, the proof of the
weak Federer property in the previous examples relies in an
essential way on the good doubling properties of Lebesgue
measure.

However, in dimension $1$ (i.e., if $T$ is a circle map), the
iterates of $T$ are conformal, which implies that $\mu$
satisfies the weak Federer property, and our results apply.
Proofs of the Federer property in this setting have been given
by Dolgopyat or Pollicott, but with small imprecisions, so we
will give a full proof in Proposition
\ref{prop:DemontreFedererGibbs} (as a very simple consequence
of the methods we develop to treat the Farey sequence). Note
that the same results also apply in higher dimension, for
conformal uniformly expanding maps (since uniformly expanding
maps always admit Markov partitions).

\subsubsection*{Farey sequences}

The results of Paragraph \ref{subsec:ResultatsLimites} also apply to
the map \eqref{def:TTFarey}, which generates the Farey sequence.
However, the proof requires more work, since checking the weak
Federer property is not trivial. Moreover, the most interesting
results stated in Theorem \ref{thm:ThmLimiteFarey} are pointwise
results (for a random walk starting from $(1,0)$), while the
statements of Paragraph \ref{subsec:ResultatsLimites} are on average
results. To prove the pointwise statements, we will therefore need
to use more technical results, established during the course of the
proof of Theorems \ref{thm:MelangeExponentiel} and
\ref{thm:EstimeesValeurPropreFinal}. As a consequence, the results
of Paragraph \ref{subsec:Farey} will be proved at the end of the
article, in Section \ref{sec:PreuvesFarey}.

\subsection{Method of proof, and contents of the article}

In general, to prove exponential mixing and a local limit theorem,
it is very comfortable to have a \emph{spectral gap property} for a
transfer operator (the spectral perturbation methods then yield the
desired results quite automatically). The spectral gap is in general
a consequence of some expansion or contraction properties. However,
in our setting, the map $\TT$ is an isometry in the fibers, and a
spectral gap seems therefore difficult to obtain. Note that
\cite{tsujii:semiflow} manages to construct a space with a spectral
gap for such maps, but under strong assumptions: the map $T$ should
be uniformly expanding, and $\tilde\mu$ should be absolutely
continuous with respect to Lebesgue measure. These properties are
unfortunately not satisfied in our setting, and we will thus have to
work without a spectral gap (on the space $X\times \Sbb^1$).

Dolgopyat developed in \cite{dolgopyat:decay,
dolgopyat:compact_extensions} techniques which he used to prove the
exponential decay of correlations for maps $\TT$ as above, if $T$ is
uniformly expanding. His main idea is to work in Fourier
coordinates, to see that each frequency is left invariant by the
transfer operator associated to $\TT$, and to obtain explicit bounds
on the mixing speed in each frequency (by using oscillatory
integrals, which give explicit compensations). The gain is not
uniform with respect to the frequency (which accounts for the lack
of spectral gap), but the estimates are nevertheless sufficiently
good to obtain exponential mixing.

We will use in an essential way Dolgopyat's ideas in this article,
as a technical tool. This tool applies to uniformly expanding maps,
which is not the case of our map $T$, we will therefore need to
induce on the set $Y$ to get uniform expansion. To obtain
information on the initial map, we will then make use of
(elementary) ideas of generating series and renewal theory.

The real difficulty of the article lies in the local limit
theorem, since a spectral gap property seems more or less
necessary to any known proof of the local limit theorem, while
Dolgopyat's arguments do not give such a spectral gap. If we
try to work on the level of frequencies, as for the exponential
mixing, we quickly run into the following additional
difficulty: if $f$ is a function of frequency $k$, i.e.,
$f(x,\omega)=u(x)e^{ik\omega}$, then $e^{it\psi} f$ is not any
more a function of frequency $k$. In other words, the
multiplication by $e^{it\psi}$ -- which is at the heart of the
proof of the local limit theorem for the function $\psi$ --
mixes the different frequencies together. Hence, even though
Dolgopyat's techniques give a good control at high frequencies,
this control is instantaneously ruined by the multiplication by
$e^{it\psi}$, which can go back into low frequencies, where no
control is available.

The central idea for the proof of the local limit theorem is to
induce at the same time in $x$ and in $k$: we consider some
kind of random walk on the space $X\times \Z$ (where the $\Z$
factor corresponds to the space of frequencies), and we induce
on a subset $Y\times[-K,K]$ where $K$ is large enough so that
what happens outside of this set can be controlled by
Dolgopyat's tools. The main interest of this process is that
the induced operator on $Y\times[-K,K]$ has a spectral gap, and
can be studied very precisely. Using techniques of operators
renewal theory \cite{sarig:decay, gouezel:local}, we will then
use this information to obtain a global control on $X\times
\Z$, finally yielding Theorem
\ref{thm:EstimeesValeurPropreFinal}.

\begin{rmq}
The next natural question is to study maps of the form
$\TT':(x,\omega,\omega')\mapsto (Tx, \omega+\phi(x),
\omega'+\psi(x,\omega))$, where $T$ and $\phi$ are as above. If
$\psi$ is aperiodic, Theorem
\ref{thm:EstimeesValeurPropreFinal} shows that the correlations
of functions of the form $u(x,\omega)e^{ik\omega'}$ (where $u$
is $C^6$ and $k\in \Z$) tend to $0$. Since the linear
combinations of such functions are dense in $L^2$, this implies
that $\TT'$ is mixing. It is even Bernoulli, by the following
argument: first, $T$ (or rather its natural extension) is
Bernoulli since it is mixing and non-uniformly hyperbolic (see
e.g.\ \cite{OrnsteinWeiss}). Since $\TT$ is a mixing isometric
extension of $T$, it is also Bernoulli by
\cite{Rudolph:IsometricExtension}. The same argument applied to
$\TT$ then implies that $\TT'$ is Bernoulli.

However, to prove further results on $\TT'$, such as
exponential mixing or the local limit theorem (probably under
stronger assumptions on $\psi$) seems out of reach by currents
techniques. More precisely, we use Dolgopyat's techniques
(which give precise explicit estimates for the map $T$) to
study the map $\TT$ (and obtain, by an abstract compactness
argument, non-explicit estimates for $\TT$). To go one step
further and study precisely $\TT'$, we would need
\emph{explicit} estimates for $\TT$ (i.e., in
\eqref{eq:ControleHorsVois0}, we would need to control
$\bar\theta$ and $C$ in terms of $t_0$), which seems
considerably more difficult.
\end{rmq}

\medskip

The article is organized as follows: in Section \ref{sec:Outils}, we
state a theorem on transfer operators giving all the technical
estimates we shall need further on (with contraction in the
classical sense, or in Dolgopyat norms). This technical theorem will
be proved in an appendix. In Section \ref{sec:MelangeExp}, it is
used to prove Theorem \ref{thm:MelangeExponentiel}. The proof is a
baby version of the proof of the local limit theorem, introducing
some tools on renewal operators that will be used further on. In
Section \ref{sec:Strategie}, we describe in details the strategy of
the proof of the local limit theorem, and give two technical results
which are essential in its proof. The proof itself is given in
Section \ref{sec:DemontreLocal}. Finally, Section
\ref{sec:PreuvesFarey} is devoted to the proof of the results on
Farey sequences, as stated in Paragraph \ref{subsec:Farey}.

In all the following, we fix once and for all a map $T$ which is
nonuniformly expanding of base $Y$, mixing, together with a function
$\phi$ which is not cohomologous to a locally constant function.

\section{Tools on transfer operators}

\label{sec:Outils}

For $k\in \Z$ and $v\in C^1(Y)$, we set
  \begin{equation}
  \Lp_k v(x)=\sum_{h\in \HH} e^{-ik \phi_Y(hx)} J(hx)v(hx),
  \end{equation}
and we define $\Lp=\Lp_0$. This is the transfer operator
associated to $T_Y$. For $x\in Y$ and $n\in \N$, let us also
write $S_n^Y \phi_Y(x)=\sum_{i=0}^{n-1}\phi_Y(T_Y^ix)$.

For $n\in \N$ and $x\in Y$, let
$r^{(n)}(x)=\sum_{i=0}^{n-1}r(T_Y^ix)$. For $n\in \N$, $A>0$
and $\epsilon>0$, we will denote by $\CC^{A,\epsilon}_{n}$ the
set of functions $v$ from $Y$ to $\C$ which are $C^1$ on each
set $h(Y)$ for $h\in \HH_n$, and such that the quantity
  \begin{equation}
  \norm{v}_{\CC^{A,\epsilon}_{n}}= \sup_{h\in \HH_n} \sup_{x\in
  Y} \max ( |v(hx)|, \norm{D(v\circ h)(x)}/A) / e^{\epsilon
  r^{(n)}(hx)}
  \end{equation}
is finite. These are the functions we will be working with.
They can be unbounded, but their explosion speed is controlled
by the return time. Typically, if one starts from a smooth
function on $X$ and induces, the resulting function will be
unbounded but in $\CC^{A,\epsilon}_1$ for some $A,\epsilon$. In
particular, for any $A>0$ and $\epsilon>0$, we have $\sup_{n\in
\N} \norm{S_n^Y\phi_Y}_{\CC^{A,\epsilon}_n}<\infty$. Note that
the set of functions $\CC^{A,\epsilon}_n$ does not depend on
$A$, but the corresponding norm does.

Let $k\in \Z$ and $C_0>1$. We will denote by $\E_k(C_0)$ the
set of pairs $(u,v)$ of functions from $Y$ to $\C$ such that
$|v|\leq u$ and $\max(\norm{Dv},\norm{Du}) \leq C_0
\max(1,|k|)u$. This set is a cone, i.e., it is stable under
addition and multiplication by nonnegative real numbers. We
will also write $\norm{v}_{D_k(C_0)}$ (or simply
$\norm{v}_{D_k}$) for the infimum of the quantities
$\norm{u}_{L^4}$ over all functions $u$ such that $(u,v)\in
\E_k(C_0)$. Since $\E_k(C_0)$ is a cone, this is a norm,
satisfying $\norm{v}_{L^4}\leq \norm{v}_{D_k}\leq
\norm{v}_{C^1}$. The $D_k$ norm has been (implicitly) used by
Dolgopyat, and is very useful since it enjoys good contraction
properties for the action of the transfer operator $\Lp_k$.

We will freely use the following trivial inequalities: if
$|k|\leq |\ell|$, then $\norm{v}_{D_\ell}\leq \norm{v}_{D_k}$.
Moreover, for any $k$, $\norm{v}_{D_k}\leq \norm{v}_{C^1}$.
Finally, we have $\norm{v}_{\CC^{A,\epsilon'}_{n}} \leq
\norm{v}_{\CC^{A,\epsilon}_{n}}$ as soon as $\epsilon'\geq
\epsilon$.

The theorem we will use is the following. Recall that $T$ is a fixed
nonuniformly expanding transformation of base $Y$, and that $\phi$
is a $C^1$ function which is not cohomologous to a locally constant
function, also fixed once and for all.

\begin{thm}
\label{thm:MainContraction} There exist $N>0$, $C_0>1$,
$\epsilon>0$ and $\theta\in (2^{-1/(1010N)},1)$, such that, for
any $M\geq 1$, the following properties hold.
\begin{description}
\item[Classical contraction]
for any $A\geq 1$, there exists a constant $C(A)$ such that, for any
$\psi \in \CC^{A,4\epsilon}_{MN}$ and for any $v\in C^1(Y)$,
  \begin{equation}
  \label{eq:ContracteTriviale}
  \norm{\Lp^{MN}(\psi v)}_{C^1} \leq \theta^{100 MN} \left(
  \sup_{x\in Y} |\psi(x)|/ e^{4\epsilon r^{(MN)}(x)} \right)
  \norm{v}_{C^1} + C(A) \norm{\psi}_{\CC^{A,4\epsilon}_{MN}}
  \norm{v}_{C^0}.
  \end{equation}
Moreover, there exists $C>0$ satisfying: let $A\geq 1$, let
$\psi_1,\dots,\psi_n \in \CC^{A,4\epsilon}_{MN}$ and let $v\in
C^1(Y)$. Write $v^0=v$ and $v^i=\Lp^{MN}(\psi_i v^{i-1})$. Then
  \begin{equation}
  \label{eq:ContractC1Iteres}
  \norm{v^n}_{C^1} \leq C A \left(\prod_{i=1}^n \norm{\psi_i}_{\CC^{A,4\epsilon}_{MN}}\right)
  \left(\theta^{100 MNn}\norm{v}_{C^1}+\theta^{-MNn}\norm{v}_{L^2}\right).
  \end{equation}
\item[Dolgopyat's contraction]
for any $A\geq 1$, there exists $K=K(A,M)$ such that, for any
$|k|\geq K$, for any $C^1$ function $v:Y\to\C$, for any
function $\psi \in \CC^{A,4\epsilon}_{MN}$,
  \begin{equation}
  \label{eq:GagneL4}
  \norm{\Lp_k^{MN}(\psi v)}_{D_k} \leq \theta^{100 MN}
  \norm{\psi}_{\CC^{A,4\epsilon}_{MN}}  \norm{v}_{D_{2^M k}}.
  \end{equation}
Moreover, for any $|\ell|\geq |k| \geq K$, we also have
  \begin{equation}
  \label{eq:PerdPasTrop}
  \norm{\Lp_k^{MN}(\psi v)}_{D_\ell} \leq \theta^{-MN}
  \norm{\psi}_{\CC^{A,4\epsilon}_{MN}} \norm{v}_{D_{2^M \ell}}.
  \end{equation}
\end{description}
\end{thm}

The first half of the theorem is really classical (it is a
consequence of the usual contraction of transfer operators on
spaces of Lipschitz or $C^1$ functions), the second half is
less classical but should not be surprising to a reader who is
used to Dolgopyat's techniques. However, this result contains
additional technical difficulties with respect to the same kind
of results in the literature. Indeed the functions in
$\CC^{A,\epsilon}_{MN}$ are usually unbounded and have
unbounded derivatives. Moreover, the application of Dolgopyat's
arguments is problematic since the function $\phi_Y$ is also
unbounded with unbounded derivative. As a consequence, the
proof of this theorem is quite unpleasant, even though it does
not need additional conceptual ideas, only technical ones.
Therefore, the proof of Theorem \ref{thm:MainContraction} is
postponed to Appendix \ref{app:app}.

\emph{In all the rest of the article (but Appendix
\ref{app:app}), $N$, $C_0$, $\epsilon$ and $\theta$ will be
fixed once and for all, and will denote the constants given by
Theorem \ref{thm:MainContraction}.}

\begin{rmq}
Note that the bounds with $\norm{\psi}_{\CC^{A,4\epsilon}_{MN}}$
imply the same bounds with $\norm{\psi}_{\CC^{A,\epsilon}_{MN}}$.
Most of the time, we will only need this weaker version (the
inequalities with $4\epsilon$ simply give a small additional margin,
which will be useful from time to time).
\end{rmq}

\begin{rmq}
Concerning the precise formulation of Theorem
\ref{thm:MainContraction}, let us make two additional remarks which
are apparently technical but are in fact extremely important for the
forthcoming proofs.
\begin{enumerate}
\item The theorem for $M=1$ is sufficient to obtain the exponential
mixing (and to prove the theorem for $M=1$ we only need the weak
Federer property of $Y$, and no uniformity on the inverse branches).
However, to prove the local limit theorem, we will need to take
larger and larger $M$'s: since $\theta$ is independent of $M$, the
gain $\theta^{100MN}$ will enable us to control some terms which are
polynomially growing with $M$. The uniformity in $M$ in Theorem
\ref{thm:MainContraction} is therefore crucial.
\item Since $\norm{v}_{D_{2^M k}}\leq \norm{v}_{D_k}$, the
inequality \eqref{eq:GagneL4} is stronger than
  \begin{equation}
  \label{eq:kljqmsklfdjqlmsdf}
  \norm{\Lp_k^{MN}(\psi v)}\leq \theta^{100MN}
  \norm{\psi}_{\CC^{A,4\epsilon}_{MN}} \norm{v}_{D_k}.
  \end{equation}
The inequality \eqref{eq:kljqmsklfdjqlmsdf} would be sufficient to
prove the exponential mixing. However, to prove the local limit
theorem, we will jump from one frequency to another, and the
additional gain in the index given by \eqref{eq:GagneL4} will be
crucial (especially  in the proof of Lemma \ref{Lem:ContractionL2}).
\end{enumerate}
\end{rmq}

The following general lemma will also be required:

\begin{lem}
\label{lem:DecritFctPropre} Let $T_0$ be an ergodic
transformation of a probability space, with corresponding
transfer operator $\hat T_0$. Let $g$ be a nonzero integrable
function, let $f$ be a measurable function with modulus at most
$1$, and let $\lambda\in \C$ with $|\lambda|\geq 1$. We assume
that $\lambda g= \hat T_0(fg)$. Then $|\lambda|=1$, $|f|=1$
almost everywhere, and $\lambda g\circ T= fg$ almost
everywhere.
\end{lem}
\begin{proof}
We have $|\lambda||g| \leq  \hat T_0|g|$. Integrating this
equation yields $|\lambda|\norm{g}_{L^1} \leq \norm{g}_{L^1}$,
which implies $|\lambda|=1$. Moreover, the function $\hat
T_0|g|-|g|$ is nonnegative and has zero integral, hence it
vanishes almost everywhere. Since $\hat T_0 |g|=|g|$, the
measure with density $|g|$ is invariant. By ergodicity, $|g|$
is almost everywhere constant (and this constant is nonzero).
The equation $\lambda g=\hat T_0(fg)$ becomes $\hat T_0
(\lambda^{-1} f g/g\circ T_0) =1$. Therefore,
  \begin{equation}
  1=\int \lambda^{-1} f \frac{g}{g\circ
  T_0} \leq \int \left| \lambda^{-1} f \frac{g}{g\circ
  T_0}\right| \leq 1.
  \end{equation}
This shows that the function $\lambda^{-1} f \frac{g}{g\circ T}$ has
to be equal to $1$ almost everywhere.
\end{proof}

\section{Exponential mixing}
\label{sec:MelangeExp}

\subsection{A model for $\TT$}
\label{subsec:modele} For $n\in \N$, we are going to define an
artificial transformation, which will model the dynamics of
$\TT$, as follows. Let $X^{(n)}=\{ (x,i) \st x\in Y,
i<r^{(n)}(x)\}$, we define a map $U^{(n)}$ (or simply $U$ if
$n$ is implicit) on $X^{(n)}$ by $U(x,i)=(x,i+1)$ if
$i+1<r^{(n)}(x)$, and $U(x,r^{(n)}(x)-1)=(T_Y^n(x),0)$. Let
$\pi^{(n)}: X^{(n)} \to X$ be given by $\pi^{(n)}(x,i)=T^i(x)$,
we obtain $\pi^{(n)}\circ U= T\circ \pi^{(n)}$. We endow each
set $h(Y)\times \{i\}$, for $h\in \HH_n$ and $i<r^{(n)}\circ
h$, with the restriction of the measure $\mu_Y$ to $h(Y)$. This
yields a measure $\mu^{(n)}$ which is invariant under $U$ and
whose restriction to $Y\times \{0\}$ is equal to $\mu_Y$.
Strictly speaking, the map $U$ is not defined everywhere since
some points of $Y$ do not come back to $Y$. However, it is
defined $\mu^{(n)}$ almost everywhere, which will be sufficient
for our needs. The measure $\pi^{(n)}_* \mu^{(n)}$ is
absolutely continuous with respect to $\tilde \mu$ and
invariant, hence these measures are proportional by ergodicity.
In particular, setting $\tilde \mu^{(n)}= \mu^{(n)}/
\mu^{(n)}(X^{(n)})$, we have $\pi^{(n)}_* \tilde\mu^{(n)}
=\tilde \mu$.

We also endow $X^{(n)}$ with a metric, as follows. The set $Y$
is canonically embedded in $X^{(n)}$ by $y \mapsto (y,0)$, we
endow the image of this embedding by the metric of $Y$. Let
$h\in \HH_n$ and $0<i< r^{(n)}\circ h$ (this function is
constant on $Y$). The map $U^{r^{(n)}\circ h -i}$ is a
bijection between $h(Y)\times\{i\}$  and $Y\times\{0\}$, we
choose the metric on $h(Y)\times\{i\}$ so that this map is an
isometry.

With this choice of the metric, the map $U$ is very expanding
on the points of the form $(y,0)$ (it expands the metric by at
least $\kappa^n$), and it is a local isometry on the points
$(y,i)$ with $i>0$. Since $T$ satisfies the third property of
Definition \ref{def:NonUnifDilatant}, the map $\pi^{(n)}$ is
almost a contraction: there exists a constant $C$ such that
  \begin{equation}
  \label{eq:piContractePresque}
  \norm{ D\pi^{(n)}(x) \cdot v} \leq C \norm{v}
  \end{equation}
for any $x\in X^{(n)}$ and $v$ tangent at $x$. If $u:X \to \C$
is a $C^1$ function, the function $u\circ \pi^{(n)}$ is then
also $C^1$ on $X^{(n)}$, and $\norm{ u\circ
\pi^{(n)}}_{C^1}\leq C \norm{u}_{C^1}$.

We finally define a map  $\UU=\UU^{(n)}$ on $X^{(n)}\times
\Sbb^1$, by $\UU(x, \omega)=(Ux, \omega+\phi\circ
\pi^{(n)}(x))$. If we define $\tilde \pi^{(n)}: X^{(n)}\times
\Sbb^1 \to X\times \Sbb^1$ as $\pi^{(n)}\times \Ide$, then
$\UU$ is a model for $\TT$ since $\tilde \pi^{(n)} \circ
\UU=\TT \circ \tilde \pi^{(n)}$. To study the properties of
$\TT$, it will therefore be sufficient to understand
$\UU^{(n)}$ (for any conveniently chosen $n$). Abusing
notations, we will simply write $\phi$ on $X^{(n)}$ instead of
$\phi\circ \pi^{(n)}$. We will also identity $Y$ with
$Y\times\{0\}\subset X^{(n)}$.

The map $U$ is not always mixing for the measure $\tilde\mu^{(n)}$:
setting
  \begin{equation}
  d=d^{(n)}=\gcd\{ r^{(n)}(x) \st x\in Y\},
  \end{equation}
then $U$ is mixing if and only if $d=1$. If $d>1$, let us
write, for $k\in \Z/d\Z$, $\tilde\mu^{(n)}_k$ for the
probability measure induced by $\tilde \mu^{(n)}$ on the set
$\{(x,i) \st i = k \mod d\}$. Then each measure $\tilde
\mu^{(n)}_k$ is invariant under $U^d$, and mixing. The measure
$\pi^{(n)}_* \tilde\mu^{(n)}_k$ is absolutely continuous with
respect to $\tilde\mu$ and invariant under $T^d$. Since $T^d$
is ergodic (because $T$ is mixing), this yields $\pi^{(n)}_*
\tilde\mu^{(n)}_k=\tilde\mu$.

\subsection{The transfer operator associated to $\UU^{(N)}$}
\label{par:DonneOperateursRetour}

In the rest of this section, we work on $X^{(N)}$, where $N$ is
given by Theorem \ref{thm:MainContraction} (and fixed once and
for all). This theorem will make it possible to study the
transfer operator $\hat\UU$ associated to the map
$\UU=\UU^{(N)}$. Our goal in this section is to use this
information to prove Theorem \ref{thm:MelangeExponentiel}.

\emph{To keep the arguments as transparent as possible, we will
assume until the end of the proof, and without repeating it
each time, that $d^{(N)}=\gcd\{r^{(N)}(x)\}$ is equal to $1$.
At the end of the proof, we will indicate the modifications to
be done in the general case.}

Let us write a function $v$ on $X^{(N)}\times \Sbb^1$ as
$v(x,\omega)=\sum_{k\in \Z} v_k(x)e^{ik\omega}$, i.e.,
  \begin{equation}
  v_k(x)=\int v(x,\omega) e^{-ik\omega}\dd \omega,
  \end{equation}
where $\dd\omega$ denotes the normalized Lebesgue measure on
$\Sbb^1$. If $\hat\UU$ is the transfer operator associated to $\UU$,
and $\KK$ is the inverse of the jacobian of $U$ for $\mu^{(N)}$,
  \begin{align*}
  \hat\UU v(x,\omega)&=\sum_{\UU(x',\omega')=(x,\omega)} \KK(x') v(x',\omega')
  = \sum_{U(x')=x} \KK(x') v(x', \omega-\phi(x'))
  \\&
  =\sum_{k\in \Z} \sum_{U x'=x} \KK(x') v_k(x') e^{ik
  (\omega-\phi(x'))}.
  \end{align*}
In the same way, if $\KK^{(n)}$ denotes the jacobian of $U^n$,
  \begin{equation}
  \label{EtudieUUn}
  \hat\UU^n v(x,\omega)= \sum_{k\in \Z} \sum_{U^n x'=x}
  \KK^{(n)}(x')v_k(x') e^{ik (\omega-S_n\phi(x'))}.
  \end{equation}
Hence, the operator $\hat\UU^n$ acts diagonally on each frequency,
by an operator
  \begin{equation}
  \Mp^n_k v(x)=\sum_{U^n x'=x}
  \KK^{(n)}(x')v(x') e^{-ik S_n\phi(x')}.
  \end{equation}
We will understand separately the action of $\Mp_k$ for each
$k$. Using the induction process, we will be able to understand
this operator for points $x$, $x'$ belonging to the base $Y$ of
$X^{(N)}$. We will then use this information to reconstruct the
whole operator $\Mp_k$. To do so, let us define the following
operators:
  \begin{gather}
  R_{n,k}v(x)=\sum_{\substack{U^n x'=x
  \\ x'\in Y,U x',\dots, U^{n-1}x'\not\in Y,U^nx'\in Y}}
  \KK^{(n)}(x')v(x') e^{-ik S_n\phi(x')},
  \\%
  T_{n,k}v(x)=\sum_{\substack{U^n x'=x
  \\ x'\in Y,U^nx'\in Y}}
  \KK^{(n)}(x')v(x') e^{-ik S_n\phi(x')},
  \\%
  A_{n,k}v(x)=\sum_{\substack{U^n x'=x
  \\ x'\in Y, Ux',\dots, U^n x'\not\in Y}}
  \KK^{(n)}(x')v(x') e^{-ik S_n\phi(x')},
  \\%
  \label{definitBnk}
  B_{n,k}v(x)=\sum_{\substack{U^n x'=x
  \\ x',\dots, U^{n-1} x'\not\in Y, U^n x'\in Y}}
  \KK^{(n)}(x')v(x') e^{-ik S_n\phi(x')},
  \\%
  C_{n,k}v(x)=\sum_{\substack{U^n x'=x
  \\ x',\dots, U^{n} x'\not\in Y}}
  \KK^{(n)}(x')v(x') e^{-ik S_n\phi(x')}.
  \end{gather}
The main interest of these definitions is the following. First,
cutting an orbit according to the first and last time it belongs to
$Y$, we get
  \begin{equation}
  \label{eq:DonneMpnk}
  \Mp^n_k=
  C_{n,k}+\sum_{a+i+b=n} A_{a,k}T_{i,k}B_{b,k}.
  \end{equation}
Moreover, considering all the times an orbit belongs to $Y$, we
obtain
  \begin{equation}
  \label{ExprimeTnk}
  T_{n,k}=\sum_{p=1}^\infty \sum_{j_1+\dots+j_p=n}
  R_{j_1,k}\dots R_{j_p,k}.
  \end{equation}
Finally, for $z\in \C$ with modulus at most $e^{\epsilon}$, we have
  \begin{equation}
  \label{exprimeRnk}
  \sum_{n>0} z^n R_{n,k} v= \Lp^N_k( z^{r^{(N)}}v).
  \end{equation}
The restriction $|z|< e^{\epsilon}$ ensures that this operator is
well defined, by Theorem \ref{thm:MainContraction}. More precisely,
we even have:
\begin{lem}
\label{lem:BorneRn} There exists $C>0$ such that, for any $n\in \N$,
for any $k\in \Z$,
  \begin{equation}
  \norm{R_{n,k}v}_{C^1(Y)} \leq C \max(1,|k|)e^{-2n\epsilon}
  \norm{v}_{C^1(Y)}.
  \end{equation}
\end{lem}
\begin{proof}
Let $\psi_{n,k}(x)=e^{-ik S_N^Y \phi_Y(x)}$ if $r^{(N)}(x)=n$,
and $0$ otherwise, so that $R_{n,k} v=\Lp^N( \psi_{n,k} v)$. We
will show that $\norm{\psi_{n,k}}_{\CC^{1,4\epsilon}_N} \leq
C\max(1,|k|) e^{-2\epsilon n}$, which will conclude the proof
by \eqref{eq:ContractC1Iteres}.

We have $|\psi_{n,k}(x)|\leq e^{-2n\epsilon} e^{2\epsilon
r^{(N)}(x)}$. Moreover, if $h\in \HH_N$ satisfies $r^{(N)}\circ
h=n$, we have
  \begin{equation}
  \norm{D(\psi_{n,k}\circ h)(x)}\leq C |k| r^{(N)}(hx) \leq
  C|k|e^{2\epsilon r^{(N)}(hx)} \leq C |k| e^{-2\epsilon n}
  e^{4 \epsilon r^{(N)}(hx)}.
  \end{equation}
This proves the lemma.
\end{proof}

\subsection{Study of the operators $T_{n,k}$}

In Equation \eqref{eq:DonneMpnk}, the complicated part in the
expression of $\Mp_k^n$ comes from $T_{i,k}$, since the other
operators are more or less explicit. This paragraph is devoted to
the study of the operators $T_{i,k}$, by using \eqref{ExprimeTnk}.

\begin{lem}
\label{lem:TnkPourknot0} There exist $C>0$ and $\bar\theta<1$
such that, for any $k\in \Z\moins\{0\}$, for any $n\in \N$ and
for any $v\in C^1(Y)$, $\norm{T_{n,k} v}_{C^1} \leq C k^2
\bar\theta^n \norm{v}_{C^1}$.
\end{lem}
\begin{proof}
For $k\in\Z$ and $|z|\leq e^{\epsilon}$, let us write
$\Lp_{k,z} v= \Lp_k^N( z^{r^{(N)}} v) =\Lp^N(e^{-ikS_N^Y\phi_Y}
z^{r^{(N)}} v)$. Since $\Lp_{k,z}=\sum z^j R_{j,k}$ by
\eqref{exprimeRnk}, Lemma \ref{lem:BorneRn} shows that this
operator acts continuously on $C^1(Y)$, and that $z\mapsto
\Lp_{k,z}$ is holomorphic on the disk $\{|z|\leq
e^{\epsilon}\}$. Formally, we can rewrite \eqref{ExprimeTnk} as
$\sum T_{n,k}z^n = (I-\sum
R_{j,k}z^j)^{-1}=(I-\Lp_{k,z})^{-1}$. Hence, for any path
$\gamma$ in $\C$ around $0$ bounding a domain on which
$\Id-\Lp_{k,z}$ is invertible for any $z$, we have for any
$n\in \N$
  \begin{equation}
  \label{ExprimeTnkCauchy}
  T_{n,k}=\frac{1}{2i\pi}\int_\gamma z^{-n-1}
  (I-\Lp_{k,z})^{-1}\dd z.
  \end{equation}
We are going to use this equation as well as the information on
$\Lp_{k,z}$ to estimate $T_{n,k}$.

\emph{First step.} Fix $A_0=1$, and let $K_0=K(A_0,1)$ be given by
the second half of Theorem \ref{thm:MainContraction} for this value
of $A$. We will first prove the lemma for $|k|\geq K_0$. Let us fix
such a $k$.

Let $|z|\leq e^{\epsilon}$. The function $z^{r^{(N)}}$ belongs
to $\CC^{A_0,\epsilon}_N$ and its norm is bounded by $1$. For
$n\in \N$, we can iterate $n$ times \eqref{eq:GagneL4} (or
rather \eqref{eq:kljqmsklfdjqlmsdf}) (for $M=1$), to obtain
  \begin{equation}
  \label{eq:ControleL4utile}
  \norm{\Lp_{k,z}^n v}_{L^4} \leq \norm{\Lp_{k,z}^n v}_{D_k}
  \leq \theta^{100Nn} \norm{v}_{D_k}
  \leq \theta^{100Nn} \norm{v}_{C^1}.
  \end{equation}

We will then use \eqref{eq:ContractC1Iteres}. Note that the
function $\psi(x)=e^{-ik S_N^Y\phi_Y(x)} z^{r^{(N)}(x)}$ is
bounded by $e^{\epsilon r^{(N)}(x)}$, and for $h\in \HH_N$ we
have
  \begin{equation*}
  \norm{D(\psi \circ h)(x)}\leq |k|\norm{D (S_N^Y \phi_Y \circ
  h)(x)} e^{\epsilon r^{(N)}(x)}
  \leq C |k|r^{(N)}(x) e^{\epsilon r^{(N)}(x)}
  \leq C' |k| e^{2\epsilon r^{(N)}(x)}.
  \end{equation*}
Letting $A=C'|k|$, we have proved that $\psi \in
\CC^{A,2\epsilon}_N$ and $\norm{\psi}_{\CC^{A,2\epsilon}_N}\leq 1$.
Applying \eqref{eq:ContractC1Iteres} for $n$ iterates, we obtain,
for any $C^1$ function $w$,
  \begin{equation}
  \label{eq:PasMalC1}
  \norm{ \Lp_{k,z}^n w}_{C^1} \leq C|k| (\theta^{100Nn}
  \norm{w}_{C^1} + \theta^{-Nn} \norm{w}_{L^2}).
  \end{equation}
Applying this equation to $w=\Lp_{k,z}^n v$ and using
\eqref{eq:ControleL4utile}, we get
  \begin{equation}
  \norm{ \Lp_{k,z}^{2n} v}_{C^1} \leq C|k| (\theta^{100Nn}
  \norm{ \Lp_{k,z}^n v}_{C^1} + \theta^{-Nn} \theta^{100Nn}
  \norm{v}_{C^1}).
  \end{equation}
Applying once again \eqref{eq:PasMalC1} but this time to $v$, we
finally get $\norm{\Lp_{k,z}^{2n} v}_{C^1} \leq C |k|^2
\theta^{99Nn} \norm{v}_{C^1}$. We can argue in the same way for odd
times, to finally obtain the existence of $C$ such that, for any
$n\in \N$, $v\in C^1(Y)$,  $|k|\geq K_0$ and $|z|\leq e^{\epsilon}$,
  \begin{equation}
  \norm{\Lp_{k,z}^n v}_{C^1} \leq C k^2 \theta^{40Nn}
  \norm{v}_{C^1}.
  \end{equation}
This shows in particular that the operator $\Id-\Lp_{k,z}$ is
invertible on $C^1(Y)$, and that its inverse $\sum \Lp_{k,z}^n$ has
a norm which is bounded by $(C k^2)/(1-\theta^{40N})$.

We can then use Equation \eqref{ExprimeTnkCauchy} by taking for
$\gamma$ a circle of radius $e^{\epsilon}$. We obtain
  \begin{equation}
  \norm{T_{n,k}}\leq C k^2 \int_\gamma |z|^{-n}
  \leq C k^2 e^{-n \epsilon}.
  \end{equation}
This concludes the proof for $|k|\geq K_0$.

\emph{Second step.} Consider now $|k|<K_0$, $k\not=0$. We will show
that, for any $z$ with $|z|\leq 1$, the operator $\Id-\Lp_{k,z}$ is
invertible on $C^1(Y)$. Since the invertible operators form an open
set, this implies the existence of $\epsilon(k)$ such that, for
$|z|\leq e^{\epsilon(k)}$, $\Id-\Lp_{k,z}$ is invertible on
$C^1(Y)$. Using a path $\gamma$ which is a circle of radius
$e^{\epsilon(k)}$, we can then conclude as above (without explicit
control, but since there are only finitely many values of $k$ to
deal with this is not a problem).

Thus, consider $z$ with $|z|\leq 1$. The inequality
\eqref{eq:PasMalC1} still holds (its proof does not use
$|k|\geq K_0$). Therefore, there exists $C>0$ such that, for
any $n\in \N$, $\norm{\Lp_{k,z}^n v}_{C^1} \leq
C\theta^{100Nn}\norm{v}_{C^1} + C(n) \norm{v}_{L^2}$. Since the
injection of $C^1(Y)$ in $L^2(Y)$ is compact, this is a
Lasota-Yorke inequality. Hennion's Theorem \cite{hennion}
therefore shows that the essential spectral radius of
$\Lp_{k,z}$ is $<1$. If $\Id-\Lp_{k,z}$ is not invertible,
there must therefore exist $v\in C^1(Y)$ nonzero such that
$\Lp_{k,z}v=v$, i.e., $\Lp^N( e^{-ik S_N^Y\phi_Y} z^{r^{(N)}}
v)=v$. The operator $\Lp^N$ is the transfer operator associated
to the map $T_Y^N$, which is ergodic on $Y$. Lemma
\ref{lem:DecritFctPropre} applies and shows on the one hand
that $|z|^{r^{(N)}}$ is almost everywhere equal to $1$ (hence
$|z|=1$) and on the other hand that $v\circ T_Y^N=z^{r^{(N)}}
e^{-ik S_N^Y \phi_Y}v$ almost everywhere. Raising this equation
to the power $K_0$, we obtain that $v^{K_0}$ is invariant under
the operator $\Lp_{kK_0, z^{K_0}}$. But we have already proved
that $\Id-\Lp_{kK_0, z^{K_0}}$ is invertible on $C^1(Y)$. As a
consequence, $v^{K_0}=0$, and $v=0$, which is a contradiction.
This concludes the proof for $|k| \in [1,K_0)$.
\end{proof}

To obtain an estimate on $T_{n,0}$, we must also take into account
the fact that $\Id-\Lp_{0,1}$ is not invertible (its kernel
corresponds to constant functions), which will add a residue in the
integral calculus of the previous proof. In the following
definition, we introduce a tool which makes the computation of this
residue possible. We will write $\D$ for the open unit disk in $\C$,
and $\overline{\D}$ for its closure.

\begin{definition}
\label{def:renouvellementSimple} Let $\B$ be a Banach space, and let
$R_{j}$ be operators acting on $\B$, for $j>0$. We say that they
form a \emph{renewal sequence of operators with exponential decay}
if
\begin{enumerate}
\item There exist $\delta>0$ and $C>0$ such that $\norm{R_j}\leq C e^{-\delta
j}$. We can thus define an operator $R(z)=\sum R_jz^j$ for
$|z|<e^{\delta}$.
\item
For any $z\in \overline{\D} \moins \{1\}$, the operator $I-R(z)$ is
invertible on $\B$.
\item The operator $R(1)$ has a simple isolated eigenvalue at $1$. Let
$P=P(1)$ be the corresponding spectral projection, and $R'(1)=\sum j
R_{j}$. We assume that there exists $\mu>0$ such that $P R'(1)P=\mu
P$.
\end{enumerate}
\end{definition}

\begin{prop}
\label{prop:RenouvSimple} Let $R_j$ be a renewal sequence of
operators with exponential decay, on a Banach space $\B$. Let us
define an operator $T_n$ by $T_n=\sum_{p=1}^\infty
\sum_{j_1+\dots+j_p=n} R_{j_1}\dots R_{j_p}$. Then there exist $C>0$
and $\bar\theta<1$ such that, for any $n\in \N$, $\norm{T_n - P/\mu}
\leq C \bar\theta^n$.
\end{prop}
\begin{proof}
For $z$ close to $1$, the operator $R(z)$ is close to $R(1)$.
Hence, it has an eigenvalue $\lambda(z)$ close to $1$, with a
corresponding spectral projection $P(z)$ (and all these
quantities depend holomorphically on $z$). Let us compute the
derivative $\lambda'(1)$.

We will denote with a prime the derivative with respect to $z$. For
any $x\in \B$,  $R(z)P(z)x=\lambda(z)P(z)x$. Differentiating with
respect to $z$ and then multiplying on the left by $P(z)$, we get
(omitting the variable $z$)
  \begin{equation}
  PR'Px+ PRP'x=\lambda' Px+\lambda PP'x.
  \end{equation}
Moreover, $PRP'=P^2 RP'=PRPP'=\lambda PP'$. After simplification, we
obtain $PR'P x= \lambda' Px$. For $z=1$, $PR'P=\mu P$. Choosing $x$
such that $Px\not=0$, we finally get
  \begin{equation}
  \lambda'(1)=\mu\not=0.
  \end{equation}
In particular, on a small enough disk $O$ around $1$, the function
$z\mapsto \lambda(z)$ is injective, and takes the value $1$ only for
$z=1$.

The operators $I-R(z)$ are invertible for $z\in
\overline{\D}-O$, hence also for $z$ in a neighborhood of this
compact set. We can therefore choose a path $\gamma$ around $0$
going along an arc of a circle of radius $>1$, and the inner
part of $\partial O$. It satisfies the equation
  \begin{equation}
  \label{eq:ExprimeTnSimple}
  T_{n}=\frac{1}{2i\pi}\int_\gamma z^{-n-1} (I-R(z))^{-1}\dd z.
  \end{equation}
We modify $\gamma$ into a new path $\tilde\gamma$ which runs
along the same arc of circle of radius $>1$, and the outer part
of $\partial O$. To obtain an analogue of
\eqref{eq:ExprimeTnSimple}, we need to add the residue of
$z^{-n-1}(I-R(z))^{-1}$ inside $O$. We have
$(I-R(z))^{-1}=(1-\lambda(z))^{-1} P(z)+ Q(z)$ where $Q(z)$ is
holomorphic inside $O$ (whence without residue). The only pole
is thus at $1$, and we get
  \begin{equation}
  T_{n}=\frac{1}{2i\pi}\int_{\tilde \gamma}z^{-n-1} (I-R(z))^{-1}\dd z +
  \frac{1}{\lambda'(1)} P.
  \end{equation}
On $\tilde \gamma$, $|z|\geq e^{\delta'}$ for some $\delta'>0$.
As $\norm{(I-R(z))^{-1}}$ is uniformly bounded along
$\tilde\gamma$, the integral term is therefore $O(e^{-n
\delta'})$. The remaining term gives the conclusion of the
proposition.
\end{proof}

We can now come back to the study of the transfer operator
associated to $\UU$, and more precisely to the operators $T_{n,0}$,
which have not yet been estimated.

\begin{cor}
\label{cor:ExprimeTn0} For any $C^1$ function $v$ on $Y$, let
$Pv=\int v\dd\mu_Y$. Then there exist $C>0$ and $\bar\theta<1$ such
that, for any $n\in \N$ and any $v\in C^1(Y)$,
  \begin{equation}
  \norm{T_{n,0}v - \frac{1}{\mu^{(N)}(X^{(N)})}Pv}_{C^1} \leq C
  \bar\theta^n \norm{v}_{C^1}.
  \end{equation}
\end{cor}
\begin{proof}
We will use the fact that the Markov transformations $T_Y$ and
$U$ are mixing. Since these transformations are topologically
mixing (by the equality $\gcd\{r^{(N)}(x)\}=1$ for $U$), the
mixing in measure results e.g.\ from \cite[Theorem
4.4.7]{aaronson:book}.

Let us show that $R_{n,0}$ is a renewal sequence of operators
with exponential decay, on the Banach space $\B=C^1(Y)$. The
exponential decay of $\norm{R_{n,0}}$ is given by Lemma
\ref{lem:BorneRn}. Let $\Lp_{0,z} v = \Lp^N( z^{r^{(N)}} v)=
\sum z^n R_{n,0}=R(z)$.

Let us check that $I-R(z)=I-\Lp_{0,z}$ is invertible for $z\in
\overline{\D}\moins \{1\}$. As in the proof of Lemma
\ref{lem:TnkPourknot0}, the operators $\Lp_{0,z}$ (for $|z|\leq
1$) have an essential spectral radius $<1$ on $C^1$. If
$I-\Lp_{0,z}$ were not invertible, there would exist a nonzero
$C^1$ function $v$ such that $\Lp_{0,z}v=v$. Lemma
\ref{lem:DecritFctPropre} implies that $|z|=1$ and $v\circ
T_Y^N= z^{r^{(N)}}v$. Let us extend $v$ to the whole space
$X^{(N)}$ by setting $v(x,i)=z^{i} v(x,0)$. Thus, the function
$v$ is bounded (and therefore integrable), and satisfies
$v\circ U=zU$. This is a contradiction since $U$ is mixing.

For $z=1$, $R(1)=\Lp_{0,1}$ simply is the transfer operator
associated to $T_Y^N$. It has a simple eigenvalue at $1$ (the
corresponding spectral projection being $P$), and no other
eigenvalue of modulus $1$. Let us compute $PR'(1)P$. We have
  \begin{equation}
  PR_{n,0}P u = \mu_Y\{ r^{(N)}=n\} Pu.
  \end{equation}
As a consequence, Kac's Formula gives $PR'(1)P=\left(\sum
n\mu_Y\{ r^{(N)}=n\}\right)P= \mu^{(N)}(X^{(N)}) P$.

We can then apply Proposition \ref{prop:RenouvSimple} and get the
conclusion of the corollary.
\end{proof}

\subsection{The exponential mixing}

The estimates on $T_{n,k}$ given in the previous paragraph will
enable us to describe $\Mp_k^n$ for any $k$, and then the full
transfer operator $\hat \UU$.

For $x\in X^{(N)}$, denote by $h(x)$ its height in the tower
(i.e., if $x=(y,i)$ with $y\in Y$ and $i<r^{(N)}(x)$, let
$h(x)=i$). We will write $C^{5,1}(X^{(N)}\times \Sbb^1)$ for
the set of functions $v : X^{(N)}\times \Sbb^1\to \C$ such that
$\partial^i v/
\partial \omega^i$ is $C^1$ for $0\leq i\leq 5$, with its canonical
norm.

\begin{thm}
\label{thm:MelangeExponentielPrecis} There exist constants
$C>0$ and $\bar \theta<1$ such that, for any $C^{5,1}$ function
$v:X^{(N)}\times\Sbb^1 \to \C$, for any $n\in \N$ and any
$(x,\omega)\in X^{(N)}\times \Sbb^1$ with $h(x)\leq n/2$,
  \begin{equation}
  \left| \hat\UU^n v(x,\omega) -\int v \dd(\tilde
  \mu^{(N)}\otimes \Leb) \right| \leq C
  \bar\theta^{n} \norm{v}_{C^{5,1}}.
  \end{equation}
\end{thm}

For the proof, we will need information on the operators $T_{i,k}$,
but we also need to describe precisely the operators $B_{i,k}$
(defined in \eqref{definitBnk}).

\begin{lem}
\label{lem:EstimeBnk} There exist $\bar\theta<1$ and $C>0$ such
that, for any $k\in \Z$, $v\in C^1(X^{(N)})$ and $n\in \N$,
  \begin{equation}
  \label{Bnkbasique}
  \norm{B_{n,k}v}_{C^1}\leq C (1+|k|)\bar \theta^n
  \norm{v}_{C^1}.
  \end{equation}
Moreover,
  \begin{equation}
  \label{Bnkintegral}
  \left| \int_{X^{(N)}} v \dd \mu^{(N)} -\sum_{j=0}^n \int_Y B_{j,0}v
  \dd\mu^{(N)} \right|
  \leq C \bar\theta^n \norm{v}_{C^1}.
  \end{equation}
\end{lem}
\begin{proof}
For $y\in Y$, let $v_n(y)=0$ if $r^{(N)}(y)\leq n$, and
  \begin{equation}
  v_n(y)=v(y, r^{(N)}(y)-n) \exp\left(-ik
  \sum_{j=r^{(N)}(y)-n}^{r^{(N)}(y)-1}\phi(y,j)\right)
  \end{equation}
otherwise. For $x\in Y$, we then have $B_{n,k} v(x)= \Lp^N
v_n(x)$ since $B_{n,k}v(x)$ takes into account the values of
$v$ on the set $Z_n$ of points that enter $Y$ after exactly $n$
iterations, i.e., points of the form $(y,r^{(N)}(y)-n)$ with
$r^{(N)}(y)>n$.

Let us check that the function $v_n$ belongs to
$\CC^{1,\epsilon}_N$. First, since $v_n$ vanishes for
$r^{(N)}\leq n$, we have
  \begin{equation}
  |v_n(x)|\leq 1_{ r^{(N)}(x)>n}\norm{v}_{C^0}
  \leq e^{-\epsilon n}e^{\epsilon
  r^{(N)}(x)}\norm{v}_{C^0}.
  \end{equation}
Moreover, if $h\in \HH_N$,
  \begin{equation}
  \norm{D(v_n\circ h)(x)} \leq 1_{r^{(N)}\circ h>n}
  (\norm{v}_{C^1}+ kn \norm{v}_{C^0})
  \leq C(1+|k|) n e^{-\epsilon n} e^{\epsilon r^{(N)}(hx)}
  \norm{v}_{C^1}.
  \end{equation}
Hence, $v_n$ belongs to $\CC^{1,\epsilon}_N$ and its norm is
bounded by $C (1+|k|)\bar\theta^n \norm{v}_{C^1}$. Applying
\eqref{eq:ContractC1Iteres}, this yields \eqref{Bnkbasique}.

For \eqref{Bnkintegral}, note that $\sum_{j=0}^\infty \int_Y
B_{j,0}v = \int v$ since $\int_Y B_{j,0}v$ is the integral of $v$ on
$Z_j$. Therefore,
  \begin{equation}
  \left| \int v -\sum_{j=0}^n \int_Y B_{j,0}v
  \right|
  \leq \sum_{j=n+1}^\infty \left| \int_Y B_{j,0}v\right|
  \leq \sum_{j=n+1}^\infty \norm{B_{j,0}v}_{C^1}
  \leq C \bar\theta^n \norm{v}_{C^1}
  \end{equation}
by \eqref{Bnkbasique}.
\end{proof}

\begin{cor}
\label{cor:Mktend0} There exist $C>0$ and $\bar\theta<1$ such that,
for any $k\in \Z$, any $n\in \N$, any $x\in X^{(N)}$ with $h(x)\leq
n/2$, and any $v\in C^1(X^{(N)})$,
  \begin{equation}
  \label{eq:EstimeMpkn}
  \left| \Mp_k^n v(x) - 1_{k=0} \int v\dd\tilde
  \mu^{(N)}\right| \leq C (1+|k|^3)\bar\theta^{n}
  \norm{v}_{C^1}.
  \end{equation}
\end{cor}
\begin{proof}
Assume first that $x\in Y$. Then \eqref{eq:DonneMpnk} simply
becomes
  \begin{equation}
  \Mp^n_kv(x)=\sum_{i=0}^n T_{n-i,k} B_{i,k} v(x).
  \end{equation}
If $k\not=0$, then
  \begin{equation}
  \norm{T_{n-i,k} B_{i,k} v}_{C^1} \leq C k^2 \bar\theta^{n-i}
  \norm{B_{i,k}v}_{C^1}
  \leq C |k|^3 \bar\theta^{n-i} \bar \theta^i \norm{v}_{C^1},
  \end{equation}
by Lemmas \ref{lem:TnkPourknot0} and \ref{lem:EstimeBnk}. Summing
over $i$, we obtain the desired bound.

If $k=0$, Corollary \ref{cor:ExprimeTn0} gives an additional
term
  \begin{multline*}
  \sum_{i=0}^n P B_{i,0}v / \mu^{(N)}(X^{(N)})
  =\sum_{i=0}^n \int_Y B_{i,0}v
  \dd\mu^{(N)}/ \mu^{(N)}(X^{(N)})
  \\
  = \int v \dd\mu^{(N)} / \mu^{(N)}(X^{(N)}) + O(\bar\theta^n)
  = \int v \dd\tilde\mu^{(N)} + O(\bar\theta^n)
  \end{multline*}
by \eqref{Bnkintegral}. This proves \eqref{eq:EstimeMpkn} for $x\in
Y$.

If $x$ has height $j\in (0,n/2]$, let us write $x=U^j(x')$, so that
  \begin{equation}
  \Mp_k^n u(x)=e^{-ikS_j \phi(x')}\Mp_k^{n-j}u(x').
  \end{equation}
The estimate for $x'$ gives the desired conclusion (after replacing
$\bar\theta$ with $\bar\theta^{1/2}$).
\end{proof}

\begin{proof}[Proof of Theorem
\ref{thm:MelangeExponentielPrecis}] Let $v:X^{(N)}\times \Sbb^1 \to
\R$ be a $C^{5,1}$ function. We decompose it as
$v(x,\omega)=\sum_{k\in \Z} v_k(x)e^{ik \omega}$. Then
  \begin{equation}
  \hat \UU^n v(x,\omega)= \sum_{k\in \Z} \Mp_k^n v_k(x) \cdot
  e^{ik\omega},
  \end{equation}
by \eqref{EtudieUUn}. Therefore, if $h(x)\leq n/2$, Corollary
\ref{cor:Mktend0} gives
  \begin{align*}
  \left| \hat\UU^n v(x,\omega) -\int v \dd(\tilde
  \mu^{(N)}\otimes \Leb) \right|&
  \leq \left| \Mp_0^n v_0(x) -\int v_0 \dd\tilde \mu^{(N)}\right|+
  \sum_{k\not=0} \left|\Mp_k^n v_k(x)\right|
  \\&
  \leq C \sum_{k\in \Z} (1+|k|^3) \bar\theta^{n}
  \norm{v_k}_{C^1}.
  \end{align*}
With 5 integrations by parts with respect to $\omega$, we show
that $\norm{v_k}_{C^1} \leq C \norm{v}_{C^{5,1}}/ (1+|k|^5)$.
This implies the theorem after summation.
\end{proof}

\begin{proof}[Proof of Theorem
\ref{thm:MelangeExponentiel} (under the assumption $d^{(N)}=1$)]
\label{proof:DemontreMelangeExponentiel}
Let us first show that, on $X^{(N)}\times \Sbb^1$,
  \begin{equation}
  \label{eq:hatUUpetitL1}
  \norm{ \hat\UU^n v -\int v \dd(\tilde
  \mu^{(N)}\otimes \Leb)}_{L^1} \leq C \bar\theta^n
  \norm{v}_{C^{5,1}}
  \end{equation}
for some constants $C>0$ and $\bar\theta<1$. To do this, we
decompose $X^{(N)}$ as $\{ x \st h(x)> n/2 \}$ and $\{x \st h(x)\leq
n/2\}$. The first set has an exponentially small measure, its
contribution is therefore exponentially small. If $x$ belongs to the
second set, $\left|\hat \UU^n v(x,\omega)-\int v\right| \leq C
\bar\theta^{n} \norm{v}_{C^{5,1}}$ by Theorem
\ref{thm:MelangeExponentielPrecis}. This proves
\eqref{eq:hatUUpetitL1}.

This implies that, for any functions $v\in C^{5,1}$ and $u\in
L^\infty$,
  \begin{multline}
  \label{eq:MelangeEnHaut}
  \left| \int u\circ \UU^n \cdot v \dd(\tilde \mu^{(N)} \otimes
  \Leb) - \left( \int u \dd(\tilde \mu^{(N)} \otimes
  \Leb)\right) \left( \int v \dd(\tilde \mu^{(N)} \otimes
  \Leb)\right) \right|
  \\
  \leq C \bar\theta^n \norm{u}_{L^\infty} \norm{v}_{C^{5,1}}.
  \end{multline}
Take now $f\in L^\infty(X\times \Sbb^1)$ and $g\in C^6(X\times
\Sbb^1)$. The functions $u=f\circ \tilde\pi^{(N)}$ and
$v=g\circ \tilde \pi^{(N)}$ are defined on $X^{(N)}\times
\Sbb^1$, respectively bounded and in $C^{5,1}$. Moreover,
\eqref{eq:piContractePresque} shows that
$\norm{v}_{C^{5,1}}\leq C \norm{g}_{C^6}$. Since $\pi^{(N)}_*
\tilde \mu^{(N)}=\tilde \mu$, \eqref{eq:MelangeEnHaut} implies
  \begin{equation}
  \left| \int f\circ \TT^n \cdot g \dd(\tilde \mu \otimes
  \Leb) - \left( \int f \dd(\tilde \mu \otimes
  \Leb)\right) \left( \int g \dd(\tilde \mu \otimes
  \Leb)\right) \right|
  \leq C \bar\theta^n \norm{f}_{L^\infty} \norm{g}_{C^{6}}.
  \end{equation}

Let $n\in \N$ and $f\in L^\infty$. The linear operator
  \begin{equation}
  g\mapsto \int f\circ \TT^n \cdot g \dd(\tilde \mu \otimes
  \Leb) - \left( \int f \dd(\tilde \mu \otimes
  \Leb)\right) \left( \int g \dd(\tilde \mu \otimes
  \Leb)\right)
  \end{equation}
is then bounded by $2\norm{f}_{L^\infty}$ in $C^0$ norm, and by
$C \bar\theta^n \norm{f}_{L^\infty}$ in $C^6$ norm. For any
noninteger $\alpha\in (0,6)$, interpolation theory on the
compact manifold $X\times \Sbb^1$ (possibly with boundary)
shows that there exists a constant $C_\alpha$ such that any
operator which is bounded by $A$ in $C^0$ norm and by $B$ in
$C^6$ norm is then bounded by $C_\alpha A^{1-\alpha/6}
B^{\alpha/6}$ in $C^\alpha$ norm (see \cite[p.
200]{triebel:interpolation}). As a consequence, we get
  \begin{multline*}
  \left| \int f\circ \TT^n \cdot g \dd(\tilde \mu \otimes
  \Leb) - \left( \int f \dd(\tilde \mu \otimes
  \Leb)\right) \left( \int g \dd(\tilde \mu \otimes
  \Leb)\right) \right|
  \\
  \leq C_\alpha \norm{f}_{L^\infty}
  2^{1-\alpha/6} (C\bar\theta^n)^{\alpha/6}
  \norm{g}_{C^\alpha}.
  \end{multline*}
This concludes the proof of the theorem for noninteger
$\alpha$. The general case follows readily. The interpolation
argument can also be replaced by an elementary (but less
synthetic) convolution argument. The idea of using
interpolation theory in this kind of setting was suggested by
Dinh and Sibony.
\end{proof}

\begin{proof}[Proof of Theorem
\ref{thm:MelangeExponentiel} in the general case] If $d=d^{(N)}>1$,
the transformation $U$ is not mixing, and the arguments used above
(especially in the proof of Corollary \ref{cor:ExprimeTn0}) do not
apply any more.

However, they can be applied to the transformation $U^d$ and its
invariant measure $\tilde\mu^{(N)}_0$ (defined in Paragraph
\ref{subsec:modele}). As $\pi^{(N)}_* \tilde\mu^{(N)}_0= \tilde\mu$,
this implies Theorem \ref{thm:MelangeExponentiel} for times $n$ of
the form $kd$. To deduce the general case, one writes $n=kd+r$ with
$0\leq r<d$ and applies the theorem to the time $kd$ and to the
functions $f\circ \TT^r$ and $g$ (which are respectively bounded and
$C^\alpha$).
\end{proof}

\subsection{Proof of one implication in Proposition
\ref{prop:CaracteriseSigma2}}

\begin{prop}
\label{prop:CaracteriseSigma2alpha} Let $\psi: X\times \Sbb^1 \to
\R$ be a H\"{o}lder continuous function of $0$ average, and define
$\sigma^2$ by \eqref{eq:DefinitSigma2}. Then $\sigma^2\geq 0$.
Moreover, if $\sigma^2=0$, there exists a measurable function
$f:X\times \Sbb^1$, continuous on $Y\times \Sbb^1$, belonging to
$L^p$ for any $p<\infty$, such that $\psi=f-f\circ \TT$ almost
everywhere.
\end{prop}
This is one of the implications in Proposition
\ref{prop:CaracteriseSigma2}. Theorem \ref{thm:LimiteTCL} will be
required for the other half, hence its proof is postponed to
Paragraph \ref{subsec:Cohom}.

\begin{proof}
We have
  \begin{equation}
  \int_{X\times \Sbb^1}
  \left(\sum_{i=0}^{n-1} \psi\circ \TT^i
  \right)^2=n\int \psi^2+2\sum_{i=0}^{n-1} (n-i) \int \psi\cdot
  \psi\circ \TT^i.
  \end{equation}
Since $\sum_{i>0} i\left|\int \psi\cdot \psi\circ
\TT^i\right|<\infty$ by Theorem \ref{thm:MelangeExponentiel}, this
yields
  \begin{equation}
  \int_{X\times \Sbb^1}
  \left(\sum_{i=0}^{n-1} \psi\circ \TT^i
  \right)^2 = n\sigma^2 + O(1).
  \end{equation}
As a consequence, $\sigma^2\geq 0$. Moreover, if $\sigma^2=0$,
the Birkhoff sums of $\psi$ are uniformly bounded in $L^2$. By
\cite{vitesse_birkhoff}, there exists an $L^2$ function $f$
with zero average such that $\psi=f-f\circ \TT$ almost
everywhere. We have to prove that $f$ is continuous on $Y\times
\Sbb^1$ and belongs to every $L^p$, $p<\infty$.

Theorem \ref{thm:MelangeExponentielPrecis} implies that there
exist $\bar\theta<1$ and $C>0$ such that, for any $C^6$
function $v:X\times \Sbb^1 \to \C$, for any $n\in \N$, for any
$x\in X^{(N)}$ with $h(x)\leq n/2$,
  \begin{equation}
  \left| \hat \UU^n (v\circ \tilde\pi^{(N)})(x,\omega) -\int v\right| \leq C
  \bar\theta^n \norm{v}_{C^6}.
  \end{equation}
Since $|\hat \UU^n (v\circ \tilde\pi^{(N)})(x,\omega) -\int v|
\leq 2\norm{v}_{C^0}$, interpolation theory as above implies
that, for any $\alpha>0$, there exist $C_\alpha>0$ and
$\bar\theta_\alpha<1$ such that, for any $x\in X^{(N)}$ with
$h(x)\leq n/2$,
  \begin{equation}
  \left| \hat \UU^n (v\circ \tilde\pi^{(N)})(x,\omega) -\int v\right| \leq
  C_\alpha  \bar\theta_\alpha^n \norm{v}_{C^\alpha}.
  \end{equation}

As $\psi$ belongs to $C^\alpha$ and has vanishing integral, we can
therefore define a function $g$ on $X^{(N)}\times \Sbb^1$ by
  \begin{equation}
  g(x,\omega)=-\sum_{n=1}^\infty \hat\UU^n(\psi\circ
  \tilde\pi^{(N)})(x,\omega).
  \end{equation}
This function is continuous on $Y\times \Sbb^1$, and belongs to
$L^p$ for any $p<\infty$ (since $|g(x,\omega)|\leq C(1+h(x))$,
this last function belonging to any $L^p$ because $\mu^{(N)}\{
h(x)\geq n\}$ decays exponentially with $n$). Moreover, by
construction, $\hat\UU g-g=\hat\UU(\psi \circ
\tilde\pi^{(N)})$.

We know that $\psi=f-f\circ \TT$ where $f\in L^2$. As a
consequence, $\psi\circ \tilde\pi^{(N)}=f\circ \tilde\pi^{(N)}-
f\circ \tilde\pi^{(N)} \circ \UU$, whence
$\hat\UU(\psi\circ\tilde\pi^{(N)})= \hat\UU(f\circ
\tilde\pi^{(N)})-f\circ \tilde\pi^{(N)}$. We get
  \begin{equation}
  g-f\circ \tilde\pi^{(N)}=\hat\UU(g-f\circ \tilde\pi^{(N)}).
  \end{equation}
In particular, for any $n\in \N$, $g-f\circ
\tilde\pi^{(N)}=\hat\UU^n(g-f\circ \tilde\pi^{(N)})$.

Theorem \ref{thm:MelangeExponentielPrecis} shows that, for any
function $v\in C^{5,1}(X^{(N)}\times \Sbb^1)$ with zero
integral, $\hat\UU^n v$ converges to $0$ in $L^2$. By density,
this convergence holds for any function $v\in L^2$ with zero
integral. In particular, $\hat\UU^n(g-f\circ \tilde\pi^{(N)})$
converges to $0$, hence $g-f\circ\tilde\pi^{(N)}=0$. As $g$ is
continuous on $Y\times\Sbb^1$ and belongs to all spaces $L^p$,
$p<\infty$, this concludes the proof.
\end{proof}

\section{Strategy and tools for the local limit theorem}
\label{sec:Strategie}

\subsection{Description of the strategy of the proof}
\label{subsec:Strategie} Let us fix an integer $M$. We work with the
transformation $U=U^{(MN)}$ on $X^{(MN)}$ (hence also with
$\UU^{(MN)}$ on $X^{(MN)}\times \Sbb^1$).

Let $\psi:X\times \Sbb^1 \to\R$ be a $C^6$ function with $0$
average. We will also write $\psi$ instead of $\psi\circ
\tilde\pi^{(MN)}$ on $X^{(MN)}\times \Sbb^1$. To prove the
local limit theorem for $\psi$, we consider for $t\in \R$ the
operator $\hat\UU_t (v):= \hat\UU(e^{it\psi} v)$. If we
understand well the iterates of $\hat\UU_t$, we will deduce the
asymptotic behavior of $\int e^{it S_n\psi}$, since this
quantity is equal to $\int \hat\UU_t^n(1)$.

Instead of working with functions on $X^{(MN)}\times \Sbb^1$, we
have seen in the proof of the exponential mixing that it is
worthwhile to use Fourier series, and work on $X^{(MN)}\times \Z$.
If $v$ is a function and $(v_k)_{k\in\Z}$ denote its Fourier
coefficients, then the Fourier coefficients of $e^{it\psi}v$ are
given by
  \begin{equation}
  (e^{it\psi} v)_k = \sum_{a+b=k} (e^{it\psi})_a v_b.
  \end{equation}
Applying then the operator $\hat\UU$ (which acts at the level of the
$k$ frequency by the operator $\Mp_k$), we obtain
  \begin{equation}
  (\hat \UU_t v)_k (x)= \sum_{l\in \Z} \sum_{Ux'=x} \KK(x') e^{-ik
  \phi(x')} (e^{it\psi})_{k-l}(x') v_l(x').
  \end{equation}
This is some kind of Markov operator on $X^{(MN)}\times \Z$, for the
``transition probability''
  \begin{equation}
  \K^t_{(x,k)\to (x',l)}:= 1_{Ux'=x} \KK(x') e^{-ik
  \phi(x')} (e^{it\psi})_{k-l}(x').
  \end{equation}
The equality $\sum_{(x',l)} \K_{ (x,k)\to (x',l)}=1$ does not
hold, so this is not a real transition kernel, but we will
nevertheless use the intuition of random walks. Let us in
particular write, for $n\in \N$,
  \begin{equation}
  \label{DonneKn1}
  \K^{t,n}_{(x,k)\to (x',l)}=\sum_{\substack{k_0=l,k_1,\dots, k_{n-1},k_n=k\\
  x_0=x',x_1,\dots,x_{n-1}, x_n=x}} \K^t_{(x_n,k_n)\to (x_{n-1},k_{n-1})}
  \dots \K^t_{(x_2,k_2)\to (x_1,k_1)} \K^t_{(x_1, k_{1})\to
  (x_0, k_0)}.
  \end{equation}
In this expression, we consider trajectories of the random walk
$x_n,x_{n-1},\dots, x_0$. It may seem unnatural to write things
in that direction, but it is designed to give the ``good''
order when we express things in terms of transfer operators.
Let $\hat\K^t$ be the operator with kernel $\K^t$, acting on
bounded functions on $X^{(MN)}\times \Z$, by
  \begin{equation}
  \hat \K^tv(x,k)=\sum_{ (x',l)}\K^t_{(x,k)\to (x',l)} v(x',l).
  \end{equation}
By construction, the powers $\hat \K^{t,n}$ of $\hat \K^t$ have
kernels $\K^{t,n}$. Moreover, $\hat \UU_t$ corresponds to the
operator $\hat\K^t$ at the level of frequencies, i.e., if $v$
is a smooth function on $X^{(MN)}\times \Sbb^1$ with Fourier
coefficients $(v_k)_{k\in\Z}$,
  \begin{equation}
  (\hat \UU_t ^n v)_k (x)=
  \sum_{(x',l)} \K^{t,n}_{(x,k)\to (x',l)}
  v_l(x').
  \end{equation}
To see that this expression and these computations are correct, we
should check that
  \begin{equation}
  \sup_{(x,k)\in X^{(MN)}\times \Z} \sum_{(x',l)} \left|
  \K^t_{(x,k)\to (x',l)} \right|<\infty,
  \end{equation}
which is always the case if $\psi$ is $C^2$ in the direction of
$\Sbb^1$ (by two integrations by parts), and will always be
satisfied in the following. A priori, this does not prevent
$\K^{t,n}_{(x,k)\to (x',l)}$ from blowing up exponentially fast
with $n$. However, $\K^{t,n}_{(x,k)\to (x',l)}$ is also the
kernel of the operator obtained by multiplying $v$ with $e^{it
S_n \psi}$, and then applying $\hat \UU^n$. Therefore,
  \begin{equation}
  \label{DonneKn}
  \K^{t,n}_{(x,k)\to (x',l)}=1_{U^n x'=x} \KK^{(n)}(x') e^{-ik S_n
  \phi(x')} (e^{it S_n \psi})_{k-l}(x'),
  \end{equation}
and this quantity is bounded by $\KK^{(n)}(x')\leq 1$. Note that
\eqref{DonneKn} can also be checked directly from the formula
\eqref{DonneKn1}, with several successive integrations.

We will let different operators (with kernels related to $\K^{t,n}$)
act on spaces of functions from $X^{(MN)} \times \Z$ to $\C$  (or
$Y\times \Z$ to $\C$ if we only consider trajectories starting from
$Y\times \Z$ or ending in $Y\times \Z$). If $\BB$ is such a
functional space, and $v\in \BB$, we will sometimes write $v_k(x)$
instead of $v(x,k)$.

To understand the previous ``random walk'', we will study its
successive returns to the set $Y \times [-K,K]$ where $K$ is
large enough. Indeed, outside of this set, we have a strong
contraction (by Theorem \ref{thm:MainContraction}) hence
excursions can be controlled. Only what happens inside $Y\times
[-K,K]$ can therefore be problematic, and we will use there an
abstract compactness argument. Let us denote by
$\K^{t,n,exc}_{(x,k)\to (x',l)}$ the ``probability'' of an
excursion, i.e., of starting from $(x,k)\in Y\times [-K,K]$,
and coming back to $(x',l)\in Y\times [-K,K]$ after a time
exactly $n$, without entering $Y\times [-K,K]$ in between.
Formally, for $(x,k) \in Y\times [-K,K]$ and $(x',l)\in Y\times
[-K,K]$,
  \begin{equation*}
  \K^{t,n,exc}_{(x,k)\to (x',l)}=\sum_{\substack{k_0=l,\dots, k_n=k\\
  x_0=x',x_1,\dots,x_{n-1}\in X, x_n=x\\ (x_i, k_i)\not\in Y\times [-K,K]\text{ for }0<i<n}}
  \K^t_{(x_n,k_n)\to (x_{n-1},k_{n-1})}\dots
  \K^t_{(x_2,k_2)\to (x_1,k_1)} \K^t_{(x_1, k_1)\to
  (x_0, k_0)}.
  \end{equation*}

Let $\B_K=\bigoplus_{|k|\leq K} C^1(Y)$. An element of $\B_K$ can
therefore be seen as a function $v$ on $X\times \Z$ such that $v_k$
is $C^1$ for $|k|\leq K$, and $v_k=0$ for $|k|>K$. We define then an
operator $R^t_n$ on $\B_K$ by
  \begin{equation}
  \label{eq:DefinitRtn}
  (R^t_n v)_k(x)=\sum_{(x',l)} \K^{t,n,exc}_{(x,k)\to (x',l)}
  v_l(x').
  \end{equation}
For $x\in Y$ and $|k|\leq K$, let also $(T^t_n v)_k (x)=
\sum_{(x',l)\in Y\times [-K,K]} \K^{t,n}_{(x,k)\to (x',l)}
v_l(x')$, i.e., we consider all the returns of the ``random
walk'' to $Y\times[-K,K]$ and not only the first ones. This
means that $T^t_n v= 1_{Y\times [-K,K]} \hat\K^{t,n}(1_{Y\times
[-K,K]} v)$ for $v\in \BB_K$. By construction,
  \begin{equation}
  \label{eq:DefinitTtn}
  T^t_n=\sum_{p=1}^\infty \sum_{j_1+\dots+j_p=n} R^t_{j_1}\dots
  R^t_{j_p}.
  \end{equation}
This is a renewal equation, that we already met in the course of the
proof of exponential mixing. The main difference is that, for the
mixing, each frequency was left invariant by the transfer operator,
which means we only had to consider random walks on $X^{(N)}$ and
excursions outside $Y$. Here, since there is also some interaction
between the frequencies, we have to localize spatially (i.e., on
$Y$), but also on the space of frequencies since the estimates given
by Theorem \ref{thm:MainContraction} are not uniform in $k$.

The proof will consist in understanding precisely the $R^t_n$'s,
deducing from that good estimates on $T^t_n$'s, and using these to
reconstruct precisely enough $\hat \UU_t^n$. We will thus need two
technical tools: on the one hand, a tool on perturbations of renewal
sequences of operators (we want estimates which are precise both
with respect to $n$ and $t$), and on the other hand good estimates
on the excursions outside of  $Y\times [-K,K]$.

Before going on, let us give another expression of $\K^{t,n,exc}$
that will be needed later on, by considering the successive returns
to $Y\times \Z$. Let us define a function $\psi_Y : Y\times
\Sbb^1\to \R$ by
  \begin{equation}
  \label{DefinitPsiY}
  \psi_Y(x,\omega)=\sum_{i=0}^{r(x)-1} \psi\left(T^i x, \omega+
  \sum_{j=0}^{i-1} \phi(T^j x)\right).
  \end{equation}
It is the function induced by $\psi$ and $\TT$ on the set $Y\times
\Sbb^1$. Let us denote by $S^Y_n \psi_Y$ the Birkhoff sums of
$\psi_Y$ for the map induced by $\TT$ on $Y\times \Sbb^1$. For
$x,x'\in Y$ and $k,l\in \Z$, let $\K^{t,Y}_{(x,k)\to
(x',l)}=1_{T^{MN}x'=x} J^{(MN)}(x') e^{-ik S^Y_{MN} \phi_Y(x')}
(e^{it S^Y_{MN} \psi_Y})_{k-l}(x')$, which corresponds to the
``probability'' (for the above random walk) of the first return in
$Y\times \Z$. Considering the successive returns to $Y\times (\Z
\moins [-K,K])$, we get for $x,x'\in Y$ and $k,l\in [-K,K]$,
  \begin{equation}
  \label{eq:RedonneKtnexc}
  \K^{t,n,exc}_{(x,k)\to (x',l)}
  =\sum_{p\geq 0}\ \sum_{\substack{k_0=l, k_1,\dots,k_{p-1} \not\in [-K,K], k_p=k \\
  x_0=x',x_1,\dots, x_{p-1}\in Y,x_p=x\\
  \sum_{i=0}^{p-1} r^{(MN)}(x_i)=n}} \K^{t,Y}_{(x_p,k_p)\to
  (x_{p-1},k_{p-1})} \dots \K^{t,Y}_{(x_1, k_1)\to (x_0,k_0)}.
  \end{equation}

\subsection{Perturbed renewal sequences of operators}

\begin{definition}
\label{def:renouvellement} Let $\B$ be a Banach space, and let
$R_j^t$ be operators acting on $\B$, for $j>0$ and $t\in [-t_0,t_0]$
for some $t_0>0$. These operators form a \emph{perturbed sequence of
renewal operators with exponential decay} if
\begin{enumerate}
\item The operators $R^0_j$ form a renewal sequence of operators with
exponential decay. We will in particular write $P$ and $\mu$ for the
associated spectral projection and coefficient, as in Definition
\ref{def:renouvellementSimple}.
\item There exist $\delta>0$ and $a,C>0$ such that, for all
$t,t'\in[-t_0,t_0]$ with $|t-t'|\leq a$, for any $j>0$,
$\norm{R^t_j-R^{t'}_j}\leq C |t-t'|e^{-\delta j}$.
\item Let us write $R(z,t)=\sum z^j R^t_j$ for $|z|<e^{\delta}$.
For $(z,t)$ close to $(1,0)$, the operator $R(z,t)$ is a small
perturbation of $R(1,0)$. Therefore, it has an eigenvalue
$\lambda(z,t)$ close to $1$. We assume that, for some
$\alpha>0$, $\lambda(1,t)=1-\alpha t^2+O(|t|^3)$.
\end{enumerate}
We say that this sequence if \emph{aperiodic} if, for any $(z,t)\in
(\overline{\D}\times [-t_0,t_0]) \moins \{(1,0)\}$, the operator
$I-R(z,t)$ is invertible on $\B$.
\end{definition}

\begin{thm}
\label{thm:renouvellement} Let $R^t_j$ be a perturbed sequence of
renewal operators with exponential decay. Let
  \begin{equation}
  T^t_n=\sum_{p=1}^\infty \sum_{j_1+\dots+j_p=n}R^t_{j_1}\dots
  R^t_{j_p}.
  \end{equation}
Then there exist $\eps\in (0,t_0)$, $\bar\theta<1$ and $c,C>0$
such that, for $t\in[-\eps,\eps]$, for $n>0$,
  \begin{equation}
  \label{eq:ControlePres0}
  \norm{T^t_n -
  \frac{1}{\mu}\left(1-\frac{\alpha t^2}{\mu}\right)^n P} \leq
  C\bar\theta^n+ C |t|(1-ct^2)^n.
  \end{equation}
Moreover, if $R^t_j$ is aperiodic, one also has, for $|t| \in
[\eps,t_0]$ and $n>0$,
  \begin{equation}
  \label{eq:ControleLoin0}
  \norm{T^t_n}\leq C \bar\theta^n.
  \end{equation}
\end{thm}
\begin{proof}
If $\gamma$ is a path around $0$ in $\C$, close enough to $0$,
  \begin{equation}
  \label{eq:ExprimeTjt}
  T_j^t=\frac{1}{2i\pi}\int_{\gamma} z^{-j-1} (I-R(z,t))^{-1} \dd z.
  \end{equation}
By analyticity, this equality holds true for any path $\gamma$
around $0$ bounding a domain on which $I-R(z,t)$ is invertible for
any $z$.

Let us first show \eqref{eq:ControleLoin0} in the aperiodic
case. Let $t\not=0$. The operators $I-R(z,t)$ are invertible
for any $z\in \overline{\D}$. Since invertible operators form
an open set, there exists an open neighborhood $I_t$ of $t$,
and $\epsilon_t>0$, such that $I-R(z,t')$ is invertible for
$t'\in I_t$ and $|z|\leq e^{\epsilon_t}$. Taking for $\gamma$
the circle of radius $e^{\epsilon_t}$, we obtain
$\norm{T_j^{t'}} \leq C(t) e^{-j\epsilon_t}$. If $\tau>0$, the
compact set $[-t_0,-\tau]\cup [\tau,t_0]$ can be covered by a
finite number of the intervals $I_t$, and we get the following:
there exist $\delta_\tau>0$ and $C_\tau>0$ such that, for any
$|t|\in [\tau,t_0]$, for any $j>0$, $\norm{T_j^t} \leq C_\tau
e^{-j \delta_\tau}$. This proves \eqref{eq:ControleLoin0}, if
we can choose $\tau$ so that \eqref{eq:ControlePres0} is
satisfied.

For \eqref{eq:ControlePres0}, we work in a neighborhood of
$(z,t)=(1,0)$. There exist an open disk $O$ around $1$, and
$\eps>0$, such that, for $(z,t)\in O\times [-\eps,\eps]$, the
operator $R(z,t)$ has a unique eigenvalue $\lambda(z,t)$ close
to $1$. Let us also denote by $P(z,t)$ the corresponding
spectral projection. These functions depend holomorphically on
$z$, and in a Lipschitz way on $t$.

We saw in the proof of Proposition \ref{prop:RenouvSimple} that
$\lambda'(1,0)=\mu\not=0$. Reducing $O$ if necessary, we can
therefore assume that $z\mapsto \lambda(z,0)$ is injective on $O$
(and takes the value $1$ only at $z=1$).

When $t$ converges to $0$, the function $z\mapsto \lambda(z,t)$
converges uniformly to $z\mapsto \lambda(z,0)$ (with a speed
$O(t)$). Since all these functions are holomorphic, the derivatives
converge uniformly with the same speed. In particular, $z\mapsto
\lambda(z,t)$ takes the value $1$ at a unique point $\gamma(t)$ in
$O$, if $t$ is small enough, by Rouch\'{e}'s Theorem. Moreover,
$\gamma(t)\to 1$ when $t\to 0$.

Let us establish an asymptotic expansion of $\gamma(t)$. We have
  \begin{align*}
  \lambda(\gamma(t),t)-\lambda(1,t)=\int_{1}^{\gamma(t)}
  \lambda'(z,t) \dd z
  =\int_1^{\gamma(t)} (\lambda'(z,t)-\lambda'(1,0))\dd z +
  \lambda'(1,0)(\gamma(t)-1).
  \end{align*}
Moreover, $|\lambda'(z,t)-\lambda'(1,0)|\leq C (|z-1|+|t|) \leq
C( |\gamma(t)-1|+|t|)$. As $\lambda(\gamma(t),t)-\lambda(1,t)=1
- \lambda(1,t)=\alpha t^2+O(|t|^3)$, we obtain
  \begin{equation}
  \lambda'(1,0) (\gamma(t)-1) = \alpha t^2+O(t^3) +O( |t|
  |\gamma(t)-1|) + O(|\gamma(t)-1|^2).
  \end{equation}
As $\lambda'(1,0)=\mu \not=0$, this yields $\gamma(t)-1 \sim
\alpha t^2/\mu$. In particular, $\gamma(t)-1=O(t^2)$. Putting
this information back in the equation, we finally obtain
  \begin{equation}
  \label{exprimegammat}
  \gamma(t)=1+\alpha t^2/\mu + O(t^3).
  \end{equation}

The operators $I-R(z,0)$ are invertible for $z\in
\overline{\D}-O$. By continuity, $I-R(z,t)$ is invertible for
any $z$ in a neighborhood of this compact set, and $t$ close
enough to $0$, say $t\in[-\eps,\eps]$. We can therefore choose
a path $\gamma$ around $0$ made of an arc of circle of radius
$>1$, and the inner part of $\partial O$, satisfying
\eqref{eq:ExprimeTjt} for $|t|\leq \eps$. We modify $\gamma$
into a new path $\tilde\gamma$ by replacing the inner part of
$\partial O$ with its outer part. To obtain an analogue of
\eqref{eq:ExprimeTjt}, we should add the residue of
$z^{-j-1}(I-R(z,t))^{-1}$ inside $O$. We have
$(I-R(z,t))^{-1}=(1-\lambda(z,t))^{-1} P(z,t)+ Q(z,t)$ where
$Q(z,t)$ is holomorphic inside $O$ (whence without residue).
The only pole is located at $\gamma(t)$, and we obtain
  \begin{equation}
  T_j^t=\frac{1}{2i\pi}\int_{\tilde \gamma}z^{-j-1} (I-R(z,t))^{-1}\dd z +
  \frac{1}{\lambda'(\gamma(t),t)} P(\gamma(t),t)
  \gamma(t)^{-j-1}.
  \end{equation}
On $\tilde \gamma$, we have $|z|\geq e^{\delta_0}$ for some
$\delta_0>0$. As $\norm{(I-R(z,t))^{-1}}$ is uniformly bounded
on $\tilde\gamma$, the integral term is $O(e^{- \delta_0 j})$.
For the remaining term, we have
$\frac{1}{\lambda'(\gamma(t),t)} P(\gamma(t),t) =
\frac{1}{\lambda'(1,0)} P(1,0)+O(t)$. Making this substitution
gives an error of $O(|t| |\gamma(t)|^{-j})=O(|t| (1-ct^2)^j)$,
by \eqref{exprimegammat}. We get
  \begin{equation}
  \norm{T_j^t- \frac{1}{\mu} P \gamma(t)^{-j-1}} \leq C
  e^{-j\delta_0}+ C|t| (1-ct^2)^j.
  \end{equation}
Finally, if we replace $\gamma(t)^{-j-1}$ with $(1-\alpha
t^2/\mu)^j$, the error is bounded, thanks to
\eqref{exprimegammat}, by
  \begin{equation*}
  C (1-ct^2)^j \bigl( (1+C|t|^3)^j-1\bigr)\leq C(1-ct^2)^j (1+C|t|^3)^j j|t|^3.
  \end{equation*}
If $t$ is small enough, $(1-ct^2)(1+C|t|^3) \leq (1-ct^2/2)$.
Finally,
  \begin{equation}
  \label{eq:TechniqueControlejttrois}
  \begin{split}
  j|t|^3 (1-ct^2/2)^j &\leq j|t|^3 (1-ct^2/4)^j (1-ct^2/4)^j \leq
  |t| (1-ct^2/4)^j \cdot jt^2 \exp(-cjt^2/4) \\&\leq C |t|
  (1-ct^2/4)^j,
  \end{split}
  \end{equation}
since the function $x\mapsto xe^{-cx^2/4}$ is bounded on
$\R_+$.
\end{proof}

\subsection{Estimates on the excursions}
\label{sec:EstimeExcursion}

In this whole paragraph, we fix an integer $M$, a constant $A>1$ and
a sequence $(\gamma_d)_{d\in \Z}$ with $\gamma_d\in (0,1]$ and
$\gamma_d=O(1/|k|^4)$ when $d\to \pm \infty$.

We then choose an integer $K$ such that
  \begin{equation}
  \label{Kgrand1}
  \forall |d|> K/2,\quad \gamma_d \leq \frac{1}{(1+|d|)^{60/17}},
  \end{equation}
and
  \begin{equation}
  \label{Kgrand2}
  K \geq K(A,M) \text{ given by Theorem
  \ref{thm:MainContraction}}.
  \end{equation}
and
  \begin{equation}
  \label{Kgrand3}
  \forall n\geq 1,\quad
  2^{Mn} - 1- n/2 \geq 2^{Mn} / K.
  \end{equation}

Let $\kk=(k_0,k_1,\dots,k_j)$ be a sequence of integers. We say that
this sequence is admissible if $|k_i|>K$ for any $i\in (0,j)$. We
say that it is strongly admissible if, additionally, $|k_j|>K$. We
will denote by $d_i=k_i-k_{i-1}$ the successive differences.

\begin{lem}
\label{Lem:ContractionL2} Let $\kk=(k_0,k_1,\dots,k_{j_0})$ be a
strongly admissible sequence. Let $\psi_1,\dots,\psi_{j_0}$ be
functions from $Y$ to $\C$, and let
$\epsilon_1,\dots,\epsilon_{j_0}$ belong to $[0,1]$. Assume that
$\norm{\psi_i}_{\CC^{A,3\epsilon}_{MN}}\leq \epsilon_i\gamma_{d_i}$.

Let $v^0:Y\to \C$, define a sequence of functions $v^i$ by
induction, by $v^i=\Lp^{MN}_{k_i}(\psi_i v^{i-1})$. Then
  \begin{equation}
  \label{eq:EstimeeL2}
  \norm{v^{j_0}}_{L^2} \leq \left(\prod_{i=1}^{j_0} \epsilon_i \gamma_{d_i}^{9/10}\right)
  \theta^{100MNj_0} \norm{v^0}_{C^1}.
  \end{equation}
\end{lem}
\begin{proof}
We will use the following ``virtual heights''
  \begin{equation}
  \beta_i=\max(|k_i|,
  |k_{i-1}|/2^M, \dots, |k_0|/2^{Mi}).
  \end{equation}
Their interest is that we will be able to control by induction the
Dolgopyat norms $\norm{ v^i}_{D_{\beta_i}}$ (while this would not be
possible for the norm $D_{k_i}$ if the jumps $d_i$ are too large).

If $|k_i|\geq \beta_{i-1}/2^M$, we have $\beta_i=|k_i|$. Then, by
Theorem \ref{thm:MainContraction} (and more precisely
\eqref{eq:GagneL4}),
  \begin{equation*}
  \norm{ v^i}_{D_{\beta_i}}
  =\norm{ \Lp_{k_i}^{MN}( \psi_i v^{i-1})}_{D_{k_i}}
  \leq \theta^{100 MN} \norm{ \psi_i}_{\CC^{A,3\epsilon}_{MN}}
  \norm{v^{i-1}}_{D_{2^Mk_i}}
  \leq \theta^{100 MN} \epsilon_i \gamma_{d_i} \norm{v^{i-1}}_{D_{\beta_{i-1}}}.
  \end{equation*}
Otherwise, $\beta_i=\beta_{i-1}/2^M > |k_i|$, and (using
\eqref{eq:PerdPasTrop})
  \begin{equation}
  \label{AppliqueL2}
  \norm{ v^i}_{D_{\beta_i}}
  = \norm{ \Lp^{MN}_{k_i}( \psi_i v^{i-1})}_{D_{\beta_i}}
  \leq \theta^{-MN}\norm{ \psi_i}_{\CC^{A,3\epsilon}_{MN}}
  \norm{v^{i-1}}_{D_{2^M\beta_i}}
  \leq \theta^{-MN} \epsilon_i \gamma_{d_i} \norm{v^{i-1}}_{D_{\beta_{i-1}}}.
  \end{equation}
In both cases, we have similar equations, with a large gain or a
small loss.

Let us show by induction on $i$ that
  \begin{equation}
  \label{ToBeProvedRec}
  \norm{v^i}_{D_{\beta_i}} \leq \theta^{100 MNi} \epsilon_1\dots\epsilon_i(\gamma_{d_1}\dots
  \gamma_{d_i})^{9/10} \norm{v^0}_{D_{k_0}},
  \end{equation}
the result being clear for $i=0$.

Assume that the result is proved up to $i-1$, and let us prove it
for $i$. If $\beta_i=|k_i|$,
  \begin{equation}
  \norm{v^i}_{D_{\beta_i}} \leq \theta^{100 MN} \epsilon_i\gamma_{d_i}
  \norm{v^{i-1}}_{D_{\beta_{i-1}}}
  \leq \theta^{100 MN}\epsilon_i (\gamma_{d_i})^{9/10} \norm{v^{i-1}}_{D_{\beta_{i-1}}}
  \end{equation}
since $\gamma_d \leq 1$ for any $d\in\Z$. The inductive
assumption concludes the proof.

If $\beta_i>|k_i|$, consider $\iota$ the last time before $i$
for which $\beta_\iota=|k_\iota|$. Iterating \eqref{AppliqueL2}
up to $\iota$, we get
  \begin{equation}
  \label{BaseRecurr}
  \norm{v^i}_{D_{\beta_i}} \leq \epsilon_i\dots\epsilon_{\iota+1}\gamma_{d_i}\dots
  \gamma_{d_{\iota+1}}\theta^{-MN(i-\iota)}
  \norm{ v^\iota}_{D_{\beta_\iota}}.
  \end{equation}
Moreover, $\beta_i=\beta_\iota / 2^{M(i-\iota)}$, and
$\beta_i>K$ since $\kk$ is strongly admissible. Hence,
  \begin{equation}
  |d_{\iota+1}|+\dots+|d_i| \geq |k_\iota-k_i| \geq (2^{M(i-\iota)}-1)
  \beta_i \geq (2^{M(i-\iota)}-1)K.
  \end{equation}
Write $J$ for the set of indexes $a\in(\iota,i]$ for which
$|d_a|> K/2$. Then $\sum_J |d_a| \geq
(2^{M(i-\iota)}-1-(i-\iota)/2)K$. By \eqref{Kgrand3}, we
therefore get $\sum_J |d_a| \geq 2^{M(i-\iota)}$. By
\eqref{Kgrand1}, $\gamma_d \leq 1/(1+|d|)$ for any $|d|>K/2$.
We obtain
  \begin{align*}
  (\gamma_{d_i}\dots \gamma_{d_{\iota+1}})^{1/10} & \leq \prod_{a\in J}
  \gamma_{d_a}^{1/10} \leq \prod_{a\in J} \frac{1}{(1+|d_a|)^{1/10}}
  = \left( \frac{1}{\prod_{a\in J} (1+|d_a|)}\right)^{1/10}
  \\&
  \leq \left( \frac{1}{\sum_{a\in J} |d_a| }\right)^{1/10}
  \leq 2^{- M(i-\iota)/10}.
  \end{align*}
By Theorem \ref{thm:MainContraction}, $\theta^{101 N} \geq
2^{-1/10}$. As a consequence, $2^{- M(i-\iota)/10} \leq \theta^{101
MN(i-\iota)}$. Hence, we obtain from \eqref{BaseRecurr}
  \begin{align*}
  \norm{v^i}_{D_{\beta_i}} &\leq \theta^{-MN(i-\iota)}(\gamma_{d_i}\dots
  \gamma_{d_{\iota+1}})^{1/10} \cdot \epsilon_i\dots\epsilon_{\iota+1}(\gamma_{d_i}\dots
  \gamma_{d_{\iota+1}})^{9/10}  \norm{ v^\iota}_{D_{\beta_\iota}}
  \\&
  \leq \theta^{100 MN(i-\iota)} \cdot \epsilon_i\dots\epsilon_{\iota+1}(\gamma_{d_i}\dots
  \gamma_{d_{\iota+1}})^{9/10}  \norm{ v^\iota}_{D_{\beta_\iota}}.
  \end{align*}
Using the induction assumption at $\iota$, we get
\eqref{ToBeProvedRec} at $i$. This concludes the induction and the
proof of \eqref{ToBeProvedRec}.

From \eqref{ToBeProvedRec} at $j_0$, we obtain in particular
  \begin{equation}
  \norm{ v^{j_0}}_{L^2} \leq \theta^{100 MNj_0} \epsilon_1\dots\epsilon_{j_0}(\gamma_{d_1}\dots
  \gamma_{d_{j_0}})^{9/10} \norm{ v^0}_{D_{k_0}}.
  \end{equation}
As $\norm{v^0}_{D_{k_0}}\leq \norm{v^0}_{C^1}$, this concludes
the proof.
\end{proof}

\begin{lem}
\label{lem:ItereC1} There exists a constant $C$ (depending on
$M,A,\{\gamma_d\}, K)$ satisfying the following property. Let
$(k_0,k_1,\dots,k_j)$ be an admissible sequence. Let
$\psi_1,\dots,\psi_{j}$ be functions from $Y$ to $\C$, and let
$\epsilon_1,\dots,\epsilon_{j}$ belong to $[0,1]$. We assume
that $\norm{\psi_i}_{\CC^{A,3\epsilon}_{MN}}\leq
\epsilon_i\gamma_{d_i}$.

Let $v^0:Y\to \C$, define a sequence of functions $v^i$ by
induction, by $v^i=\Lp^{MN}_{k_i}(\psi_i v^{i-1})$. Then
  \begin{equation}
  \norm{v^{j}}_{C^1} \leq C (1+k_0^2)\left(\prod_{i=1}^{j} \epsilon_i \gamma_{d_i}^{1/3}\right)
  \theta^{30MNj} \norm{v^0}_{C^1}.
  \end{equation}
\end{lem}
\begin{proof}
We write $j_0=j/2$ or $(j-1)/2$, depending on whether $j$ is even or
odd.

Let $\pshi_i=e^{-i k_i S_{MN}^Y\phi_Y}\psi_i$, so that
$v^i=\Lp^{MN}(\pshi_i v^{i-1})$. We have $|\pshi_i(x)|\leq
\epsilon_i\gamma_{d_i}e^{3\epsilon r^{(MN)}(x)}$ and, for $h\in
\HH_{MN}$,
  \begin{align*}
  \norm{D(\pshi_i\circ h)(x)}&\leq  \norm{D (\psi_i\circ h)(x)}+|k_i|
  \norm{D(S^Y_{MN}\phi_Y\circ h)(x)} |\psi_i(hx)|
  \\&
  \leq C\epsilon_i \gamma_{d_i} e^{3\epsilon r^{(MN)}(hx)}
  + C |k_i| r^{(MN)}(hx) \epsilon_i \gamma_{d_i}
  e^{3\epsilon r^{(MN)}(hx)}
  \\&
  \leq
  C |k_i| \epsilon_i \gamma_{d_i} e^{4\epsilon r^{(MN)}(hx)}
  \end{align*}
for some constant $C\geq 1$ depending only on $M$ and $A$. Let
$B=C \max|k_i|$, this shows that
$\norm{\pshi_i}_{\CC^{B,4\epsilon}_{MN}}\leq
\epsilon_i\gamma_{d_i}$.

We can apply \eqref{eq:ContractC1Iteres} between the indexes $1$ and
$j_0$, to get
  \begin{align*}
  \norm{v^{j_0}}_{C^1} &
  \leq C (\max|k_i|) \left(\prod_{i=1}^{j_0} \epsilon_i\gamma_{d_i}\right)
  \left( \theta^{100MNj_0}\norm{v^0}_{C^1} + \theta^{-MNj_0}\norm{v^0}_{L^2}\right)
  \\&
  \leq C \theta^{-MNj_0} \left(\prod_{i=1}^{j_0} \epsilon_i\gamma_{d_i}\right) (\max |k_i|)
  \norm{v^0}_{C^1}.
  \end{align*}

Applying \eqref{eq:ContractC1Iteres} between the indexes $j_0+1$ and
$j$, we obtain
  \begin{equation*}
  \norm{v^j}_{C^1}\leq C (\max|k_i|) \left(\prod_{i=j_0+1}^{j} \epsilon_i\gamma_{d_i}\right)
  \left( \theta^{100MN(j-j_0)}\norm{v^{j_0}}_{C^1} +
  \theta^{-MN(j-j_0)}\norm{v^{j_0}}_{L^2}\right).
  \end{equation*}
We will use the bound on $\norm{v^{j_0}}_{C^1}$ given by the
previous equation, and the bound on $\norm{v^{j_0}}_{L^2}$ from
Lemma \ref{Lem:ContractionL2} (if $j_0=0$, this lemma does not
apply since the sequence $(k_0)$ is not necessarily strongly
admissible, but the estimate \eqref{eq:EstimeeL2} is trivial in
this case). We obtain:
  \begin{align*}
  \norm{v^{j}}_{C^1}&
  \leq C \left( \prod_{i=1}^{j} \epsilon_i \gamma_{d_i}\right) \theta^{40MNj}
  (\max |k_i|)^2 \norm{v^0}_{C^1}
  \\&\ \ \ \ \ \
  + C \left(\prod_{i=1}^{j_0} \epsilon_i
  \gamma_{d_i}^{9/10}\right) \left( \prod_{i=j_0+1}^{j}
  \epsilon_i \gamma_{d_i}\right) (\max |k_i|)\theta^{40MNj}\norm{v^0}_{C^1}
  \\&
  \leq C \theta^{40MNj} (\max |k_i|)^2 \left( \prod_{i=1}^{j}
  \epsilon_i \gamma_{d_i}^{9/10}\right) \norm{v^0}_{C^1}.
  \end{align*}
Assume first that $\max |k_i| \leq 2(|k_0|+jK)$. As $\theta^{40MNj}
j^2 \leq C \theta^{30MNj}$, we obtain the conclusion of the lemma
(by bounding directly $\left( \prod_{i=1}^{j}
\gamma_{d_i}\right)^{9/10}$ by $\left( \prod_{i=1}^{j}
\gamma_{d_i}\right)^{1/3}$).

Assume now that $\max |k_i| > 2 (|k_0|+jK)$. We have $|k_0|+\sum
|d_i| \geq \max |k_i|$. Denote by $J$ the set of indexes $\geq 1$
for which $|d_i|>K$. Then
  \begin{equation}
  \sum_{i\in J} |d_i| \geq \max |k_i| - |k_0|-jK \geq \max
  |k_i|/2.
  \end{equation}
By \eqref{Kgrand1}, $\gamma_d \leq 1/(1+|d|)^{60/17}$ for any $|d|>
K$. We get
  \begin{equation*}
  \left( \prod \gamma_{d_i}\right)^{17/30}
  \leq \left(\frac{1}{\prod_{i\in J} (1+|d_i|)^{60/17}}\right)^{17/30}
  \leq \left( \frac{1}{\sum_{i\in J} |d_i|} \right)^2
  \leq 4/(\max |k_i|)^2.
  \end{equation*}
Finally,
  \begin{equation*}
  (\max |k_i|)^2 \left( \prod_{i=1}^{j}
  \gamma_{d_i}\right)^{9/10}
  = (\max |k_i|)^2 \left( \prod_{i=1}^{j}
  \gamma_{d_i}\right)^{17/30}\cdot \left( \prod_{i=1}^{j}
  \gamma_{d_i}\right)^{1/3}
  \leq 4  \left( \prod_{i=1}^{j}
  \gamma_{d_i}\right)^{1/3}.
  \end{equation*}
This yields again the conclusion of the lemma.
\end{proof}

\section{Proof of the local limit theorem}
\label{sec:DemontreLocal}

We fix a $C^6$ function $\psi : X\times \Sbb^1 \to \R$ with
vanishing average, and a real number $t_0>0$. We will study the
operators $\hat \TT_{t}:=\hat \TT( e^{it\psi}\cdot) $ for
$|t|\leq t_0$. We will first choose $M$, $A$, a sequence
$\gamma_d$ and an integer $K$ so that the results of Paragraph
\ref{sec:EstimeExcursion} apply. All these choices will depend
on $\psi$ and $t_0$.

\subsection{Choosing the constants}
Let $\psi_Y$ be the function defined in \eqref{DefinitPsiY}.
There exists a constant $C(\psi)$ such that $|S^Y_n
\psi_Y(x,\omega)| \leq C(\psi) r^{(n)}(x)$. More generally, as
$\TT$ is an isometry in the fiber direction $\Sbb^1$, we even
have
  \begin{equation}
  \left| \frac{\partial^4}{\partial \omega^4} S^Y_n\psi_Y(x,\omega)\right|
  \leq C(\psi) r^{(n)}(x).
  \end{equation}
In particular, for any $|t|\leq t_0$,
  \begin{equation}
  \left| \frac{\partial^4}{\partial \omega^4} e^{it
  S^Y_n\psi_Y(x,\omega)}\right| \leq C(t_0,\psi)r^{(n)}(x)^4.
  \end{equation}
Let us denote by $F^{(n,t)}_d$ the $d$-th Fourier coefficient
of $e^{it S^Y_n\psi_Y}$ in the circle direction. Making 4
integrations by parts in the circle direction and using the
previous equation yields
  \begin{equation}
  \label{psinpetit}
  |F^{(n,t)}_d(x)| \leq \frac{C(t_0,\psi)r^{(n)}(x)^4}{1+|d|^4}
  \leq \frac{C'(t_0,\psi) e^{\epsilon r^{(n)}(x)}}{1+|d|^4}.
  \end{equation}
There also exists $C(n,t_0,\psi)$ such that, for any $h\in \HH_n$,
  \begin{equation}
  \label{DeriveePasTropGrande}
  \norm{ D (F^{(n,t)}_d\circ h)(x)} \leq C(n,t_0,\psi)
  \frac{e^{\epsilon r^{(n)}(hx)}}{1+|d|^4}.
  \end{equation}

We fix once and for all an integer $M$ such that
  \begin{equation}
  \label{eq:Mgrand}
  \theta^{20MN} \sum_{d\in \Z} \min\left(1, \frac{C'(t_0,\psi)}
  {1+|d|^4}\right)^{1/3} <
  \theta^{10MN}
  \end{equation}
and
  \begin{equation}
  \label{eq:Mgrand2}
  \theta^{100MN} \sum_{d\in \Z} \min\left(1, \frac{C'(t_0,\psi)}
  {1+|d|^4}\right) < 1/4.
  \end{equation}

Let $\gamma_d= \min\left(1, \frac{C'(t_0,\psi)} {1+|d|^4}\right)$.
By \eqref{DeriveePasTropGrande}, we can then choose a constant $A$
such that
  \begin{equation}
  \norm{F^{(MN,t)}_d}_{\CC^{A,\epsilon}_{MN}} \leq \gamma_d
  \end{equation}
for any $d\in \Z$. Finally, we choose $K$ satisfying
\eqref{Kgrand1}--\eqref{Kgrand3}.

All the constants $C$ we will consider until the end of this section
may depend on $M,A,\{\gamma_d\}, K$. We will work on the space
$X^{(MN)}$, with the map $U=U^{(MN)}$, to prove Theorem
\ref{thm:EstimeesValeurPropreFinal} for $t\in [-t_0,t_0]$. We will
freely use all the results that we proved in Section
\ref{sec:MelangeExp}. Formally, we proved these results for
$X^{(N)}$, but the same arguments hold verbatim in $X^{(MN)}$.

\emph{As in the proof of Theorem \ref{thm:MelangeExponentiel}, we
will assume until the end of the proof that $d^{(MN)}=1$, i.e.,
$U^{(MN)}$ is mixing. Only at the end of the proof will we give the
modifications to be done to handle the general case.}

\subsection{The renewal process}
\label{subsec:DefinitBarQ}

As in Paragraph \ref{subsec:Strategie}, let us define a space
$\B_K=\bigoplus_{|k|\leq K} C^1(Y)$, endowed with the norm of the
supremum of the $C^1$ norms of the different components. We will see
an element $v$ of $\B_K$ as a set of functions $(v_k)_{|k|\leq K}$
where $v_k$ corresponds to frequency $k$, and then
$\norm{v}_{\B_K}=\sup_{|k|\leq K} \norm{v_k}_{C^1}$. We will also
write $\norm{v}_{C^0}=\sup \norm{v_k}_{C^0}$.

For $z\in \C$, $t\in [-t_0,t_0]$ and $\kk=(k_0,\dots,k_j)$ an
admissible sequence, we formally define an operator $Q_{\kk}^t(z)$
on $C^1(Y)$, by
  \begin{equation}
  Q_{\kk}^t(z) v= \Lp^{MN}_{k_j}( z^{r^{(MN)}} F^{(MN,t)}_{d_j}
  \Lp^{MN}_{k_{j-1}} z^{r^{(MN)}}\dots \Lp^{MN}_{k_1}( z^{r^{(MN)}}
  F^{(MN,t)}_{d_1} v)\dots).
  \end{equation}
Intuitively, this operator applies to a function of frequency $k_0$,
and gives a function of frequency $k_j$. If $\BB$ is a Banach space
of functions from $Y \times \Z$ to $\C$, it is therefore more
natural to consider an operator $\bar Q_{\kk}^t(z)$ from $\B$ to
$\B$, defined by $( \bar Q_{\kk}^t(z)v)_k=0$ if $k\not=k_j$, and $(
\bar Q_{\kk}^t(z)v)_{k_j}=Q_{\kk}^t(z) v_{k_0}$. This applies for
instance if $\BB=\BB_K$ (and $|k_0|\leq K$, $|k_j|\leq K$). We will
occasionally use the operators $\bar Q^t_{\kk}(z)$, but the
technical estimates will be formulated in terms of $Q_{\kk}^t(z)$.

\begin{lem}
\label{lem:ControleQkt} The operator $Q_{\kk}^t(z)$ acts
continuously on $C^1(Y)$ for any $t\in[-t_0,t_0]$ and any $|z|\leq
e^{2\epsilon}$, and its norm is bounded by $C(1+k_0^2)\theta^{20MNj}
\prod_{i=1}^j \gamma_{d_i}^{1/3}$. Moreover, the map $z\mapsto
Q_{\kk}^t(z)$ is holomorphic from $\{|z|< e^{2\epsilon}\}$ to
$\End(C^1(Y))$ the set of continuous linear operators on $C^1(Y)$.

There exist $a>0$ and $C>0$ such that, for all $|t-t'|\leq a$, for
any admissible sequence $\kk$,
  \begin{equation}
  \label{eq:ControlePerturb}
  \norm{Q_{\kk}^t(z) -Q_{\kk}^{t'}(z)}_{\End(C^1(Y))} \leq
  C |t-t'| (1+k_0^2)\theta^{20MNj} \prod_{i=1}^j
  \gamma_{d_i}^{1/3}.
  \end{equation}
Finally, if $|t|\leq a$,
  \begin{equation}
  \label{eq:ControlePerturb0}
  \norm{Q_{\kk}^t(z)} \leq C(1+k_0^2) (C|t|)^{ \#\{i\st d_i\not=0\}}
  \theta^{20MNj} \prod_{i=1}^j \gamma_{d_i}^{1/3}.
  \end{equation}
\end{lem}
\begin{proof}
To estimate the norm of $Q_{\kk}^t(z)$, we use the estimate
given by Lemma \ref{lem:ItereC1}, taking $\epsilon_i=1$ and
$\psi_i=z^{r^{(MN)}} F^{(MN,t)}_{d_i}$. If $|z|\leq
e^{2\epsilon}$, we have
$\norm{\psi_i}_{\CC^{A,3\epsilon}_{MN}}\leq
\norm{F^{(MN,t)}_{d_i}}_{\CC^{A,\epsilon}_{MN}} \leq
\gamma_{d_i}$. We obtain
  \begin{equation}
  \norm{Q^t_{\kk}(z)}_{\End(C^1(Y))}
  \leq C (1+k_0^2)\left(\prod_{i=1}^j \gamma_{d_i} ^{1/3} \right)
  \theta^{30 MNj}.
  \end{equation}
If $|z|<e^{2\epsilon}$, each function $\psi_i 1_{r^{(MN)}>n}$
tends to $0$ in $\CC^{A,3\epsilon}_{MN}$ when $n$ tends to
infinity. As a consequence, $z\mapsto Q_{\kk}^t(z)$ is a
uniform limit of polynomials on any compact subset of $\{|z|<
e^{2\epsilon}\}$, and is therefore holomorphic there.

To prove the rest of the lemma, we will use the following inequality
(which can easily be proved by $4$ integrations by parts): there
exists $C>0$ such that, for any $t,t'\in [-t_0,t_0]$ and for any
$d\in \Z$,
  \begin{equation}
  \label{eq:ttprimeproche}
  \norm{ F^{(MN,t)}_d - F^{(MN,t')}_d}_{\CC^{A,
  \epsilon}_{MN}} \leq C|t-t'|\gamma_d.
  \end{equation}

To prove \eqref{eq:ControlePerturb}, let us write $Q_{\kk}^t(z)v
-Q_{\kk}^{t'}(z)v$ as
  \begin{multline*}
  \sum_{b=0}^j \Lp^{MN}_{k_j}(z^{r^{(MN)}} F^{(MN,t)}_{d_j} \Lp^{MN}_{k_{j-1}} \dots
  \Lp^{MN}_{k_b}(
  z^{r^{(MN)}} (F^{(MN,t)}_{d_b}-F^{(MN,t')}_{d_b})\Lp^{MN}_{k_{b-1}}
  (
  \\ z^{r^{(MN)}}F^{(MN,t')}_{d_{b-1}}
  \Lp^{MN}_{k_{b-2}}(\dots \Lp^{MN}_{k_1}(z^{r^{(MN)}}F^{(MN,t')}_{d_1}
  v)\dots).
  \end{multline*}
Fix $b$. To estimate the corresponding term in this equation,
we will again use Lemma \ref{lem:ItereC1}. Let
$\psi_i=z^{r^{(MN)}}F^{(MN,t)}_{d_i}$ for $i>b$,
$\psi_i=z^{r^{(MN)}}F^{(MN,t')}_{d_i}$ for $i<b$ and
$\psi_b=z^{r^{(MN)}}(F^{(MN,t)}_{d_b}-F^{(MN,t')}_{d_b})$. Let
also $\epsilon_i=1$ for $i\not=b$. Then $\psi_i,\epsilon_i$
satisfy the assumptions of Lemma \ref{lem:ItereC1} for
$i\not=b$. Let finally $\epsilon_b=C|t'-t|$ (where $C$ is as in
\eqref{eq:ttprimeproche}). If $t'$ is close enough to $t$, we
have $\epsilon_b\leq 1$, and the assumptions of Lemma
\ref{lem:ItereC1} are again satisfied by
\eqref{eq:ttprimeproche}.

Using this lemma, we obtain (after summation over $b$)
  \begin{equation}
  \norm{Q_{\kk}^t(z)v
  -Q_{\kk}^{t'}(z)v}_{C^1} \leq C(j+1)|t'-t| (1+k_0^2)\left(\prod_{i=1}^j
  \gamma_{d_i}^{1/3}\right) \theta^{30 MNj}\norm{v_{k_0}}_{C^1}.
  \end{equation}
As $(j+1) \theta^{30MNj} \leq C \theta^{20MNj}$, we get
\eqref{eq:ControlePerturb}.

Finally, to prove \eqref{eq:ControlePerturb0}, note that
$F^{(MN,0)}_{d}=0$ if $d\not=0$. As a consequence,
\eqref{eq:ttprimeproche} applied to $t'=0$ gives $\norm{
F^{(MN,t)}_d }_{\CC^{A, \epsilon}_{MN}} \leq C|t|\gamma_d$. We can
therefore apply Lemma \ref{lem:ItereC1} to $\epsilon_i=1$ if
$d_i=0$, and $\epsilon_i=C|t|$ if $d_i\not=0$, to obtain
\eqref{eq:ControlePerturb0}.
\end{proof}

Let us then define formally an operator $R(z,t)$ on $\B_K$ by
$R(z,t)=\sum \bar Q_{\kk}^t(z)$, where we sum over all admissible
sequences $\kk$ with $|k_0|\leq K$ and $|k_j|\leq K$, i.e.,
  \begin{equation}
  (R(z,t)v)_k = \sum_{j=1}^\infty \sum_{\substack{k_0,k_1,\dots,k_{j-1}\\ |k_0|\leq K\\
  \kk=(k_0,k_1,\dots,k_{j-1},k) \text{ admissible}}}
  Q_{\kk}^t(z) v_{k_0}.
  \end{equation}
The coefficient of $z^n$ corresponds to considering the first
returns to $Y\times[-K,K]$ after a time exactly $n$. By
\eqref{eq:RedonneKtnexc}, this is exactly the operator $R^t_n$
defined in \eqref{eq:DefinitRtn}. Using the estimates in Lemma
\ref{lem:ControleQkt}, our next goal is to prove that the
operators $R^t_n$ satisfy the assumptions of Theorem
\ref{thm:renouvellement}. Indeed, this theorem will thus
provide us with a good estimate for $T^t_n$ (defined in
\eqref{eq:DefinitTtn}), which is the main building block of
$\hat\UU_t^n$.

\begin{lem}
The formal series $R(z,t)$ defines an holomorphic function on
the disk $|z|< e^{2\epsilon}$, uniformly bounded in $t\in
[-t_0,t_0]$. In particular, there exists $C>0$ such that, for
any $t\in [-t_0,t_0]$, for any $n\in \N$, for any $v\in \B_K$,
$\norm{R_n^tv}_{\B_K} \leq C e^{-n\epsilon} \norm{v}_{\B_K}$.

Moreover,
  \begin{equation}
  \norm{R(z,t)v -R(z,t')v}_{\B_K} \leq C|t-t'| \norm{v}_{\B_K}.
  \end{equation}
In particular, for any $n\in \N$, for any $v\in \B_K$,
$\norm{R_n^tv-R_n^{t'}v }_{\B_K} \leq C |t-t'| e^{-n\epsilon}
\norm{v}_{\B_K}$.
\end{lem}
\begin{proof}
As $\theta^{20MN} \sum_{d\in\Z}\gamma_d^{1/3}<1$, the estimates
given by Lemma \ref{lem:ControleQkt} are summable. This directly
implies the lemma.
\end{proof}

\begin{lem}
There exists a constant $C$ such that, for any $z$ with
$|z|\leq e^{2\epsilon}$, for any $t\in [-t_0,t_0]$, for any
$v\in \B_K$,
  \begin{equation}
  \norm{R(z,t) v}_{\B_K} \leq \frac{1}{2}\norm{v}_{\B_K} +
  C \norm{v}_{C^0}.
  \end{equation}
\end{lem}
\begin{proof}
Fix an integer $P$. We define a truncated series $R(z,t,P)$ by
summing as in $R(z,t)$ along admissible sequences
$\kk=(k_0,k_1,\dots,k_j)$, but with the additional restrictions
$\sup |k_i|\leq P$ and $j\leq P$. When $P$ tends to infinity,
$R(z,t,P)$ converges (in norm) to $R(z,t)$, uniformly for $(z,t)\in
\{|z|\leq e^{2\epsilon}\}\times [-t_0,t_0]$. We will show that, for
any $P\in \N$, there exists $C(P)$ such that
  \begin{equation}
  \label{eq:PetitBorne}
  \norm{R(z,t,P)v}_{\B_K} \leq
  \frac{1}{3}\norm{v}_{\B_K}+ C(P) \norm{v}_{C^0}.
  \end{equation}
This implies the desired result, by choosing a large enough $P$.

Let $\kk$ be an admissible sequence of length $j>0$. Iterating
$j$ times the equation \eqref{eq:ContracteTriviale} (applied to
the functions $\psi_i=z^{r^{(MN)}} e^{-i k_i S^Y_{MN}\phi_Y}
F^{(MN,t)}_{d_i}$), we obtain a constant $C(\kk)$ such that,
for any $v\in C^1(Y)$,
  \begin{equation}
  \norm{ Q^t_{\kk}(z) v}_{C^1} \leq \theta^{100MN j}
  \left(\prod_{i=1}^j \gamma_{d_i}\right) \norm{v}_{C^1} +
  C(\kk) \norm{v}_{C^0}.
  \end{equation}
The operator $R(z,t,P)$ involves only a finite number of
admissible sequences. Denoting by $C(P)$ the sum of $C(\kk)$
over these admissible sequences, we obtain for any $v\in \BB_K$
  \begin{align*}
  \norm{ R(z,t,P) v}_{\B_K} &\leq \sum_{j=1}^P \theta^{100MNj}\left(\sum_{d\in
  \Z} \gamma_d\right)^j \norm{v}_{\B_K} + C(P) \norm{v}_{C^0}
  \\&
  \leq \frac{ \theta^{100MN} \sum \gamma_d}{1-\theta^{100MN}
  \sum \gamma_d} \norm{v}_{\B_K} + C(P) \norm{v}_{C^0}
  \leq \frac{1}{3}\norm{v}_{\B_K} + C(P) \norm{v}_{C^0},
  \end{align*}
by \eqref{eq:Mgrand2}.
\end{proof}

\begin{cor}
\label{cor:TrouSpectral} For any $t\in [-t_0,t_0]$ and for any
$|z|\leq e^{2\epsilon}$, the operator $R(z,t)$ acting on $\B_K$
has an essential spectral radius bounded by $1/2$.
\end{cor}
\begin{proof}
This is a consequence of Hennion's Theorem \cite{hennion} (or more
precisely of the version without iteration of this theorem given in
\cite[Lemma 2.2]{BGK:coupling}, since the operator $R(z,t)$ is
\emph{a priori} not continuous for the $C^0$ norm).
\end{proof}

\begin{definition}
Let $\psi:X\times\Sbb^1 \to\R$ be a $C^6$ function. We say that it
is continuously periodic if there exist $a>0$, $\lambda>0$ and
$f:X\times \Sbb^1 \to \R/\lambda \Z$ measurable such that
$\psi=f-f\circ \TT+a \mod \lambda$ almost everywhere, and $f$ is
continuous on $Y\times \Sbb^1$. Otherwise, we say that $\psi$ is
continuously aperiodic.
\end{definition}
Proposition \ref{prop:CaracteriseAperiodique} says that
aperiodicity and continuous aperiodicity are equivalent.
However, we will be able to prove this equivalence only at the
complete end of our arguments. Until then, it will be more
convenient to work with the notion of continuous aperiodicity.

\begin{prop}
\label{prop:ConstruitFonctionPropre} For any $z\in\overline{\D}
\moins \{1\}$, the operator $I-R(z,0)$ is invertible on $\B_K$.
Moreover, if the function $\psi$ is continuously aperiodic, the
operator $I-R(z,t)$ is invertible on $\B_K$ for any $(z,t)\in
(\overline{\D}\times [-t_0,t_0])\moins \{(1,0)\}$.
\end{prop}
\begin{proof}
Let $|z|\leq 1$ and $t\in[-t_0,t_0]$. If the operator $I-R(z,t)$ is
not invertible, its kernel contains a nonzero function
$v=(v_{-K},\dots,v_{K})$ by Corollary \ref{cor:TrouSpectral}. Let us
define a function $v_k$, for $|k|>K$, by
  \begin{equation*}
  v_k=\sum_{p=1}^\infty \sum_{\substack{\kk=(k_0,k_1,\dots,k_{j-1},k) \text{
  admissible}\\ |k_0|\leq K}}Q_{\kk}^t(z) v_{k_0}.
  \end{equation*}
Lemma \ref{lem:ControleQkt} implies (after summation over the
admissible sequences) that $\sum_{k\in \Z}\norm{v_k}_{C^1}<\infty$.
Moreover, for any $k\in \Z$,
  \begin{equation}
  \label{eq:Invariancevk}
  v_k = \sum_{l\in \Z} \Lp_k^{MN}(z^{r^{(MN)}}
  F^{(MN,t)}_{k-l} v_l).
  \end{equation}
This equation is indeed a consequence of the construction of the
$v_k$'s if $|k|>K$, and of the fact that $v$ is a fixed point of
$R(z,t)$ if $|k|\leq K$.

Let us define a continuous function $g$ on $Y \times \Sbb^1$ by
$g(x,\omega) = \sum_{k\in \Z} v_k(x) e^{ik\omega}$. As $v$ is
nonzero, $g$ is also nonzero. The invariance equation
\eqref{eq:Invariancevk} translates into the following for $g$:
  \begin{equation}
  \hat\UU_Y( z^{r^{(MN)}} e^{it S^Y_{MN}\psi_Y} g)=g,
  \end{equation}
where $\hat\UU_Y$ is the transfer operator associated to the
map which is induced by $\UU=\UU^{(MN)}$ on $Y$. Lemma
\ref{lem:DecritFctPropre} yields $|z|=1$ and $g\circ \UU_Y=
e^{it S^Y_{MN}\psi_Y} z^{r^{(MN)}}g$. Let us extend $g$ to the
whole space $X^{(MN)}\times \Sbb^1$ by setting
  \begin{equation}
  g(x,i,\omega)=z^{i} g(x,0,\omega)\exp\left(it
  \sum_{j=0}^{i-1}\psi\circ\UU^j(x,\omega)\right).
  \end{equation}
This function is bounded (since $g$ is bounded on $Y$), nonzero, and
satisfies $g\circ \UU= z e^{it\psi} g$.

If $t=0$, we obtain $g\circ \UU= z g$. But the map $\UU$ is mixing
(this was proved in Theorem \ref{thm:MelangeExponentielPrecis} and
in \eqref{eq:MelangeEnHaut} for $\UU^{(N)}$, the same proof holds
for $\UU^{(MN)}$). As a consequence, $z=1$.

If $t\not=0$, let $f:X^{(MN)}\times\Sbb^1 \to \R/2\pi\Z$ be the
logarithm of $g$, and let $a$ be such that $z=e^{-ia}$. Then
$t\psi\circ \tilde \pi^{(MN)}=f\circ \UU-f+a \mod 2\pi$, and
$f$ is continuous on $Y\times \Sbb^1 \subset X^{(MN)}\times
\Sbb^1$ (we have reintroduced the projection $\tilde\pi^{(MN)}$
in the notations since we will soon be confronted to lifting
problems). In general, $f$ is not constant on the fibers of
$\tilde\pi^{(MN)}$, and can therefore not be written as $\tilde
f \circ \tilde\pi^{(MN)}$ in $\R/2\pi\Z$. However, since the
fibers of $\tilde\pi^{(MN)}$ are countable, \cite[Theorem
1.4]{gouezel:local} shows that there exist $\lambda$ of the
form $2\pi/n$ for some integer $n$, and $\tilde f : X\times
\Sbb^1 \to \R/\lambda \Z$, such that $f=\tilde f \circ
\tilde\pi^{(MN)} \mod \lambda$ almost everywhere. As a
consequence, $t\psi=\tilde f\circ \TT -\tilde f+a \mod
\lambda$, and $\tilde f$ has a continuous version on $Y\times
\Sbb^1$ (since this is the case for $f$). Hence, $\psi$ is
continuously periodic.
\end{proof}

\begin{lem}
\label{lem:DonneProjecteur} The operator $R(1,0)$ has a simple
eigenvalue at $1$. The corresponding spectral projection is given by
$(Pv)_0=\int_Y v_0 \dd\mu_Y$, and $(Pv)_k=0$ if $k\not=0$. Denoting
by $R'(z,t)$ the derivative with respect to $z$ of $R(z,t)$, we have
$P R'(1,0)P=\mu^{(MN)}(X^{(MN)}) P$.
\end{lem}
\begin{proof}
We have $(R(1,0) v)_k=\Lp^{MN}_{k} v_k$, it is therefore sufficient
to know the spectral properties of the operators $\Lp^{MN}_k$ (for
$|k|\leq K)$ to conclude. For $k\not=0$, there operators have a
spectral radius $<1$, while for $k=0$ there is a simple eigenvalue
at $1$, the corresponding eigenprojection being given by integration
(as we saw in the proofs of Lemma \ref{lem:TnkPourknot0} and
Corollary \ref{cor:ExprimeTn0}). This yields the desired formula for
$P$.

As $PR^0_{j}P=\mu_Y\{ r^{(MN)}=j\}P$ for $j\geq 1$, we have
  \begin{equation}
  PR'(1,0)P= \sum j\mu_Y\{r^{(MN)}=j\} P =\mu^{(MN)}(X^{(MN)}) P,
  \end{equation}
by Kac's Formula.
\end{proof}

\subsection{Estimate of the perturbed eigenvalue}

In this paragraph, we prove the following estimate (which is
necessary to apply Theorem \ref{thm:renouvellement}).
\begin{thm}
\label{thm:Controlelambdat} Denote by $\lambda(1,t)$ the eigenvalue
close to $1$ of $R(1,t)$, for small $t$. Then
  \begin{equation}
  \label{eq:ExprimeLambda}
  \lambda(1,t)=1-\mu^{(MN)}(X^{(MN)})\frac{ \sigma^2 t^2}{2} +O(t^3),
  \end{equation}
where $\sigma^2$ is given by \eqref{eq:DefinitSigma2}.
\end{thm}

The proof will take the rest of this paragraph. We will write $R(t)$
and $\lambda(t)$ instead of $R(1,t)$ and $\lambda(1,t)$, since we
will only consider $z=1$.

Let $f^t$ be the eigenfunction (in $\B_K$) of $R(t)$ for the
eigenvalue $\lambda(t)$, normalized so that $\int f^t_0=1$
(this is possible since $\int f_0^0=1$ and $f^t$ converges to
$f^0$ in $\B_K$). Note that $f^t=f^0+O(t)$ and
$\lambda(t)=1+O(t)$ (since $R(t)=R(0)+O(t)$ and the simple
isolated eigenvalues, as well as the corresponding
eigenfunctions, depend in a Lipschitz way on the operator).
Moreover, $f^0_0=1$, and $f^0_k=0$ for $k\not=0$.

\begin{lem}
We have $\lambda(t)=1+O(t^2)$.
\end{lem}
\begin{proof}
We have $(R(t)f^t)_0=\sum Q_{\kk}^t(1) f^t_{k_0}$ where the
summation is over the admissible sequences $\kk=(k_0,\dots,k_j)$
with $|k_0|\leq K$ and $k_j=0$. If $j\geq 2$, there are at least two
nonzero differences $d_i=k_i-k_{i-1}$, and the sum of the
corresponding terms is therefore bounded by $Ct^2$, by
\eqref{eq:ControlePerturb0}. If $j=1$ but $k_0\not=0$, the
difference is nonzero, which gives a $O(t)$ factor. As
$f^t_{k_0}=O(t)$, the resulting term is therefore also $O(t^2)$. It
remains $(R(t)f^t)_0=Q^t_{(0,0)}(1)f^t_0+O(t^2)$. As
$R(t)f^t=\lambda(t)f^t$ and $\int f^t_0=1$, we obtain after
integration
  \begin{align*}
  \lambda(t)&=\int_Y Q^t_{(0,0)}(1)f^t_0 + O(t^2)
  =\int_Y \Lp^{MN}( F^{(MN,t)}_0 f^t_0) + O(t^2)
  \\&
  =\int_{Y\times \Sbb^1} e^{itS_{MN}^Y\psi_Y(x,\omega)}f^t_0(x)
  + O(t^2).
  \end{align*}
As $\int f^t_0=1$, we get
  \begin{equation}
  \lambda(t)=1+ \int (e^{itS_{MN}^Y\psi_Y}-1) (f^t_0-1) + \int
  (e^{itS_{MN}^Y\psi_Y}-1) + O(t^2).
  \end{equation}
Since $f^t_0 = f^0_0+O(t)=1+O(t)$, the first integral is $O(t^2)$.
For the second one,
  \begin{equation}
  \int (e^{itS_{MN}^Y\psi_Y}-1)=
  it\int S_{MN}^Y\psi_Y+O(t^2)=MNit\int_{X\times \Sbb^1} \psi +
  O(t^2)  =O(t^2)
  \end{equation}
since $\int \psi=0$. This finally yields $\lambda(t)=1+O(t^2)$.
\end{proof}

Define a function $g_k$ on $Y$ by $g_k(x)=\int  S_{MN}^Y
\psi_Y(x,\omega)e^{-ik\omega} \dd\omega$.
\begin{lem}
\label{lem:GoodEstimateFd}
The function $g_k$ belongs to
$\CC^{1,\epsilon}_{MN}$. Moreover, there exists a constant
$C>0$ such that, for any small enough $t$ and for any $k\in\Z$,
  \begin{equation}
  \label{eq:EstimeFMNtprecis}
  \norm{F^{(MN,t)}_k- 1_{k=0}-it
  g_k}_{\CC^{1,\epsilon}_{MN}}\leq \frac{Ct^2}{1+k^4}.
  \end{equation}
\end{lem}
\begin{proof}
Write
  \begin{align*}
  F^{(MN,t)}_k(x)- 1_{k=0}-it
  g_k(x)&=\int_{\Sbb^1} \left(e^{itS_{MN}^Y\psi_Y(x,\omega)}
  -1-itS_{MN}^Y\psi_Y(x,\omega)\right)e^{-ik\omega}\dd\omega
  \\
  &=-t^2 \int_{v=0}^1 (1-v) \left(\int_{\Sbb^1} S_{MN}^Y
  \psi_Y(x,\omega)^2
  e^{itS_{MN}^Y\psi_Y(x,\omega)v}e^{-ik\omega}\dd\omega\right)
  \dd v.
  \end{align*}
This gives \eqref{eq:EstimeFMNtprecis} after $4$ integrations
by parts with respect to $\omega$.
\end{proof}

\begin{lem}
\label{lem:Decritftk} For any $|k|\leq K$, we have in $C^1(Y)$
  \begin{equation}
  f^t_k=f^0_k + it \sum_{n=1}^\infty \Lp_k^{MNn}(g_k)+O(t^2).
  \end{equation}
\end{lem}
Note that $g_k$ belongs to $\CC^{1,\epsilon}_{MN}$, which implies
that $\Lp_k^{MN} g_k\in C^1(Y)$ by Theorem
\ref{thm:MainContraction}. The series $\sum_{n\in \N} \Lp_k^{MNn}
\Lp_k^{MN} g$ is therefore convergent in $C^1(Y)$: for $k\not=0$,
the spectral radius of $\Lp_k^{MN}$ on $C^1(Y)$ is $<1$ and the
convergence is trivial. For $k=0$, there is still exponential
convergence for functions with zero average, which is the case of
$g_0$ because $\int\psi=0$.

\begin{proof}[Proof of Lemma \ref{lem:Decritftk}]
As $\lambda(t)=1+O(t^2)$, we have
  \begin{align*}
  \frac{f^t-f^0}{t}&=\frac{\lambda(t)f^t-f^0}{t}+O(t)
  =\frac{R(t)f^t-R(0)f^0}{t}+O(t)
  \\&
  =(R(t)-R(0))\frac{f^t-f^0}{t} + R(0) \frac{f^t-f^0}{t}+
  \frac{R(t)-R(0)}{t}f^0 +O(t).
  \end{align*}
Since $R(t)-R(0)=O(t)$ and $f^t-f^0=O(t)$, we obtain
  \begin{equation}
  \label{eq:mlqjsdfkljqsdmf}
  (I-R(0)) \frac{f^t-f^0}{t}=\frac{R(t)-R(0)}{t}f^0+O(t).
  \end{equation}
The operator $R(0)$ simply acts by $(R(0)v)_k=\Lp_k^{MN}v_k$. Let us
study $(R(t)f^0)_k=\sum_{\kk}Q_{\kk}^t(1) 1$, where $\kk$ is an
admissible sequence beginning by $0$ and ending by $k$. If the
length of this admissible sequence is at least $2$, there are two
nonzero differences, and we obtain a term bounded by $O(t^2)$.
Hence,
  \begin{equation}
  (R(t)f^0)_k=Q_{(0,k)}^t(1)1+O(t^2)=\Lp_k^{MN}(
  F^{(MN,t)}_k)+O(t^2).
  \end{equation}
Applying Lemma \ref{lem:GoodEstimateFd} and using the fact that
$\Lp_k^{MN}$ is continuous from $\CC^{1,\epsilon}_{MN}$ to
$C^1(Y)$, we get in $C^1(Y)$
  \begin{equation}
  \label{eq:qlsdkfjlk}
  (R(t)f^0)_k= 1_{k=0}+it \Lp_k^{MN}g_k+O(t^2)
  = (R(0)f^0)_k + it\Lp_k^{MN}g_k+O(t^2).
  \end{equation}
Let $h_k=\sum_{n>0} \Lp_k^{MNn}g_k$. Denote by $h$ the corresponding
element in $\B_K$, so that the $k$-th component of $(I-R(0))h$ is
equal to $\Lp_k^{MN}g_k$. The equations \eqref{eq:mlqjsdfkljqsdmf}
and \eqref{eq:qlsdkfjlk} imply that
  \begin{equation}
  (I-R(0)) \left(\frac{f^t-f^0}{t}-ih\right)=O(t).
  \end{equation}
As $I-R(0)$ is invertible on the set of elements $v$ of $\B_K$ with
$\int v_0=0$, this shows that $(f^t-f^0)/t -ih=O(t)$, which is the
desired conclusion.
\end{proof}

Let $\UU_Y$ be the map induced by $\UU=\UU^{(MN)}$ on $Y\times
\Sbb^1$. The associated transfer operator $\hat\UU_Y$ acts on
each frequency $k$ by $\Lp_k^{MN}$. From the spectral
properties of the operators $\Lp_k^{MN}$, we obtain the
convergence of the series
  \begin{equation}
  \begin{split}
  \tilde\sigma^2&=\int_Y (S_{MN}^Y\psi_Y)^2+2\sum_{n=1}^\infty
  \int_Y S_{MN}^Y\psi_Y\cdot S_{MN}^Y\psi_Y\circ \UU_Y^n
  \\&
  = \int_Y (S_{MN}^Y\psi_Y)^2 +2\sum_{n=1}^\infty
  \int_Y \hat\UU_Y^n S_{MN}^Y\psi_Y\cdot S_{MN}^Y\psi_Y.
  \end{split}
  \end{equation}

\begin{lem}
\label{eq:DonneLambdattilde} We have
$\lambda(t)=1-\tilde\sigma^2 t^2/2+O(t^3)$.
\end{lem}
\begin{proof}
Let us estimate $(R(t)f^t)_0$. We have
  \begin{equation*}
  (R(t)f^t)_0=\sum_{1\leq |k|\leq K}
  \sum_{\kk=(k,k_1,\dots,k_{j-1},0) \text{ admissible}}
  Q_{\kk}^t(1)f^t_k + \sum_{\kk=(0,k_1,\dots,k_{j-1},0) \text{ admissible}}
  Q_{\kk}^t(1)f^t_0.
  \end{equation*}
In the first sum, $f^t_k=O(t)$. If there are two nonzero
differences in the admissible sequence $\kk$, we therefore
obtain terms bounded by $O(t^3)$ by
\eqref{eq:ControlePerturb0}. In the second sum, we also get
$O(t^3)$ unless there are at most two nonzero differences,
which is possible only for the sequences $\kk=(0,0)$ and
$\kk=(0,\ell,\dots,\ell,0)$, where $\ell$ is repeated a number
of times, say $j$, and $|\ell|>K$. Hence,
  \begin{equation*}
  (R(t)f^t)_0=\sum_{1\leq |k|\leq K} \Lp^{MN}( F^{(MN,t)}_{-k}
  f^t_k) + \Lp^{MN}(F^{(MN,t)}_0 f^t_0) + \sum
  Q_{(0,\ell,\dots, \ell,0)}^t(1)f^t_0 + O(t^3).
  \end{equation*}
We have
  \begin{equation}
  Q_{(0,\ell,\dots,\ell,0)}^t(1) v= \Lp^{MN}(
  F^{(MN,t)}_{-\ell} \Lp_{\ell}^{MN} F^{(MN,t)}_{0}
  \Lp_{\ell}^{MN}\dots
  \Lp_{\ell}^{MN}(F^{(MN,t)}_{\ell}f^t_0)\dots).
  \end{equation}
As there are two nonzero differences in these admissible sequences,
the contribution of these terms to $R(t)f^t_0$ is $O(t^2)$.
Moreover, $F^{(MN,t)}_0=1+O(t)$. If we replace $F^{(MN,t)}_0$ by
$1$, we get an additional error of $O(t)$ in each term. It can be
checked as in the proof of \eqref{eq:ControlePerturb} that these
errors are summable. In the same way, $f^t_0$ may be replaced by $1$
since the error is $O(t)$. We get
  \begin{multline*}
  (R(t)f^t)_0=\sum_{1\leq |k|\leq K} \Lp^{MN}( F^{(MN,t)}_{-k}
  f^t_k) + \Lp^{MN}(F^{(MN,t)}_0 f^t_0)\\ +
  \sum_{j>0}\sum_{|\ell|>K}\Lp^{MN}(F^{(MN,t)}_{-\ell}
  \Lp_\ell^{MNj} F^{(MN,t)}_\ell )+O(t^3).
  \end{multline*}

For $|\ell|>K$ and $j>0$, we have $\norm{\Lp_\ell^{MNj}
v}_{C^1}\leq C (1+\ell^2)\theta^{30MNj}
\norm{v}_{\CC^{1,\epsilon}_{MN}}$ for any function $v$, by
Lemma \ref{lem:ItereC1}. Hence, \eqref{eq:EstimeFMNtprecis}
enables us to replace $F_\ell^{(MN,t)}$ and
$F_{-\ell}^{(MN,t)}$ respectively with $itg_\ell$ and
$itg_{-\ell}$, the additional errors being summable and giving
a term of order $O(t^3)$. Using also the estimates on $f^t_k$
of Lemma \ref{lem:Decritftk}, we obtain
  \begin{multline*}
  (R(t)f^t)_0=-t^2 \sum_{1\leq |k|\leq
  K}\sum_{n>0}\Lp^{MN}(g_{-k}\Lp_k^{MNn} g_k) +
  \Lp^{MN}(F^{(MN,t)}_0 f^t_0)\\ -t^2 \sum_{ |\ell|>K} \sum_{j>0}
  \Lp^{MN}( g_{-\ell} \Lp_\ell^{MNj} g_\ell)+O(t^3).
  \end{multline*}
To estimate $\Lp^{MN}(F^{(MN,t)}_0 f^t_0)$, we write, in
$\CC^{1,\epsilon}_{MN}$,
  \begin{equation}
  F^{(MN,t)}_0(x)=1+itg_0(x)-\frac{t^2}{2}\int_{\Sbb^1}
  S_{MN}^Y\psi_Y(x,\omega)^2 \dd\omega + O(t^3).
  \end{equation}
Consequently, by Lemma \ref{lem:Decritftk} and since $\int f^t_0=1$,
$\int g_0=0$,
  \begin{align*}
  \int_Y \Lp^{MN}( F^{(MN,t)}_0 f^t_0)&
  =\int_Y F^{(MN,t)}_0 f^t_0
  \\&
  =1+ \int_Y itg_0 f^t_0 -
  \frac{t^2}{2}\int_Y \int_{\Sbb^1}S_{MN}^Y\psi_Y(x,\omega)^2
  f^t_0(x)\dd\omega+O(t^3)
  \\&
  =1-t^2 \sum_{n=1}^\infty \int_Y g_0 \Lp^{MNn}g_0-\frac{t^2}{2}
  \int_{Y\times \Sbb^1} S_{MN}^Y\psi_Y(x,\omega)^2 + O(t^3).
  \end{align*}
Finally, as $\lambda(t)=\int_Y \lambda(t)f^t_0=\int_Y (R(t)f^t)_0$,
we obtain
  \begin{equation}
  \lambda(t)=1-\frac{t^2}{2}
  \int_{Y\times \Sbb^1} S_{MN}^Y\psi_Y(x,\omega)^2-t^2
  \sum_{k\in \Z} \sum_{n>0}\int_Y g_{-k} \Lp_k^{MNn} g_k +
  O(t^3),
  \end{equation}
and the sum is absolutely converging. To conclude the proof, it is
therefore sufficient to show that, for any $n>0$,
  \begin{equation}
  \label{eq:posiufmlkj}
  \sum_{k\in \Z} \int_Y g_{-k} \Lp_k^{MNn} g_k = \int_{Y\times
  \Sbb^1} S_{MN}^Y \psi_Y \cdot S_{MN}^Y \psi_Y\circ \UU_Y^n.
  \end{equation}
We have
  \begin{align*}
  \int_Y g_{-k} \Lp_k^{MNn} g_k &= \int g_{-k}
  \Lp^{MNn}(e^{-itk\sum_{j=0}^{n-1}S_{MN}^Y\phi_Y \circ U_Y^j}
  g_k)
  =\int_Y g_{-k}\circ U_Y^n e^{-itk\sum_{j=0}^{n-1}S_{MN}^Y\phi_Y \circ U_Y^j}
  g_k
  \\&
  =
  \int_Y \left(\int_{\Sbb^1} S_{MN}^Y \psi_Y(U_Y^n x,
  \tilde\omega)e^{ik\tilde\omega}\dd\tilde\omega\right)
  e^{-itk\sum_{j=0}^{n-1}S_{MN}^Y\phi_Y \circ U_Y^j(x)}\times
  \\&\ \ \ \ \ \quad\quad\quad\times
  \left(\int_{\Sbb^1}S_{MN}^Y \psi_Y(x,
  \omega)e^{-ik\omega}\dd\omega\right) \dd\mu_Y(x).
  \end{align*}
Let $\omega'=\tilde\omega-\sum_{j=0}^{n-1}S_{MN}^Y\phi_Y \circ
U_Y^j(x)$, so that the previous formula becomes
  \begin{equation}
  \label{eq:opaizuklj}
  \int_Y g_{-k} \Lp_k^{MNn} g_k
  =\int_{Y} \left(\int_{\Sbb^1} S_{MN}^Y
  \psi_Y\circ\UU_Y^n(x,\omega') e^{ik\omega'}\right)
  \left(\int_{\Sbb^1} S_{MN}^Y
  \psi_Y(x,\omega)e^{-ik\omega}\dd\omega\right)\dd\mu_Y(x).
  \end{equation}
For any $u,v\in L^2(Y\times \Sbb^1)$, we have
  \begin{equation}
  \int_{Y\times\Sbb^1} uv=\sum_{k\in \Z} \int_Y \left( \int_{\Sbb^1}
  u(x,\omega')e^{ik\omega'}\dd\omega'\right)
  \left(\int_{\Sbb^1}
  v(x,\omega)e^{-ik\omega}\dd\omega\right)\dd\mu_Y(x),
  \end{equation}
where the series on the right converges absolutely.
%Indeed, if $u_k$
%and $v_k$ are the Fourier coefficients of $u$ and $v$, the
% right hand term is bounded by
%  \begin{equation}
%  \sum_{k\in \Z}\left| \int_Y u_k v_{-k}\right| \leq \sum_{k\in
%  \Z}\norm{u_k}_{L^2}\norm{v_k}_{L^2} \leq \sqrt{\sum_{k\in \Z}
%  \norm{u_k}_{L^2}^2} \sqrt{\sum_{k\in \Z}
%  \norm{v_k}_{L^2}^2}
%  =\norm{u}_{L^2}\norm{v}_{L^2}.
%  \end{equation}
This is simply Parseval's equality in each fiber $\Sbb^1$,
integrated with respect to $x$. Together with
\eqref{eq:opaizuklj}, this yields \eqref{eq:posiufmlkj} and
concludes the proof of the lemma.
\end{proof}

\begin{lem}
We have $\tilde\sigma^2=\mu^{(MN)}(X^{(MN)})\sigma^2$.
\end{lem}
Together with Lemma \ref{eq:DonneLambdattilde}, this concludes the
proof of Theorem \ref{thm:Controlelambdat}.
\begin{proof}
We will show that
  \begin{equation}
  \label{eq:pouisqfpoj}
  \tilde\sigma^2=\int_{X^{(MN)}\times \Sbb^1}
  \psi^2\dd(\mu^{(MN)}\otimes\Leb) + 2\sum_{n=1}^\infty
  \int_{X^{(MN)}\times \Sbb^1}\psi\cdot \psi\circ \UU^n
  \dd(\mu^{(MN)}\otimes\Leb).
  \end{equation}
Since $\mu^{(MN)}$ projects on $\mu^{(MN)}(X^{(MN)}) \tilde\mu$,
this will imply the result of the lemma.

It is easy to convince oneself of \eqref{eq:pouisqfpoj} by expanding
the expression of $S_{MN}^Y\psi_Y$ in $\tilde\sigma^2$ and then
gluing back together the different pieces to get the right member of
\eqref{eq:pouisqfpoj}. However, this process involves series which
are \emph{a priori} not convergent, which is a problem. We will
therefore do the computation in a different way, inspired by
\cite[Proposition 4.8]{gouezel:stable}.

Let us define a function $c$ on $X^{(MN)}\times \Sbb^1$ by
$c=\sum_{n=1}^\infty \hat\UU^n(\psi)$. This series converges by
Theorem \ref{thm:MelangeExponentielPrecis}, and defines a function
belonging to $L^p(X^{(MN)}\times \Sbb^1)$ for any $p$. Moreover,
$c=\hat\UU\psi+\hat\UU c$. Let $a$ be the restriction of $c$ to $Y$.
The previous equation implies that $a=\hat\UU_Y S_{MN}^Y\psi_Y +
\hat\UU_Y a$. As a consequence, the function $\tilde a=a-\int a$ is
equal to $\sum_{n=1}^\infty \hat\UU_Y^n(S_{MN}^Y\psi_Y)$ (and this
series is indeed converging, since $\int S_{MN}^Y\psi_Y=0$). In
particular,
  \begin{equation}
  \tilde\sigma^2=\int_{Y\times\Sbb^1} (S_{MN}^Y\psi_Y)^2 +
  2\int_{Y\times\Sbb^1}
  S_{MN}^Y\psi_Y \cdot \tilde a
  = \int_{Y\times\Sbb^1} (S_{MN}^Y\psi_Y)^2 + 2\int_{Y\times\Sbb^1}
  S_{MN}^Y\psi_Y \cdot a.
  \end{equation}
The explicit relationship between $a$ and $c$ then makes it
possible to show (as in the proof of \cite[Proposition
4.8]{gouezel:stable}) that this quantity is equal to $
\int_{X^{(MN)}\times \Sbb^1}( \psi^2 + 2 \psi c)$, which proves
\eqref{eq:pouisqfpoj} given the definition of $c$.
\end{proof}

\subsection{Reconstruction of $\hat\UU_t^n$}

Let us assume from now on that $\sigma^2>0$.

We proved in the previous paragraphs that the sequence $R_n^t$
is a perturbed renewal sequence of operators with exponential
decay, in the sense of Definition \ref{def:renouvellement}, and
that it is aperiodic if the function $\psi$ itself is
continuously aperiodic. We can therefore apply Theorem
\ref{thm:renouvellement} and get the following estimate on
$T_n^t$ (defined in \eqref{eq:DefinitTtn}):
\begin{prop}
\label{prop:EstimeTntPrecis} Let $P$ be the operator on $\B_K$
defined in Lemma \ref{lem:DonneProjecteur}. There exist
$\tau_0>0$, $c>0$, $C>0$ and $\bar\theta<1$ such that, for any
$n\in \N$, $t\in[-\tau_0,\tau_0]$ and $v\in \B_K$,
  \begin{equation}
  \norm{T_n^t v - \frac{1}{\mu^{(MN)}(X^{(MN)})}
  \left(1-\frac{\sigma^2 t^2}{2}\right)^n Pv}_{\B_K}
  \leq C (\bar \theta^n + |t| (1-ct^2)^n) \norm{v}_{\B_K}.
  \end{equation}
Moreover, if $\psi$ is continuously aperiodic, we also have for
any $|t|\in [\tau_0,t_0]$,
  \begin{equation}
  \label{eq:EstimeTntHors0}
  \norm{T_n^t v}_{\B_K} \leq C \bar\theta^n \norm{v}_{\B_K}.
  \end{equation}
\end{prop}

We recall that $T^t_n$ is also given by $T^t_n v = 1_{Y\times
[-K,K]} \hat \K^{t,n} (1_{Y\times [-K,K]}v)$. As we have a good
control on $\hat \K^t$ outside $Y\times[-K,K]$, the information
given by Proposition \ref{prop:EstimeTntPrecis} will therefore
make it possible to reconstruct precisely $\hat \K^{t,n}$. As a
first step, we will estimate $P_n^t v:=1_{Y\times \Z} \hat
\K^{t,n} (1_{Y\times \Z}v)$. As in Paragraph
\ref{par:DonneOperateursRetour}, we thus define operators
$A_n^t$, $B_n^t$ and $C_n^t$ using the kernel $\K^t$ along
trajectories of the ``random walk'' of length $n$, starting and
ending in $Y\times \Z$, with the following additional
restrictions. For the operator $A_n^t$, we only sum over the
trajectories that enter in $Y\times [-K,K]$ after a time
exactly $n$, for the operator $B_n^t$ over the trajectories
starting in $Y\times[-K,K]$ and staying out of it for the next
$n$ iterates , and for the operator $C_n^t$ over the
trajectories spending all their iterates outside of
$Y\times[-K,K]$. Formally, for $n>0$,
  \begin{equation*}
  A_n^t v(x,k)
  =\sum_{p\geq 0}\ \sum_{\substack{k_0\in[-K,K], k_1,\dots,k_{p-1},k_p=k \not\in [-K,K]\\
  x_0,x_1,\dots, x_{p-1},x_p=x\\
  \sum_{i=0}^{p-1} r^{(MN)}(x_i)=n}} \K^{t,Y}_{(x_p,k_p)\to
  (x_{p-1},k_{p-1})} \dots \K^{t,Y}_{(x_1, k_1)\to (x_0,k_0)}
  v(x_0,k_0),
  \end{equation*}
and $B_n^t$, $C_n^t$ are defined in an analogous way.

By construction, the operator $P_n^t$ satisfies:
  \begin{equation}
  P_n^t = C_n^t  + \sum_{a+i+b=n} A_a^t T_i^t B_b^t,
  \end{equation}
as long as this expression makes sense. We therefore need to
introduce different Banach spaces of functions from $Y\times
\Z$ to $\C$ such that the operators $A_n^t$, $B_n^t$ and
$C_n^t$ are well defined between these spaces. In addition to
$\BB_K$, let us denote by $\B^1$ the set of functions $v$ from
$Y\times \Z$ to $\C$ such that $\sum_{k\in \Z} (1+k^2)
\norm{v_k}_{C^1(Y)}<\infty$, with its canonical norm, and by
$\B^2$ the set of functions $v$ from $Y\times \Z$ to $\C$ such
that $\sum_{k\in \Z} \norm{v_k}_{C^1(Y)} <\infty$. We will
consider $A_a^t$ as an operator from $\B_K$ to $\B^2$, $B^t_b$
as an operator from $\B^1$ to $\B_K$, and $C^t_n$ as an
operator from $\B^1$ to $\B^2$. It should of course be checked
that these operators are bounded for these respective norms.
This is done in the following lemma.
\begin{lem}
\label{lem:EstimeABC} There exists $C>0$ such that, for any $n\in
\N^*$ and any $t\in [-t_0,t_0]$,
  \begin{equation}
  \norm{A_n^t}_{\B_K\to \B^2} \leq C|t| e^{-\epsilon n},\quad
  \norm{B^t_n}_{\B^1 \to \B_K} \leq C|t| e^{-\epsilon n},\quad
  \norm{C^t_n}_{\B^1\to \B^2} \leq C e^{-\epsilon n}.
  \end{equation}
\end{lem}
\begin{proof}
Let us start with $A_n^t$. If $\kk=(k_0,\dots,k_j)$ is an
admissible sequence, we have defined an operator $\bar
Q_{\kk}^t(z)$ in Paragraph \ref{subsec:DefinitBarQ}, by $( \bar
Q_{\kk}^t(z)v)_k=0$ if $k\not=k_j$, and $( \bar
Q_{\kk}^t(z)v)_{k_j}=Q_{\kk}^t(z) v_{k_0}$. We define an
operator $A(z,t)$ from $\B_K$ to $\B^2$ by
  \begin{equation}
  A(z,t)=\sum_{j=1}^\infty \sum_{\substack{\kk=(k_0,k_1,\dots,k_{j-1},k_j)
  \text{ admissible}\\ |k_0|\leq K, |k_j|>K}}
  \bar Q_{\kk}^t(z).
  \end{equation}
By construction, $A_n^t$ is the coefficient of $z^n$ in this
series. Moreover, summing the estimates of Lemma
\ref{lem:ControleQkt} over admissible sequences with $|k_0|\leq
K$ and $|k_j|>K$, we obtain that $A(z,t)$ is holomorphic on the
disk $\{|z|< e^{2\epsilon}\}$ (as a function from $\B_K$ to
$\B^2$). Summing the estimates \eqref{eq:ControlePerturb0} for
small $t$, we also get that $A(z,t)$ is bounded by $C|t|$
(since the number of differences in such an admissible sequence
is at least $1$). As a consequence, $A(z,t)$ is bounded by
$C|t|$ for $t\in[-t_0,t_0]$ since this inequality is trivial
outside of a neighborhood of $0$. Thus, the coefficient of
$z^n$ in $A(z,t)$ decays at least like $C|t|e^{-\epsilon n}$.
This concludes the proof of the estimate of $A_n^t$.

For $B_n^t$, we argue in the same way, using the fact that it is the
coefficient of $z^n$ in the series
  \begin{equation}
  \label{lkqjsfdlkjqsf}
  \sum_{j=1}^\infty \sum_{\substack{\kk=(k_0,k_1,\dots,k_{j-1},k_j)
  \text{ admissible}\\ |k_0|> K, |k_j|\leq K}}
  \bar Q_{\kk}^t(z).
  \end{equation}
As $\norm{Q_{\kk}^t(z)}_{C^1(Y) \to C^1(Y)} \leq C|t|(1+k_0^2)
\theta^{20MNj} \prod_{i=1}^j \gamma_{d_i}^{1/3}$ by Lemma
\ref{lem:ControleQkt}, we also have
  \begin{equation}
  \norm{\bar Q_{\kk}^t(z)}_{\B^1 \to \B_K} \leq C|t| \theta^{20MNj}
  \prod_{i=1}^j \gamma_{d_i}^{1/3}.
  \end{equation}
Since this quantity is summable with respect to $\kk$, the
series \eqref{lkqjsfdlkjqsf} is holomorphic on the disk $\{|z|<
e^{2\epsilon}\}$ and bounded by $C|t|$. We conclude as above.

Finally, $C_n^t$ is the coefficient of $z^n$ in the series
  \begin{equation}
  \sum_{j=1}^\infty \sum_{\substack{\kk=(k_0,k_1,\dots,k_{j-1},k_j)
  \text{ admissible}\\ |k_0|> K, |k_j|> K}}
  \bar Q_{\kk}^t(z),
  \end{equation}
which defines an holomorphic function from $\B^1$ to $\B^2$ in
the disk $\{|z|< e^{2\epsilon}\}$ (by summing the estimates of
Lemma \ref{lem:ControleQkt}). This yields the desired estimate
for $C_n^t$.
\end{proof}

We have defined a projection $P$ on $\B_K$, which can be extended to
an operator from $\B^1$ to $\B^2$, as follows: $(Pv)_k=0$ if $k=0$,
and $(Pv)_0= \int_Y v_0 \dd\mu_Y$.
\begin{cor}
\label{cor:ControleErreur} There exist constants $\eps>0$,
$c>0$, $C>0$ and $\bar\theta<1$ such that, for any $n\in \N$,
$t\in [-\eps,\eps]$ and $v\in \B^1$,
  \begin{equation}
  \norm{P_n^t v - \frac{1}{\mu^{(MN)}(X^{(MN)})}
  \left(1-\frac{\sigma^2 t^2}{2}\right)^n Pv
  }_{\B^2}
  \leq C (\bar \theta^n + |t| (1-ct^2)^n) \norm{v}_{\B^1}.
  \end{equation}
Moreover, if $\psi$ is continuously aperiodic, one also has for any
$|t|\in [\eps,t_0]$
  \begin{equation}
  \label{eq:PntHorsVoisinage0}
  \norm{P_n^t v}_{\B^2} \leq C \bar\theta^n \norm{v}_{\B^1}.
  \end{equation}
\end{cor}
\begin{proof}
We write $P_n^t=A_0^t T_n^t B_0^t + C_n^t + \sum_{a+i+b=n,\ i<n}
A_a^t T_i^tB_b^t$, as an operator from $\B^1$ to $\B^2$. The term
$A_0^t T_n^t B_0^t$ gives the desired asymptotics, by Proposition
\ref{prop:EstimeTntPrecis} (and since $A_0^t$ and $B_0^t$ are simply
trivial extension and restriction operators). The term $C_n^t$ is
$O(\bar\theta^n)$ by Lemma \ref{lem:EstimeABC}. Hence, we should
estimate the sum $\sum_{a+i+b=n,\ i<n} A_a^t T_i^tB_b^t$, whose norm
is bounded by
  \begin{equation}
  C|t| \sum_{a+i+b=n} e^{-\epsilon a} ( \bar \theta^i + (1-ct^2)^i)
  e^{-\epsilon b},
  \end{equation}
again by Lemma \ref{lem:EstimeABC} and Proposition
\ref{prop:EstimeTntPrecis}. The term $\sum e^{-\epsilon a}
\bar\theta^i e^{-\epsilon b}$ is exponentially small in $n$, while
the remaining term is bounded by
  \begin{equation*}
  |t| \sum_{i+j=n} (j+1)e^{-\epsilon j} (1-ct^2)^i
  \leq C|t| (1-ct^2)^n \sum_{j=0}^n
  \bigl((1-ct^2)^{-1}e^{-\epsilon}\bigr)^j
  \leq \frac{C|t| (1-ct^2)^n}{1- (1-ct^2)^{-1}e^{-\epsilon}}.
  \end{equation*}
This is bounded by $C|t| (1-ct^2)^n$ if $t$ is small enough.

When $\psi$ is continuously aperiodic, the equation
\eqref{eq:PntHorsVoisinage0} is proved in the same way by combining
\eqref{eq:EstimeTntHors0} and Lemma \ref{lem:EstimeABC}.
\end{proof}

The next step in the reconstruction of $\hat\K^{t,n}$ is to
understand $\tilde P_n^t v:= 1_{Y\times \Z} \hat \K^{t,n}(v)$. We
will let this operator act on the space $\B^0$ of functions $v$ from
$X^{(MN)}\times \Z$ to $\C$ such that $\sum_{k\in \Z} (1+|k|^3)
\norm{v_k}_{C^1(X^{(MN)})} <\infty$, and take its values in $\B^2$.
Let us also define an operator $\tilde P$ from $\B^0$ to $\B^2$ by
$(\tilde P v)_k=0$ for $k\not=0$, and $(\tilde P
v)_0=\int_{X^{(MN)}} v_0 \dd\tilde \mu^{(MN)}$ (recall that $\tilde
\mu^{(MN)}$ is a probability measure on $X^{(MN)}$, whose
restriction to $Y$ is $\mu_Y / \mu^{(MN)}(X^{(MN)})$).

\begin{prop}
\label{prop:ControleTildeP} There exist constants $\eps>0$, $c>0$,
$C>0$ and $\bar\theta<1$ such that, for any $n\in \N$,
$t\in[-\eps,\eps]$ and $v\in \B^0$,
  \begin{equation}
  \label{eq:tildePntPres0}
  \norm{\tilde P_n^t v -
  \left(1-\frac{\sigma^2 t^2}{2}\right)^n \tilde P v
  }_{\B^2}
  \leq C (\bar \theta^n + |t| (1-ct^2)^n) \norm{v}_{\B^0}.
  \end{equation}
Moreover, if $\psi$ is continuously aperiodic, one also has for any
$|t|\in [\eps,t_0]$
  \begin{equation}
  \label{eq:tildePntHorsVoisinage0}
  \norm{\tilde P_n^t v}_{\B^2} \leq C \bar\theta^n \norm{v}_{\B^0}.
  \end{equation}
\end{prop}
\begin{proof}
Let us define an operator $D_n^t$, which corresponds to considering
the trajectories of the ``random walk'' starting from $Y\times \Z$
and staying outside of $Y\times \Z$ during a time $n$, so that
$\tilde P_n^t=\sum_{i+j=n} P_i^t D_j^t$. Formally, for $x\in Y$,
  \begin{equation}
  D_n^t v(x,k)=
  \sum_{\substack{k_0,\dots, k_n=k\\
  x_0,\dots, x_n=x\\ x_i \not\in Y\text{ for }0\leq i<n}}
  \K^t_{(x_n,k_n)\to (x_{n-1},k_{n-1})}\dots
  \K^t_{(x_1, k_1)\to
  (x_0, k_0)} v(x_0,k_0).
  \end{equation}
We will first study $D_n^t$, as an operator from $\B^0$ to
$\B^1$. As the dynamics of $U$ between two returns to $Y$ is
trivial, $D_n^t$ can be explicitly described as follows. Recall
that a point $x$ in $X^{(MN)}$ is a pair $(y,i)$ where $y\in Y$
and $i<r^{(MN)}(y)$. The preimages of $(x,0)$ under $U^n$ which
do not enter $Y$ in between are exactly the points $(hx,
r^{(MN)}(hx)-n)$ where $h\in \HH_{MN}$ is an inverse branch of
$T_Y^{MN}$ whose return time $r^{(MN)}\circ h$ is $>n$. Let
$v\in \B^0$. For $k,l\in \Z$, let us define a function
$v^n_{k,l}$ on $Y$ by
  \begin{equation*}
  v^n_{k,l}(y)=1_{r^{(MN)}(y)>n} v_l(y, r^{(MN)}(y)-n)e^{-ik S_n
  \phi(y, r^{(MN)}(y)-n)} (e^{it S_n \psi})_{k-l}(y, r^{(MN)}(y)-n).
  \end{equation*}
Here, $(y,r^{(MN)}(y)-n)$ is a point in $X^{(MN)}$,
$e^{-ikS_n\phi}$ is a function on $X^{(MN)}$ and
$(e^{itS_n\psi})_{k-l}$ is the $k-l$-th Fourier coefficient (in
the $\omega$ direction) of the function $e^{itS_n\psi}$ on
$X^{(MN)}\times \Sbb^1$, so it is also a function on
$X^{(MN)}$. We have defined $v^n_{k,l}$ so that $D_n^t v(x,k)=
\sum_l \Lp^{MN} v^n_{k,l} (x)$.

Let us now estimate $\norm{D_n^t v}_{\B^1}$ in terms of
$\norm{v}_{\B^0}$. As $\psi$ belongs to $C^{5,1}$, the $k-l$-th
Fourier coefficient of $e^{it S_n \psi}$ is bounded by $C n^5 /
(1+|k-l|^5)$. As $r^{(MN)}(x)>n$, we get
  \begin{equation}
  |v^n_{k,l}(x)|\leq C \norm{v_l}_{C^0} \frac{n^5}{1+|k-l|^5} \leq C
  \norm{v_l}_{C^0}e^{-\epsilon n} \frac{e^{2\epsilon
  r^{(MN)}(x)}}{1+|k-l|^5}
  \end{equation}
and, for any inverse branch $h$,
  \begin{equation}
  \norm{D (v^n_{k,l}\circ h)}_{C^0} \leq C \norm{v_l}_{C^1}
  (1+|k|)n \frac{n^5}{1+|k-l|^5} \leq C \norm{v_l}_{C^1}(1+|k|)e^{-\epsilon
  n}\frac{e^{2\epsilon
  r^{(MN)}(x)}}{1+|k-l|^5}.
  \end{equation}
As a consequence,
  \begin{equation}
  \norm{v^n_{k,l}}_{\CC^{1,2\epsilon}_{MN}}\leq
  \frac{C(1+|k|)}{1+|k-l|^5} \norm{v_l}_{C^1}e^{-\epsilon
  n}.
  \end{equation}
By Theorem \ref{thm:MainContraction}, $\norm{\Lp^{MN}
v^n_{k,l}}_{C^1(Y)}\leq C
\|v^n_{k,l}\|_{\CC^{1,2\epsilon}_{MN}}$. Finally,
  \begin{equation}
  \norm{D_n^t v}_{\B^1} = \sum_k (1+|k|^2) \norm{(D_n^t v)_k}_{C^1(Y)}
  \leq Ce^{-\epsilon n} \sum_{k,l} \frac{ 1+|k|^3}{1+|k-l|^5}
  \norm{v_l}_{C^1}.
  \end{equation}
If $l$ is fixed,
  \begin{equation}
  \sum_k \frac{ 1+|k|^3}{1+|k-l|^5}=\sum_{j} \frac{ 1+|j+l|^3}{1+|j|^5}
  \leq C \sum_j \frac{1+|j|^3+|l|^3}{1+|j|^5}
  \leq C (1+|l|^3).
  \end{equation}
Consequently,
  \begin{equation}
  \norm{D_n^t v}_{\B^1} \leq C e^{-\epsilon n} \norm{v}_{\B^0}.
  \end{equation}

In $\tilde P_n^t v= \sum_{i+j=n} P_i^t D_j^t v$, let us replace
$P_i^t$ with $(1-\sigma^2 t^2/2)^i P
/\mu^{(MN)}(X^{(MN)})+E_i^t$, where $E_i^t$ is an error term.
The control of $E_i^t$ given by Corollary
\ref{cor:ControleErreur}, combined with the computation made at
the end of the proof of this lemma, gives
  \begin{equation}
  \sum_{i+j=n} \norm{E_i^t D_j^t}_{\B^0\to\B^2} \leq C \sum_{i+j=n} (\bar
  \theta^i + |t| (1-ct^2)^i) e^{-\epsilon j}
  \leq C' ( \bar\theta^i+ |t| (1-ct^2)^n).
  \end{equation}
Hence, there is only one term left to be estimated in $\tilde P_n^t
v$, with frequency $0$, given by
  \begin{equation}
  I^t_n:=\frac{1}{\mu^{(MN)}(X^{(MN)})}\sum_{i+j=n}
  \left(1-\frac{\sigma^2 t^2}{2}\right)^i \int_Y (D_j^t v)_0
  \dd\mu_Y.
  \end{equation}
For all $u,v\in \R$ holds $|e^u-e^v|\leq |u-v|e^{\max(u,v)}$.
As $\left|\int_Y (D_j^t v)_0\right| \leq C e^{-\epsilon
j}\norm{v}_{\B^0}$, we obtain
  \begin{align*}
  \left| \sum_{j=0}^n
  \left(1-\frac{\sigma^2 t^2}{2}\right)^{n-j} \int_Y (D_j^t v)_0
  \dd\mu_Y -
  \left(1-\frac{\sigma^2 t^2}{2}\right)^n \sum_{j=0}^n \int_Y (D_j^t
  v)_0\right|
  \!\!\!\!\!\!\!\!\!\!\!\!\!\!\!\!\!\!\!\!\!\!\!\!\!\!\!\!
  \!\!\!\!\!\!\!\!\!\!\!\!\!\!\!\!\!\!\!\!\!\!\!\!\!\!\!\!
  \!\!\!\!\!\!\!\!\!\!\!\!\!\!\!\!\!\!\!\!\!\!\!\!\!\!\!\!
  \!\!\!\!\!\!\!\!\!\!\!\!\!\!\!\!\!\!\!\!\!\!\!\!\!\!\!\!
  &
  \\&
  \leq C \left(1-\frac{\sigma^2 t^2}{2}\right)^n \sum_{j=0}^n j
  \left|\log\left(1-\frac{\sigma^2 t^2}{2}\right)\right|
  \left(1-\frac{\sigma^2 t^2}{2}\right)^{-j} e^{-\epsilon j}\norm{v}_{\B^0}
  \\&
  \leq C t^2 \left(1-\frac{\sigma^2 t^2}{2}\right)^n
  \norm{v}_{\B^0}.
  \end{align*}
Let us define a function $f$ on $X^{(MN)}\times \Sbb^1$ by
$f(x,\omega)=\sum_k v_k(x) e^{ik\omega}$. If $Z_j \subset
X^{(MN)}$ denotes the set of points in $X^{(MN)}$ which enter
into $Y$ after exactly $j$ iterates, we have
  \begin{equation}
  \int_Y (D_j^t v)_0 \dd\mu_Y = \int_{Z_j \times \Sbb^1} f e^{itS_j
  \psi} \dd(\mu^{(MN)} \otimes \Leb).
  \end{equation}
Since the measure of $Z_j$ decays exponentially fast,
  \begin{equation}
  \left| \int_Y (D_j^t v)_0 \dd\mu_Y - \int_{Z_j\times \Sbb^1} f \dd(\mu^{(MN)} \otimes
  \Leb) \right| \leq C \int_{Z_j \times \Sbb^1} |t| j \norm{f}_{C^0}
  \leq C |t| \bar\theta^j \norm{v}_{\B^0}.
  \end{equation}
Finally,
  \begin{multline*}
  \left|\sum_{j=0}^n \int_{Z_j\times \Sbb^1} f \dd(\mu^{(MN)} \otimes
  \Leb) -\int_{X^{(MN)}\times \Sbb^1} f \dd(\mu^{(MN)} \otimes
  \Leb)\right|
  \\
  \leq C \norm{f}_{C^0} \sum_{j=n+1}^\infty
  \mu^{(MN)}(Z_j)
  \leq C \norm{v}_{\B^0} \bar\theta^n.
  \end{multline*}
Combining these different estimates, we obtain
  \begin{align*}
  I_n^t&=\left(1-\frac{\sigma^2 t^2}{2}\right)^n \frac{1}{\mu^{(MN)}(X^{(MN)})}
  \int_{X^{(MN)}\times \Sbb^1} f \dd(\mu^{(MN)} \otimes \Leb)
  + O ( \bar\theta^n + |t|(1-ct^2)^n)
  \\&
  = \left(1-\frac{\sigma^2 t^2}{2}\right)^n \int_{X^{(MN)}} v_0 \dd\tilde\mu^{(MN)}
  + O ( \bar\theta^n + |t|(1-ct^2)^n).
  \end{align*}
This proves \eqref{eq:tildePntPres0}. Finally,
\eqref{eq:tildePntHorsVoisinage0} is proved in the same way, by
using \eqref{eq:PntHorsVoisinage0}.
\end{proof}

Let $\hat\UU_t$ denote the operator acting on functions on
$X^{(MN)}\times\Sbb^1$ by $\hat\UU_t(v)=\hat\UU(e^{it\psi}v)$,
where $\hat\UU$ is the transfer operator associated to $\UU$.

\begin{thm}
\label{thm:MelangeExponentielPrecisLocal} Assume $\sigma^2>0$.
Then there exist constants $\eps>0$, $c>0$, $C>0$ and
$\bar\theta<1$ such that, for any $C^{5,1}$ function
$v:X^{(MN)}\times\Sbb^1 \to \C$, for any $n\in \N$, for any
$t\in[-\eps,\eps]$ and for any $(x,\omega)\in X^{(MN)}\times
\Sbb^1$ such that $h(x)\leq n/2$,
  \begin{equation}
  \label{eq:ControletPres0}
  \left| \hat\UU_t ^n v(x,\omega) -\left(1-\frac{\sigma^2 t^2}{2}\right)^n
  \int v \dd(\tilde
  \mu^{(MN)}\otimes \Leb) \right| \leq C (1+h(x))(\bar \theta^{n} + |t|(1-ct^2)^{n})\norm{v}_{C^{5,1}}.
  \end{equation}
Moreover, if $\psi$ is continuously aperiodic, we also have for
any $|t|\in [\eps,t_0]$ and for any $(x,\omega)$ with $h(x)\leq
n/2$
  \begin{equation}
  \label{eq:ControletLoinO}
  \left| \hat\UU_t ^n v(x,\omega) \right| \leq C
  \bar\theta^{n} \norm{v}_{C^{5,1}}.
  \end{equation}
\end{thm}
Note that this theorem implies Theorem
\ref{thm:MelangeExponentielPrecis}, taking simply $t=0$ (and a
different value of $\bar\theta$).
\begin{proof}
Define $w$ in $\B^0$ by $w(x,k)=\int_{\Sbb^1} v(x,\omega)
e^{-ik \omega}\dd\omega$, so that $v(x,\omega)=\sum
w(x,k)e^{ik\omega}$. As $v\in C^{5,1}$, $w$ belongs to $\B^0$
and $\norm{w}_{\B^0}\leq C \norm{v}_{C^{5,1}}$.

For $x\in Y$, we have $\hat\UU_t^n v(x,\omega)= \sum_{k\in \Z}
(\tilde P_n^t w)_k(x)e^{ik\omega}$ by construction of $\tilde
P_n^t$. Hence, Proposition \ref{prop:ControleTildeP} implies that,
for $x\in Y$ and $t\in[-\eps,\eps]$
  \begin{align*}
  \left|\hat\UU_t^n v(x,\omega) -
  \left(1-\frac{\sigma^2 t^2}{2}\right)^n\int v\right|&
  \leq \left|(\tilde P_n^t w)_0(x) -
  \left(1-\frac{\sigma^2 t^2}{2}\right)^n\int w_0\right| +
  \sum_{k\in \Z^*} |(\tilde P_n^t w)_k(x)|
  \\&
  \leq \norm{\tilde P_n^t w -
  \left(1-\frac{\sigma^2 t^2}{2}\right)^n \tilde
  Pw}_{\B^2}
  \\ &
  \leq C( \bar\theta^n + |t|(1-ct^2)^n) \norm{w}_{\B^0}
  \leq C( \bar\theta^n + |t|(1-ct^2)^n) \norm{v}_{C^{5,1}}.
  \end{align*}
This proves \eqref{eq:ControletPres0} for the points $x$ with
$h(x)=0$.

Assume now that $j=h(x)\in (0,n/2]$. Let $x'$ be such that
$U^jx'=x$, and let $\omega'=\omega-S_j\phi(x')$, so that
$\UU^j(x',\omega')=(x,\omega)$. Then $\hat \UU_t^n v(x,\omega)=
e^{itS_j \psi(x',\omega')}\hat \UU_t^{n-j} v(x', \omega')$. Using
the result for $(x',\omega')$, we get
  \begin{equation}
  \label{eq:ControlejnPresque}
  \left|\hat\UU_t^n v(x,\omega)- e^{itS_j
  \psi(x',\omega')}\left(1-\frac{\sigma^2 t^2}{2}\right)^{n-j}
  \int v \right| \leq C (\bar\theta^{n-j} +
  |t|(1-ct^2)^{n-j})\norm{v}_{C^{5,1}}.
  \end{equation}
Since $n-j\geq n/2$, this last term is bounded by $\bar\theta^{n/2}+
|t| (1-ct^2)^{n/2}$, which is compatible with
\eqref{eq:ControletPres0} (upon changing the values of $\bar\theta$
and $c$).

Moreover, $|e^{itS_j \psi(x',\omega')} -1|\leq C|t| j$.
Replacing $e^{itS_j \psi(x',\omega')}$ by $1$ in
\eqref{eq:ControlejnPresque}, we add an error which is bounded
by $C|t| h(x) (1-\sigma^2 t^2/2)^{n/2}$. This is again
compatible with \eqref{eq:ControletPres0}. Finally,
  \begin{equation*}
  \left|\left(1-\frac{\sigma^2 t^2}{2}\right)^{n-j} -
  \left(1-\frac{\sigma^2 t^2}{2}\right)^{n}\right|
  \leq j\left| \log\left(1-\frac{\sigma^2 t^2}{2}\right)\right|
  \left(1-\frac{\sigma^2 t^2}{2}\right)^{n-j}
  \leq Cjt^2 (1-ct^2)^{n/2},
  \end{equation*}
still compatible with \eqref{eq:ControletPres0}. Doing all these
substitutions, we obtain \eqref{eq:ControletPres0}.

Finally, \eqref{eq:ControletLoinO} is proved in the same way, by
using \eqref{eq:tildePntHorsVoisinage0}.
\end{proof}

\begin{proof}[Proof of Theorem
\ref{thm:EstimeesValeurPropreFinal}] Theorem
\ref{thm:MelangeExponentielPrecis} enabled us to prove Theorem
\ref{thm:MelangeExponentiel}, page
\pageref{proof:DemontreMelangeExponentiel}. The same arguments make
it possible to deduce Theorem \ref{thm:EstimeesValeurPropreFinal}
from Theorem \ref{thm:MelangeExponentielPrecisLocal}, when
$d^{(MN)}=1$.

When $d=d^{(MN)}>1$, let us show \eqref{eq:GoodEstimeFinal}
(\eqref{eq:ControleHorsVois0} is analogous). Applying the
previous arguments to the transformation $U^d$, which is
mixing, we almost obtain \eqref{eq:GoodEstimeFinal} for times
$n$ of the form $kd$, with a slight difference: since
$\sigma^2$ is replaced with
  \begin{equation}
  \int (S_d \psi)^2 + 2\sum_{j=1}^\infty (S_d\psi) (S_d\psi)\circ
  \TT^{jd}
  = d \sigma^2,
  \end{equation}
we in fact obtain
  \begin{multline*}
  \left| \int e^{itS_{kd} \psi}\cdot f\circ \TT^n \cdot g \dd(\tilde\mu\otimes
  \Leb)- \left( 1-d\frac{\sigma^2 t^2}{2}\right)^k
  \left(\int f\dd(\tilde\mu\otimes \Leb)\right)
  \left(\int g\dd(\tilde\mu\otimes \Leb)\right)\right|
  \\
  \leq C (\bar\theta^k + |t| (1-ct^2)^k)
  \norm{f}_{L^\infty}\norm{g}_{C^6}.
  \end{multline*}
To really obtain \eqref{eq:GoodEstimeFinal}, we thus have to
bound $(1-\sigma^2t^2/2)^{kd}-(1-d \sigma^2 t^2/2)^k$. We have
  \begin{align*}
  \left| \left(1-\frac{\sigma^2 t^2}{2}\right)^{kd} -
  \left(1-d\frac{\sigma^2 t^2}{2}\right)^{k} \right|
  \!\!\!\!\!\!\!\!\!\!\!\!\!\!\!\!\!\!\!\!\!\!\!\!\!\!\!\!\!\!
  \!\!\!\!\!\!\!\!\!\!\!\!\!\!\!\!\!\!\!\!\!\!\!\!\!\!\!\!\!\!\!\!\!\!\!\!
  \!\!\!\!&
  \\&
  \leq \left| kd \log\left(1-\frac{\sigma^2 t^2}{2}\right)
  - k \log\left(1-d\frac{\sigma^2 t^2}{2}\right)\right| \cdot
  \max\left(\left(1-\frac{\sigma^2 t^2}{2}\right)^{kd},
  \left(1-d\frac{\sigma^2 t^2}{2}\right)^{k} \right)
  \\&
  \leq C k |t|^4 (1-ct^2)^k.
  \end{align*}
By \eqref{eq:TechniqueControlejttrois}, this term is bounded by $C
t^2 (1-ct^2/2)^k$. This concludes the proof for times $n=kd$.

If $n$ is a general time, it can be written as $kd+r$ with
$0\leq r <d$. The theorem at time $kd$, applied to the
functions $e^{itS_r \psi} f\circ \TT^r$ and $g$ (respectively
bounded and H\"{o}lder continuous) gives almost the result, the
factor $(1-\sigma^2 t^2/2)^n$ simply being replaced with
$(1-\sigma^2 t^2/2)^{kd}$. As above, one checks that the
resulting additional error term is still compatible with
\eqref{eq:GoodEstimeFinal}.
\end{proof}

\subsection{Proof of Theorem \ref{thm:LimiteTCL}}

Assume first that $\psi$ is a $C^6$ function, with $\sigma^2>0$.
Theorem \ref{thm:EstimeesValeurPropreFinal} for $f=g=1$ shows that
the characteristic function of $S_n\psi/\sqrt{n}$ converges to
$e^{-\sigma^2 t^2/2}$, which is equivalent to the convergence of
$S_n\psi/\sqrt{n}$ towards the gaussian distribution
$\boN(0,\sigma^2)$. This concludes the proof in this case.

Assume now that $\psi$ is only $C^\alpha$, with zero average,
and with $\sigma^2>0$. Let $\psi_\epsilon$ be a $C^6$ function,
close to $\psi$ in $C^{\alpha/2}$, with corresponding
asymptotic variance $\sigma_\epsilon^2$. Theorem
\ref{thm:MelangeExponentiel} (applied in $C^{\alpha/2}$) shows
that the variance of $S_n(\psi-\psi_\epsilon)/\sqrt{n}$ is
uniformly small in $n$. This implies on the one hand that the
distributions of $S_n\psi/\sqrt{n}$ and
$S_n\psi_\epsilon/\sqrt{n}$ are close, and on the other hand
that $\sigma_\epsilon^2$ is close to $\sigma^2$. In particular,
if $\epsilon$ is small enough, $\sigma_\epsilon^2>0$. As
$S_n\psi_\epsilon/\sqrt{n}$ converges to
$\boN(0,\sigma_\epsilon^2)$, this implies that
$S_n\psi/\sqrt{n}$ is close in distribution to
$\boN(0,\sigma^2)$ if $n$ is large enough. Therefore,
$S_n\psi/\sqrt{n}$ is indeed converging to
$\boN(0,\sigma^2)$.\qed

\subsection{Regularity in the cohomological equation}
\label{subsec:Cohom}

\begin{proof}[Proof of Proposition
\ref{prop:CaracteriseSigma2}] We proved half of the proposition in
Proposition \ref{prop:CaracteriseSigma2alpha}. It remains to prove
that, if $\psi=f-f\circ \TT$ for some measurable $f$, then
$\sigma^2=0$. If $\sigma^2>0$, Theorem \ref{thm:LimiteTCL} implies
that $S_n\psi/\sqrt{n}$ converges to a gaussian distribution.
However, $S_n\psi/\sqrt{n}=(f-f\circ \TT^n)/\sqrt{n}$ converges in
distribution to $0$, which is a contradiction. Hence, $\sigma^2=0$.
\end{proof}

\begin{proof}[Proof of Proposition
\ref{prop:CaracteriseAperiodique}] Let $\psi: X\times\Sbb^1\to
\R$ be a $C^6$ function. We have to show that $\psi$ is
periodic if and only if $\psi$ is continuously periodic.

If $\psi$ is continuously periodic, it is trivially periodic.
Conversely, suppose that $\psi$ is continuously aperiodic, but it is
nevertheless possible to write $\psi=u-u\circ \TT+a\mod \lambda$,
where $u$ is measurable and $a\in\R$.

If $\sigma^2$ vanished, $\psi$ would be continuously periodic by
Proposition \ref{prop:CaracteriseSigma2}, which is a contradiction.
Hence $\sigma^2>0$. As $\psi$ is continuously aperiodic, it
satisfies Theorem \ref{thm:EstimeesValeurPropreFinal} (because
\eqref{eq:ControleHorsVois0} has been proved under the sole
assumption of continuous aperiodicity). In particular, for $t\not=0$
and for any functions $f,g$ which are respectively bounded and
$C^6$, $\int e^{itS_n \psi} f\circ \TT^n g\to 0$. By density, this
convergence to $0$ holds for any $f,g\in L^2$. However, for
$t=2\pi/\lambda$, $f=e^{itu}$ and $g=e^{-itu}$,
  \begin{equation}
  \int e^{itS_n \psi} f\circ \TT^n g
  = \int e^{it (u-u\circ \TT^n+na)} e^{itu\circ \TT^n} e^{-itu}
  =e^{itna},
  \end{equation}
which does not converge to $0$. This is a contradiction.
\end{proof}

\section{Proofs for Farey sequences}
\label{sec:PreuvesFarey}

\subsection{A general criterion for the weak Federer property}

We would like to prove that some measures $\mu$ satisfy the weak
Federer property. In the introduction, we have seen that this
property is quite easy to check for Lebesgue measure. However, in
view of the application to Farey sequences, it is desirable to have
a sufficiently simple criterion, that does not apply only to
absolutely continuous measures. In this paragraph, we describe such
a criterion.

Let us consider a riemannian manifold $Z$ endowed with a
measure $\mu$ such that, for any $\rho>0$, $\inf_{x\in Z}
\mu(B(x,\rho))>0$. We assume that $Z$ is partitioned in a
finite number of subsets $Y_1,\dots, Y_p$, and that each set
$Y_j$ admits a (finite or countable) subpartition modulo $0$,
into sets $(W_{l,j})_{l\in \Lambda(j)}$. Let also $\TY$ be a
map which sends each set $W_{l,j}$ diffeomorphically to one of
the $Y_k$. We can define $\HH_n$ as the set of inverse branches
of $\TY^n$. Such an inverse branch $h$ is not defined on the
whole space $Z$, only on one of the sets $Y_j=Y_{j(h)}$. We
assume that:
\begin{enumerate}
\item There exist $\kappa>1$ and $C_{l,j}$ such that, for any $x\in
W_{l,j}$ and $v$ tangent at $Z$ in $x$, $\kappa\norm{v}\leq
\norm{D\TY(x) v}\leq C_{l,j} \norm{v}$.
\item Let $J(x)$ be the inverse of the jacobian of $\TY$ with
respect to $\mu$. There exists $C>0$ such that, for any $h\in
\HH_1$, $\norm{D((\log J)\circ h)}\leq C$.
\item For any $\tC>1$, there exist $\tD>1$ and $\eta_0>0$ such that,
for any $\eta<\eta_0$, for any $1\leq j\leq p$, there exist
disjoint balls $B(x_1,\tC\eta),\dots, B(x_k, \tC\eta)$ which
are compactly included in $Y_j$, sets $A_1,\dots,A_k$ with
$A_i\subset B(x_i, \tD\tC\eta)\cap Y_j$ such that, for any
$x'_i\in B(x_i, (\tC-1)\eta)$, holds $\mu(B(x'_i,\eta))\geq
\mu(A_i)/\tD$, and a finite number of inverse branches
$h_1,\dots, h_\ell \in \HH_1$ defined respectively on
$Y_{j_1},\dots,Y_{j_\ell}$ such that, for any $i\in[1,\ell]$,
there exist $x\in Y_{j_i}$ and $v$ a unit tangent vector at $x$
with
  \begin{equation}
  \label{eq:TaillehiGrande}
  \norm{Dh_i(x)v}\geq \tC\eta,
  \end{equation}
such that:
  \begin{equation}
  \label{balls_included}
  \bigcup_{i=1}^k B(x_i, \tC \eta) \subset \bigcup_{i=1}^k A_i
  \end{equation}
and
  \begin{equation}
  \label{covers_all}
  Y_j= \left( \bigcup_{i=1}^k A_i \right) \sqcup
  \left(\bigsqcup_{i=1}^\ell h_i(Y_{j_i}) \right) \mod 0.
  \end{equation}
\item The transformation $\TY$ is uniformly quasi-conformal, in the
following sense: there exists $K>0$ such that, for any
$h\in\bigcup_{n\in \N}\HH_n$ defined on a set $Y_j$, for any
$x,x'\in Y_j$ and any unit tangent vectors $v$ and $v'$
respectively at $x$ and $x'$,
  \begin{equation}
  \norm{ Dh(x)v} \leq K \norm{Dh(x')v'}.
  \end{equation}
\end{enumerate}
The first two properties are uniform expansion properties,
analogous to the similar requirements on $T_Y$ in Definition
\ref{def:NonUnifDilatant}. The difference is that the full
shift structure has been replaced by a subshift of finite type,
since such a structure will naturally appear in the proofs for
Farey sequences. The third property is a kind of weak Federer
property, but not on the whole space, rather on the images of
branches whose size is at most $\tC\eta$ (by the requirement
\eqref{eq:TaillehiGrande}). It is therefore much easier to
check than the true weak Federer property. Finally, the last
property of uniform quasi-conformality will enable us to
iterate the dynamics, to get information at scales which are
not covered by the third assumption.

\begin{prop}
\label{prop:CritereFederer} Under the previous assumptions, the sets
$h(Y_{j(h)})$ (for $h\in\bigcup_{n\in\N}\HH_n$) uniformly have the
weak Federer property (for the measure $\mu$).
\end{prop}
\begin{proof}
The quasi-conformality assumption shows that it is sufficient
to prove that each set $Y_j$ satisfies the weak Federer
property: if sets $A_i$ as in the definition of the weak
Federer property can be constructed on $Y_j$, they can be
transported to $h(Y_j)$ by the map $h$. In this process, one
loses only harmless constant factors, and this implies the
uniform weak Federer property. From this point on, we shall
therefore work only on $Y_j$, for each $1\leq j\leq p$.

We want to constructs sets $A_i$ as in the definition of the
weak Federer property. The third assumption of the proposition
gives some of these sets, but to get the other ones we will
need to iterate the dynamics. Thus, the construction will be
inductive.

For any $1\leq j\leq p$, let us fix a point $a_j\in Y_j$, and a
unit tangent vector $v_j$ at $a_j$. Let also $\rho>0$ be such
that the balls $B(a_j,\rho)$ are compactly included in $Y_j$.
Fix a constant $C$ for which one wants to prove the weak
Federer property, and consider $\eta$ small enough.  We will
say that an inverse branch $h\in \HH_n$, defined on $Y_j$, is
$(C,\eta)$--\emph{good}, or simply \emph{good}, if
$\norm{Dh(a_j)v_j} \geq K C \eta / \rho$.

We will prove the following fact: \emph{there exists a constant
$M$ such that, if $h\in \HH_n$ is a good branch defined on
$Y_j$, then there exist disjoint balls
$B(x_1,C\eta),\dots,B(x_k,C\eta)$ compactly included in
$h(Y_j)$, sets $A_1,\dots,A_k$ with $A_i\subset h(Y_j)\cap
B(x_i, MC\eta)$ such that any ball $B(x'_i,\eta)$ included in
$B(x_i,C\eta)$ satisfies $\mu(B(x'_i,\eta))\geq \mu(A_i)/M$,
and good branches $h_1,\dots,h_\ell \in \HH_{n+1}$ defined
respectively on $Y_{j_1},\dots,Y_{j_\ell}$ such that
 \begin{equation}
  \bigcup_{i=1}^k B(x_i, C \eta) \subset \bigcup_{i=1}^k A_i
  \end{equation}
and
  \begin{equation}
  h(Y_j)= \left( \bigcup_{i=1}^k A_i \right) \sqcup
  \left(\bigsqcup_{i=1}^\ell h_i(Y_{j_i}) \right).
  \end{equation}
}

This fact easily implies the proposition: we first apply it to
the inverse branch $\Ide_{Y_j}$ (which is obviously good if
$\eta$ is small enough), and then by induction to the inverse
branches which are produced by the fact at the previous step.
This process terminates, since there is no good branch in
$\HH_n$ if $n$ is large enough.

To prove that fact, we will use the assumption (3) for the
constant $\tC=\max(K^2 C, K^4 C/\rho)$. Let $\eta_0$ and
$\tD>0$ be given by (3), for this value of $\tC$. Let
$\eta<\eta_0$. Let $h\in \HH_n$ be a good branch, defined on a
set $Y_j$.

\emph{First case: assume that $\eta/(K \norm{Dh(a_j)v_j}) \geq
\eta_0$.} The image of the ball $B(a_j,\rho)$ contains the ball
$B(ha_j, \rho \norm{Dh(a_j)v_j}/K)$, which itself contains $
B(ha_j, C \eta)$ since $h$ is good. Moreover, for $x,x'\in Y$
holds $d(hx,hx') \leq d(x,x') K \norm{Dh(a_j)v_j} \leq \diam Y
\frac{\eta}{\eta_0}$. In particular, if $M \geq \diam Y/ (C
\eta_0)$, we get $h(Y) \subset B(ha_j, MC \eta)$. We can thus
take a ball $B(ha_j, C \eta)$ and a set $A_1=h(Y)$. To
conclude, we should check that $\mu( B(x',\eta)) \geq M^{-1}
\mu(A_i)$ for any $x' \in B(ha_j, (C-1) \eta)$, if $M$ is large
enough. Since the iterates of $\TY$ have a uniformly bounded
distortion,
  \begin{equation}
  \frac{\mu(B(x', \eta))}{\mu(A_i)} \asymp \frac{\mu(h^{-1}
  B(x',\eta))}{\mu(Y)}.
  \end{equation}
Moreover, $h^{-1} B(x',\eta)$ contains $B(h^{-1}x', \eta/ (K
\norm{Dh(a_j)v_j}))$, which itself contains $B(h^{-1}x',
\eta_0)$. The measure of these balls is uniformly bounded from
below. This concludes the proof in this case.

\emph{Second case: assume now that $\eta/ (K \norm{Dh(a_j)v_j})
\leq \eta_0$.} Let $\eta_h= \eta/ (K \norm{Dh(a_j)v_j})$, it is
bounded by $\eta_0$. Hence, the assumption (3) gives sets
$A_1,\dots, A_k$, balls $B(x_1, \tC \eta_h), \dots, B(x_k, \tC
\eta_h)$ and inverse branches $h_1,\dots, h_\ell$ defined
respectively on $Y_{j_1},\dots, Y_{j_\ell} $. We will show that
the balls $B(hx_1, C \eta),\dots, B(hx_k, C\eta)$, the sets
$\bar A_i=h(A_i)$ and the inverse branches $h\circ h_1,\dots,
h\circ h_\ell$ satisfy the conclusion of the fact.

Let us first show that the inverse branch $h\circ h_i$ is good. By
definition of $h_i$, $\norm{Dh_i(a_{j_i}) v_{j_i}} \geq \tC \eta_h/K
\geq K^2 C \eta/ (\rho \norm{Dh(a_j)v_j})$. We have $D(h\circ
h_i)(a_{j_i})v_{j_i} = Dh( h_i a_{j_i}) Dh_i(a_{j_i})v_{j_i}$.
Moreover, $\norm{Dh(x)v} \geq K^{-1}\norm{v} \norm{Dh(a_j) v_j}$.
Therefore,
  \begin{equation*}
  \norm{ D(h\circ h_i)(a_{j_i}) v_{j_i}} \geq K^{-1} \norm{Dh_i(a_{j_i})v_{j_i}}
  \norm{Dh(a_j)v_j}
  \geq K^{-1}  \frac{K^2 C \eta}{\rho \norm{Dh(a_j)v_j}}
  \norm{Dh(a_j) v_j}
  = K C \eta/\rho.
  \end{equation*}
This shows that $h\circ h_i$ is good.

The set $h B(x_i, \tC \eta_h)$ contains the ball $B(h x_i, \tC
\eta_h \norm{Dh(a_j)v_j}/K)$, which itself contains the ball $B(h
x_i, C \eta)$ because $\tC \geq K^2 C$. Moreover, for any $x' \in
B(h x_i, (C-1)\eta)$, the set $h^{-1} B(x',\eta)$ contains the ball
$B(h^{-1}x', \eta / (K \norm{Dh(a_j)v_j})) = B(h^{-1}x',\eta_h)$. As
the distortion of the iterates of $\TY$ is uniformly bounded, we
obtain for any $x' \in B(h x_i, (C-1)\eta)$
  \begin{equation}
  \frac{\mu (B(x', \eta))}{\mu (\bar A_i)} \asymp
  \frac{ \mu (h^{-1} B(x',\eta))}{\mu (A_i)}
  \geq \frac{ \mu( B(hx',\eta_h))}{\mu (A_i)}
  \geq \tD^{-1}.
  \end{equation}
Finally, as $A_i \subset B(x_i, \tD \tC \eta_h)$, $\bar A_i$ is
contained in $B(hx_i, \tD \tC \eta_h K \norm{Dh(a_j)v_j})= B(hx_i,
\tD \tC \eta)$.
\end{proof}

The previous criterion easily implies that Gibbs measures in
dimension 1 have the uniform weak Federer property:
\begin{prop}
\label{prop:DemontreFedererGibbs} Let $T$ be a $C^2$ uniformly
expanding map on the circle $\Sbb^1$, and let $\mu$ be a Gibbs
measure corresponding to a $C^1$ potential. Then there exists a
subset $Y$ of $\Sbb^1$ such that $T$ is nonuniformly expanding
with base $Y$, for the measure $\mu$.
\end{prop}
\begin{proof}
Let $d$ be the topological degree of $T$, and let $x_0$ be a fixed
point of $T$. Let $Y=Z=\Sbb^1 \moins\{x_0\}$. Then $\Sbb^1 \moins
T^{-1}(x_0)$ it the union of $d$ intervals $W_1,\dots, W_j$, each of
them being sent by $T$ onto $Z$. These intervals form a partition
(modulo $0$) of $Z$ satisfying the first four points of Definition
\ref{def:NonUnifDilatant} (for $r_i=1$, $1\leq i\leq d$). If we can
prove that $T$ satisfies the assumptions of the previous
proposition, the proof will be complete. The assumptions (1) and (2)
are clear, the fourth is equivalent to the bounded distortion for
Lebesgue measure since we are in one dimension. Let us check (3),
for some $\tC>0$. Let $\eta_0$ be small enough so that, for any
$x\in Z$ and any inverse branch $h\in \HH$, $|h'(x)|\geq \tC
\eta_0$. We take no ball $B(x_i,\tC \eta)$, no set $A_i$, and all
the inverse branches $h\in \HH$. Then \eqref{balls_included} is
empty, hence trivial, and \eqref{covers_all} is also trivial.
\end{proof}

\subsection{Farey sequences}

Let $r>1$. Let $T$ be the map on $X=[0,1]$ given by
\eqref{def:TFarey}, and let $\TT$ be its extension to
$[0,1]\times \R/(\log r)\Z$ defined in \eqref{def:TTFarey},
using a function $\phi$. This function is not $C^1$ on $[0,1]$,
which seems to be a problem since we always worked with a
function $\phi$ of class $C^1$. To avoid this problem, we can
simply work with the disjoint union $X=[0,1/2]\sqcup [1/2,1]$,
on which $\phi$ is $C^1$. All our results in the previous
sections have been formulated for transformations on $X\times
\R/2\pi\Z$, but the same results hold verbatim on $X\times
\R/\gamma\Z$ for any $\gamma\not=0$, and in particular for
$\gamma=\log r$. Henceforth, we will simply denote $\R/(\log
r)\Z$ by $\Sbb^1$ and apply without further notice the
preceding results.

Let $x_0=1/2$, and set $x_n=h_A(x_{n-1})$, i.e., $x_n$ is the
preimage of $x_{n-1}$ under the left branch of $T$. Explicitly,
$x_n=1/(n+2)$. Let $I_j=(x_{j},x_{j-1})$. Let also $\bar
I_j=1-I_j$ be the symmetric of $I_j$ with respect to $1/2$. Let
$Y=(x_1,x_0)=(1/3,1/2)$, and denote by $T_Y$ the map induced by
$T$ on $Y$. Its combinatorics can be described as follows: a
point of $Y$ is sent by $T$ in $(1/2,1)$, it spends some time
$i>0$ there, is then sent back to $(0,1/2)$, and increases (for
$j\geq 0$ iterates) before entering back in $Y$. The points
with this combinatorics form an interval $I_{i,j}:=T^{-1}(\bar
I_i) \cap T^{-i-1}(I_{j+1})$, and $T^{i+j+1}(I_{i,j})=Y$.
Letting $r_{i,j}=i+j+1$, we thus obtain a partition of $Y$ that
satisfies the first point of Definition
\ref{def:NonUnifDilatant}.

\begin{prop}
The map $T$ is nonuniformly expanding of base $Y$, in the sense
of Definition \ref{def:NonUnifDilatant}, for the partition
$\{I_{i,j}\}_{i>0, j\geq 0}$ and Minkowski's measure $\mu$.
Moreover, it is mixing.
\end{prop}
\begin{proof}
The first point of Definition \ref{def:NonUnifDilatant} is clear.
For the second one, note that the jacobian of $T$ for Minkowski's
measure is everywhere equal to $2$ by definition. Hence, the
jacobian of $T_Y$ on $I_{i,j}$ is constant (equal to $2^{i+j+1}$),
and $D((\log J)\circ h_{i,j})=0$. The third point is trivial. For
the fourth one, we have for any $\sigma>0$
  \begin{equation}
  \int_Y e^{\sigma r} = \sum \mu(I_{i,j}) e^{\sigma (i+j+1)}
  =\sum 2^{-i-j-3} e^{\sigma (i+j+1)},
  \end{equation}
which is finite as soon as $\sigma<\log 2$. The mixing of $T$ is a
consequence of the equality $\gcd\{r_{i,j}\}=1$.

Thus, we just have to prove the uniform weak Federer property.
To do this, we will use Proposition \ref{prop:CritereFederer}.
Let $Y_0=Y$, and let $Y_1$ be its symmetric with respect to
$1/2$. Let $Z=Y_0\cup Y_1$, and let $\TY$ be the first return
map induced by $T$ on $Z$. It sends each interval $T^{-1}(\bar
I_i)\cap Y_0$ bijectively to $Y_1$, and each interval
$T^{-1}(I_i)\cap Y_1$ bijectively to $Y_0$. If we prove that
$\TY$ satisfies the assumptions of Proposition
\ref{prop:CritereFederer}, this will conclude the proof of the
uniform weak Federer property, since the inverse branches of
the iterates of $T_Y$ are in particular inverse branches of
iterates of $\TY$.

Assumptions (1) and (2) of Proposition
\ref{prop:CritereFederer} are trivial (since $J$ is constant on
each monotonicity interval of $\TY$). For the fourth point, the
quickest argument is certainly to use the fact that all the
inverse branches of the iterates of $\TY$ are homographies
(hence with vanishing schwarzian derivative) which can be
extended to the whole interval $[0,1]$. Koebe's Lemma
\cite[Theorem IV.1.2]{demelo_vanstrien} directly yields the
uniform quasi-conformality.

Hence, we just have to check point (3). It is sufficient to
check it on $Y_0$, since everything is symmetric with respect
to $1/2$. If $J$ is an interval, we will denote its length by
$|J|$. Then $|\bar I_n|$ is a decreasing sequence, with $|\bar
I_{n+1}|/|\bar I_n|\to 1$ when $n\to\infty$, since $T'(1)=1$.
As a consequence, $K_n=T^{-1}(\bar I_n)\cap Y_0$ satisfies
$|K_{n+1}|/|K_n|\to 1$, and there exists $C>0$ such that
$|K_m|\leq C |K_n|$ for all $m\geq n$. Finally, $\mu(K_n)=
2^{-n-2}$.

We will use the following fact: \emph{for any $C>0$, there exists
$D>0$ such that, for any interval $J$ included in an interval $K_n$
with $|J|\geq C^{-1}|K_n|$, then $\mu(J) \geq D^{-1}\mu(K_n)$}. To
prove this fact, we apply once the map $\TY$, which sends $K_n$ to
$Y_1$, and $J$ to an interval $J'$ satisfying $|J'|\geq C^{-1}K^{-1}
|Y_1|$ by quasi conformality. Hence, $\mu(J')$ is uniformly bounded
from below. As $\mu(J')/\mu(Y_1)=\mu(J)/\mu(K_n)$, this proves the
fact.

We can now prove the third assumption of Proposition
\ref{prop:CritereFederer}, on $Y_0$. Let $\tC>1$. We will construct
inverse branches $h_1,\dots,h_\ell$, balls $B(x_1,\tC \eta),\dots,
B(x_k,\tC\eta)$ and sets $A_1,\dots, A_k$ as follows, if $\eta$ is
small enough.

Let $N$ be maximal such that $|K_n| \geq \tC \eta$ for $n\leq
N$. We take $\ell=N$, and let $h_1,\dots,h_\ell$ be the inverse
branches of $\TY$ whose images are the intervals $K_1,\dots,
K_\ell$. Then $h_i$ is defined on $Y_1$, of length $1/6$, and
the length of its image $K_i$ is $\geq \tC\eta$. Hence, there
exists a point $y_i\in Y_1$ with $h_i'(y_i)\geq 6\tC\eta$. This
proves \eqref{eq:TaillehiGrande}.

We decompose the remaining interval as a union of intervals of
length $2\tC\eta$, excepted maybe the first one whose length
belongs to $[2\tC\eta, 4\tC \eta)$. Let us denote this
decomposition by $J_0,\dots,J_p$. Since $|K_N|=o( \sum_{n>N}
|K_n|)$ when $N\to\infty$, we have $p\geq 2$ if $\eta$ is small
enough. Let us define sets $A_1,\dots,A_{p}$ by $A_i=J_i$ for
$i>1$, and $A_1=J_0\cup J_1$. Let $B(x_i, \tC \eta)=J_{i-1}$
for $i>1$, and let $B(x_1, \tC \eta)$ be the leftmost part of
$J_0$. For $i>1$, the ball $B(x_i, \tC\eta)$ is \emph{not}
included in the set $A_i$, it is strictly to its left. The
balls are disjoint, and $A_i \subset B(x_i, 5\tC \eta)$. Let us
show that they satisfy the desired conclusion: we have to prove
that, for any interval $J$ of length $2\eta$ included in
$B(x_i, \tC \eta)$, then $\mu(J) \geq \tD^{-1} \mu(A_i)$ holds
for some constant $\tD$ (independent of $\eta$). Either $J$
contains an interval $K_n$, or it intersects such an interval
along a subinterval of length at least $\eta$. Moreover,
$|K_n|\leq C |K_{N+1}|\leq C \tC\eta$. In both cases, the fact
we proved above implies that $\mu(J) \geq D^{-1} \mu(K_n)$.

We first deal with $i=1$. As $|K_{n+1}|\sim |K_n|$, the set
$A_1$ is covered by $\bigcup_{k=1}^7 K_{N+k}$ if $N$ is large
enough (hence, if $\eta$ is small enough). These $7$ intervals
have comparable measures since $\mu(K_m)=2^{-m-2}$, hence
$\mu(A_{1})\leq C \mu(K_{N+k})$ for $1\leq k\leq 7$. As
$\mu(J)\geq D^{-1}\mu(K_n)$ for at least one these $K_n$'s, we
indeed conclude $\mu(J)\geq C^{-1}\mu(A_1)$.

Assume now $i>1$. There exists an interval $K_n$ intersecting $J$
with $\mu(J) \geq C^{-1}\mu(K_n)$. Since $A_i$ is located to the
right of $K_n$, we get
  \begin{equation}
  \mu(A_i)\leq C \sum_{m=n}^\infty \mu(K_m) = C\sum_{m=n}^\infty
  2^{-m-2} \leq C 2^{-n-2} \leq C \mu(K_n).
  \end{equation}
This also concludes the proof in this case.
\end{proof}

\begin{lem}
The function $\phi$ is not cohomologous to a locally constant
function.
\end{lem}
\begin{proof}
Assume by contradiction that there exists a $C^1$ function $f$ such
that $\phi_Y-f+f\circ T_Y$ is constant on each interval $I_{i,j}$,
equal to some number $a_{i,j}$. The interval $I_{1,1}$ contains the
point $x=3/2-\sqrt{5}/2$, with $T_Y(x)=x$. Necessarily,
$a_{1,1}=\phi_Y(x)$. In the same way, the interval $I_{2,1}$
contains $x'=1-\sqrt{3}/3$, invariant under $T_Y$, which gives
$a_{2,1}=\phi_Y(x')$.

Let now $y=1-\sqrt{6}/4$. This point belongs to $I_{1,1}$, but
$T_Y(y)\in I_{2,1}$, and $T_Y^2(y)=y$. Then
  \begin{equation}
  \phi_Y(y)+\phi_Y(T_Yy)= a_{1,1}+a_{2,1}=\phi_Y(x)+\phi_Y(x').
  \end{equation}
However, it is possible to compute explicitly
$\phi_Y(y)+\phi_Y(T_Yy)- \phi_Y(x)-\phi_Y(x')$, and check that this
quantity is nonzero (approximately equal to $-0.013$). This is a
contradiction.
\end{proof}

The previous proposition and lemma show that the results of
Paragraph \ref{subsec:ResultatsLimites} apply to $\TT$.
However, this is not sufficient to prove Theorems
\ref{thm:MelangeExponentielFarey} and \ref{thm:ThmLimiteFarey},
since these results are pointwise while the results of
Paragraph \ref{subsec:ResultatsLimites} are averaged. We will
therefore need an additional ingredient. Let $X^{(n)}$ be the
extension of $X$ defined in Paragraph \ref{subsec:modele}, and
let $\pi^{(n)}$, $\tilde\pi^{(n)}$ be the corresponding
projections.
\begin{lem}
\label{lem:DecritRelevement} For any $n\in \N$, there exists a
constant $C(n)$ such that, for any integrable function
$u:X\times \Sbb^1 \to \C$, for almost all $(x,\omega)\in
X\times\Sbb^1$ and for any $k\in \N$,
  \begin{equation}
  \hat \TT^k u(x,\omega)= C(n)\sum_{\pi^{(n)}(x')=x} 2^{-h(x')}
  \hat \UU^k (u\circ \tilde\pi^{(n)})(x',\omega).
  \end{equation}
\end{lem}
\begin{proof}
Let $\BB$ be the $\sigma$-algebra of Borel measurable subsets
of $X\times \Sbb^1$, and let
$\BB'=(\tilde\pi^{(n)})^{-1}(\BB)$. This is a
sub-$\sigma$-algebra of the Borel $\sigma$-algebra on
$X^{(n)}\times \Sbb^1$. A function $v$ on $X^{(n)}\times
\Sbb^1$ can be written as $u\circ \tilde\pi^{(n)}$ if and only
if $v$ is $\BB'$-measurable.

Let us first prove that
  \begin{equation}
  \label{eq:ReleveTransfert}
  (\hat \TT^k u)\circ \tilde \pi^{(n)}= E( \hat \UU^k (u\circ \tilde \pi^{(n)}) \given
  \BB').
  \end{equation}
To do this, let us write $E( \hat \UU^k (u \circ \tilde
\pi^{(n)}) \given \BB')=v\circ \tilde\pi^{(n)}$. As $\tilde
\mu\otimes\Leb=\tilde\pi^{(n)}_*(\tilde \mu^{(n)}\otimes
\Leb)$, we have for any measurable function $f$ on $X\times
\Sbb^1$
  \begin{equation}
  \int_{X\times \Sbb^1} vf=\int_{X^{(n)}\times \Sbb^1} v\circ
  \tilde\pi^{(n)} f\circ \tilde\pi^{(n)}
  = \int_{X^{(n)}\times \Sbb^1} E( \hat \UU^k (u\circ \tilde \pi^{(n)}) \given
  \BB') f\circ \tilde\pi^{(n)}.
  \end{equation}
As $f\circ \tilde\pi^{(n)}$ is $\BB'$-measurable, we get
  \begin{align*}
  \int_{X\times \Sbb^1} vf
  &= \int_{X^{(n)}\times \Sbb^1} \hat \UU^k (u\circ \tilde \pi^{(n)})
  f\circ \tilde\pi^{(n)}
  = \int_{X^{(n)}\times \Sbb^1} u\circ \tilde \pi^{(n)} f\circ
  \tilde\pi^{(n)} \circ \UU^k
  \\&
  = \int_{X^{(n)}\times \Sbb^1} u\circ \tilde \pi^{(n)}
  f\circ\TT^k\circ \tilde\pi^{(n)}
  =\int_{X\times \Sbb^1} u f\circ \TT^k.
  \end{align*}
This last equality shows that $v=\hat\TT^k u$, and concludes the
proof of \eqref{eq:ReleveTransfert}.

The set $X^{(n)}$ is endowed with a countable partition $\AAA$
such that $\pi^{(n)}$ is injective on each element of the
partition. Let us define a function $F$ on $X^{(n)}$ as
follows: on each set $a\in \AAA$, let $F=\dd \tilde \mu^{(n)}/
\dd (\tilde\mu\circ \pi^{(n)}_{|a})$. This is the local
Radon-Nikodym derivative of $\tilde\mu^{(n)}$ with respect to
$(\pi^{(n)})^* \tilde \mu$. As
$\pi^{(n)}_*\tilde\mu^{(n)}=\tilde\mu$, we have
$\sum_{\pi^{(n)}(x')=x} F(x')=1$ for almost every $x\in X$. Let
us show that the conditional expectation with respect to $\BB'$
is given by
  \begin{equation}
  \label{eq:FormuleEsperanceConditionnelle}
  E( v\given \BB')(x,\omega)= \sum_{ \pi^{(n)}(x')=\pi^{(n)}(x)}F(x') v(x',\omega).
  \end{equation}
Let us indeed define a function $w$ on $X\times \Sbb^1$ by
  \begin{equation}
  w(x,\omega)=\sum_{ \pi^{(n)}(x')=x} F(x')v(x',\omega)
  =\sum_{a\in \AAA} 1_{x\in \pi^{(n)}a} F( (\pi^{(n)}_{|a})^{-1}x)
  v((\pi^{(n)}_{|a})^{-1}x, \omega).
  \end{equation}
If $f$ is a measurable function on $X\times \Sbb^1$,
  \begin{align*}
  \int_{X\times \Sbb^1} fw&= \sum_{a\in \AAA} \int_{\pi^{(n)}(a)}
  f(x,\omega) F( (\pi^{(n)}_{|a})^{-1}x)
  v((\pi^{(n)}_{|a})^{-1}x, \omega)\dd\tilde\mu(x)\dd\omega
  \\&
  = \sum_{a\in \AAA} \int_{a} f(\pi^{(n)}x',\omega)
  v(x',\omega) \dd\tilde\mu^{(n)}(x')\dd\omega
  = \int_{X^{(n)}\times \Sbb^1} f\circ \tilde\pi^{(n)} v.
  \end{align*}
This proves \eqref{eq:FormuleEsperanceConditionnelle}. Together with
\eqref{eq:ReleveTransfert}, this implies the lemma if we can prove
that
  \begin{equation}
  \label{eq:DonneFormuleF}
  F(x')=C(n)2^{-h(x')}.
  \end{equation}
As $T_Y$ is the first return map to $Y$, the jacobian of $\pi^{(1)}$
for the measure $\tilde\mu^{(1)}$ on $Y$ is equal to $1$. Since
$\tilde\mu^{(n)}$ is proportional to $\tilde\mu^{(1)}$ on $Y$, this
implies that $F$ is constant on $Y$, equal to a constant $C(n)$.
This proves \eqref{eq:DonneFormuleF} for points with zero height.

The jacobian of $T$ for $\tilde\mu$ is equal to $2$, while the
jacobian of $U$ is equal to $1$ on the set of points that do not
come back to the basis. By induction over $h(x')$, this implies
\eqref{eq:DonneFormuleF}.
\end{proof}

\begin{cor}
\label{cor:OKFareyTransfert} There exist constants $C>0$ and
$\bar\theta<1$ such that, for any $C^6$ function $f:X\times
\Sbb^1 \to\C$, for any $(x,\omega)\in X\times\Sbb^1$,
  \begin{equation}
  \left| \hat\TT^n f(x,\omega)-\int f \right|\leq C \bar\theta^n
  \norm{f}_{C^6}.
  \end{equation}
\end{cor}
\begin{proof}
Since everything is symmetric with respect to $1/2$, and
continuous, it is sufficient to prove the assertion for almost
every $x\in (1/2,1)$.

We work in $X^{(N)}$, where $N$ is given by Theorem
\ref{thm:MainContraction}. Note that $d^{(N)}$ is equal to $1$,
since $r^{(N)}$ takes the values $2N$ and $2N+1$. Applying
Theorem \ref{thm:MelangeExponentielPrecis} to the function
$v=f\circ \tilde\pi^{(N)}$, we get: for any $n\in \N$, for any
$x'\in X^{(N)} $ with $h(x')\leq n/2$,
  \begin{equation}
  \left| \hat\UU^n (f\circ \tilde\pi^{(N)})(x',\omega) -\int f\right|
  \leq C \bar\theta^n \norm{f}_{C^6}.
  \end{equation}
Together with Lemma \ref{lem:DecritRelevement}, this yields
  \begin{equation*}
  \left| \hat\TT^n f(x,\omega)-\int f\right| \leq C \left(\sum_{\pi^{(N)}(x')=x,
  h(x')\leq n/2} \bar\theta^n 2^{-h(x')} + \sum_{\pi^{(N)}(x')=x, h(x')>n/2}
  2^{-h(x')}\right) \norm{f}_{C^6}.
  \end{equation*}
To conclude, it is thus sufficient to prove that, for $x\in
(1/2,1)$, the cardinality of
  \begin{equation}
  \label{eq:qsdiufiop}
  \{ x'\st \pi^{(N)}(x')=x,\ h(x')=k\}
  \end{equation}
grows at most polynomially with $k$. If we write a point of
$X^{(N)}$ as a pair $(x',j)$ with $x'\in Y$ and
$j<r^{(N)}(x')$, it is easy to check that $U^k$ induces a
bijection between the set \eqref{eq:qsdiufiop} and the set of
points in $T^{-k}(x)\cap Y$ whose first $k$ iterates under $T$
spend a time $t<N$ in $Y$. If $t$ is fixed, such a point is
determined by the combinatorics $(i_1,j_1,\dots,
i_{t},j_{t},i_{t+1})$ of times spent in $[1/2,1]$, then in
$[0,1/2]$, then in $[1/2,1]$, and so on, with the constraint
that the sum of these lengths is $k$ (we recall that we assume
$x\in (1/2,1)$). As a consequence,
  \begin{equation}
  \Card\{ x'\st \pi^{(N)}(x')=x,\ h(x')=k\}
  \leq \sum_{t=0}^{N-1} k^{2t+1} \leq C k^{2N}.
  \end{equation}
This quantity indeed grows polynomially.
\end{proof}

\begin{proof}[Proof of Theorem
\ref{thm:MelangeExponentielFarey}] If $f$ is a continuous
function on $[0,1]\times \Sbb^1$, then $\int f \dd\bar\mu_n=
\hat\TT^n f(1,0)$. Hence, Corollary \ref{cor:OKFareyTransfert}
shows the theorem for $C^6$ functions. The case of $C^\alpha$
functions is then deduced by interpolation, just like at the
end of the proof of Theorem \ref{thm:MelangeExponentiel}.
\end{proof}

\begin{proof}[Proof of Theorem \ref{thm:ThmLimiteFarey}]
If $\psi$ is a $C^6$ function which is not a coboundary, we
show like in the proof of Corollary \ref{cor:OKFareyTransfert}
(but using Theorem \ref{thm:MelangeExponentielPrecisLocal}
instead of Theorem \ref{thm:MelangeExponentielPrecis}) that,
for $|t|\leq \tau_0$,
  \begin{equation}
  \left| \hat\TT_t^n f(x,\omega)
  -\left(1-\frac{\sigma^2t^2}{2}\right)^n \int f \right| \leq
  C(\bar\theta^n+ |t|(1-ct^2)^n) \norm{f}_{C^6}.
  \end{equation}
Moreover, if $\psi$ is aperiodic, for $\tau_0\leq |t|\leq t_0$,
  \begin{equation}
  \left| \hat\TT_t^n f(x,\omega) \right| \leq
  C\bar\theta^n \norm{f}_{C^6}.
  \end{equation}
As $\hat\TT_t^n 1 (1,0)=E(e^{it\sum_{k=1}^n\psi(X_k)})$, this
implies the limit assertions in Theorem
\ref{thm:ThmLimiteFarey}.

The automatic regularity properties still have to be checked. If
$\psi=f-f\circ \TT$ with $f$ measurable, let us show that $f$ is
continuous on $[0,1]$. Proposition \ref{prop:CaracteriseSigma2}
shows that $f$ is continuous on $Y\times \Sbb^1$. As $\TT$ is an
homeomorphism between $Y\times \Sbb^1$ and $[1/2,1]\times \Sbb^1$,
we conclude from the equality $f\circ \TT=f-\psi$ that $f$ is
continuous on $[1/2,1]\times \Sbb^1$. Finally, as $\TT$ is an
homeomorphism between $[1/2,1]\times \Sbb^1$ and $[0,1]\times
\Sbb^1$, we obtain with the same argument the continuity of $f$ on
the whole space.

We argue in the same way for the cohomological equation in
$\R/\lambda\Z$, by using Proposition
\ref{prop:CaracteriseAperiodique}.
\end{proof}

\appendix

\section{Contraction properties of transfer operators}
\label{app:app} In this appendix, we prove Theorem
\ref{thm:MainContraction} on the contraction properties (in
$C^1$ norm or in Dolgopyat norm) of the transfer operator
associated to a map $T_Y$, where $T$ is a nonuniformly
expanding map of base $Y$. Henceforth, the notations and
assumptions will be those of Theorem \ref{thm:MainContraction}.

\subsection{Contraction in the $C^1$ norm}
\label{app:ContractClassique}

In this paragraph, we introduce the tools to prove the first
part of Theorem \ref{thm:MainContraction}. However, the choice
of the constants $N$ and $\theta$ of Theorem
\ref{thm:MainContraction} will only be possible at the complete
end of the proof, in the next paragraph.

We will use several times the following distortion lemma, whose
proof is completely standard and will be omitted.
\begin{lem}
\label{lem:ControleJacobien} Let $J^{(n)}(x)$ be the inverse of the
jacobian of $T_Y^n$ at the point $x$. There exists $C>0$
(independent of $n$) such that, for any $h\in \HH_n$, for any
$x,y\in Y$, $\norm{D(J^{(n)}\circ h)(x)}\leq C J^{(n)}\circ h(x)$
and $J^{(n)}\circ h(x)\leq C J^{(n)}\circ h(y)$.
\end{lem}

For small enough $\epsilon$, we define an operator
$\Lp_\epsilon$ acting on functions from $Y$ to $\C$, by
$\Lp_\epsilon u(x)=\sum J(hx)u(hx)e^{\epsilon r(x)}$. If
$H_0\subset \HH$, we will also denote by $\Lp_{\epsilon,H_0}$
the same operator but where the sum is restricted to the
inverse branches belonging to $H_0$. The following elementary
estimates will be used again and again in all the forthcoming
arguments.
\begin{lem}
\label{lem:ContractDebile} There exists a function
$\alpha(\epsilon)$ which tends to $0$ when $\epsilon\to 0$ such
that $\norm{\Lp_{\epsilon}}_{L^2\to L^2} \leq
e^{\alpha(\epsilon)}$ and $\norm{\Lp_\epsilon}_{C^0\to C^0}
\leq e^{\alpha(\epsilon)}$.

Moreover, if $\epsilon_0>0$ is small enough, for any
$\gamma>0$, there exists $H_0\subset \HH$ with a finite
complement such that $\norm{\Lp_{\epsilon_0,H_0}}_{L^2 \to
L^2}\leq \gamma$.
\end{lem}
\begin{proof}
We have
  \begin{equation*}
  (\Lp_{\epsilon,H_0}u(x))^2 = \left( \sum_{h\in H_0} J(hx) u(hx)e^{\epsilon
  r(hx)}\right)^2 \leq \left( \sum_{h\in H_0} J(hx) u(hx)^2 \right)
  \left(\sum_{h\in H_0}
  J(hx) e^{2\epsilon r(hx)}\right).
  \end{equation*}
Consequently, $\norm{\Lp_{\epsilon,H_0}u}_{L^2} \leq
\norm{u}_{L^2}\cdot \sup_{x\in Y}\left(\sum_{h\in H_0} J(hx)
e^{2\epsilon r(hx)}\right)^{1/2}$. We have $J(hx)\leq C J(hy)$
for any $h\in \HH$ and all $x,y\in Y$, hence $\sum
J(hx)e^{2\epsilon r(hx)} \leq C \sum J(hy)e^{2\epsilon r(hy)}$.
Integrating this inequality with respect to $y$, we get
  \begin{equation}
  \sum_{h\in H_0} J(hx) e^{2\epsilon r(hx)}  \leq C \sum_{h\in H_0}
  \int_Y J(hy)e^{2\epsilon r(hy)}\dd\mu_Y(y)=C \int_{H_0(Y)} e^{2\epsilon
  r(y)}\dd\mu_Y(y).
  \end{equation}
This quantity is finite if $\epsilon$ is small enough, by the fourth
assumption of Definition \ref{def:NonUnifDilatant}. Taking the
complement of $H_0$ small enough, it can even be made arbitrarily
small. This proves the second point of the lemma.

For the first point, we have to be slightly more precise. For
any $x$, we have $e^{2\epsilon r(hx)}\leq 1+2\epsilon r(hx)
e^{2\epsilon r(hx)}$. Hence, using the inequality $J(hx)\leq C
J(hy)$ for any $h\in \HH$ and $x,y\in Y$, we get
  \begin{equation*}
  \sum_{h\in \HH} J(hx)
  e^{2\epsilon r(hx)} \leq \sum_{h\in \HH} J(hx) + 2\epsilon
  \sum_{h\in \HH} J(hx) r(hx) e^{2\epsilon r(hx)}
  \leq 1 + C \epsilon \sum_{h\in \HH} J(hy) r(hy) e^{2\epsilon
  r(hy)}.
  \end{equation*}
Integrating with respect to $y$,
  \begin{equation}
  \sum_{h\in \HH} J(hx)
  e^{2\epsilon r(hx)} \leq 1+ C \epsilon \int_Y r(y)e^{2\epsilon r(y)} \dd\mu_Y(y),
  \end{equation}
and this last integral is uniformly bounded if $\epsilon$ is small
enough. This gives the desired estimate for the action of
$\Lp_\epsilon$ on $L^2$ and $C^0$.
\end{proof}

Let us prove a lemma which will easily imply
\eqref{eq:ContracteTriviale}.
\begin{lem}
\label{lem:contractSKJFDL} There exist $\epsilon_0>0$ and
$\theta_0<1$ such that, for any $A>0$, $n\in \N$ and
$\epsilon<\epsilon_0$, there exists $C>0$ such that, for any
$\psi\in \CC^{A,\epsilon}_{n}$ and $v\in C^1(Y)$,
  \begin{equation}
  \norm{ \Lp^{n} (\psi v)}_{C^1}\leq \theta_0^{n} \left(
  \sup_{x\in Y} |\psi(x)|/ e^{\epsilon r^{(n)}(x)} \right)
  \norm{v}_{C^1} + C \norm{\psi}_{\CC^{A,\epsilon}_{n}}
  \norm{v}_{C^0}.
  \end{equation}
\end{lem}
\begin{proof}
First, since $|\psi(x)| \leq \norm{\psi}_{\CC^{A,\epsilon}_{n}}
e^{\epsilon r^{(n)}(x)}$, we have
  \begin{equation}
  \norm{\Lp^n(\psi v)}_{C^0} \leq \norm{\psi}_{\CC^{A,\epsilon}_{n}}
  \norm{\Lp_{\epsilon}^n |v|}_{C^0}
  \leq \norm{\psi}_{\CC^{A,\epsilon}_{n}} e^{n\alpha(\epsilon)}
  \norm{v}_{C^0},
  \end{equation}
by Lemma \ref{lem:ContractLepsilon}. This gives the desired control
in the $C^0$ norm. For the $C^1$ norm, we differentiate $\Lp^n(\psi
v)=\sum_{h\in \HH_n} J^{(n)}(hx) \psi(hx) v(hx)$. If we
differentiate $J^{(n)}(hx)$, we use the estimate
$\norm{D(J^{(n)}\circ h)(x)} \leq C J^{(n)}(hx)$ given by Lemma
\ref{lem:ControleJacobien}, and get the same bound as for the $C^0$
norm. If we differentiate $\psi(hx)$, its derivative is bounded by
$A \norm{\psi}_{\CC^{A,\epsilon}_{n}}e^{\epsilon r^{(n)}(hx)}$, and
using the same argument as for the $C^0$ norm we obtain the same
bound (with an additional factor $A$, which is not a problem since
$C$ is allowed to depend on $A$ in the statement of the lemma).

Finally, if we differentiate $v\circ h$, we have $\norm{D(v\circ
h)(x)}\leq \kappa^{-n} \norm{Dv(hx)}$, and we therefore get a bound
  \begin{align*}
  \kappa^{-n} \norm{Dv}_{C^0} \Lp^n |\psi|
  &\leq \kappa^{-n}\norm{v}_{C^1} \left(
  \sup_{x\in Y} |\psi(x)|/ e^{\epsilon r^{(n)}(x)} \right)
  \Lp^n(e^{\epsilon r^{(n)}})
  \\&
  \leq \kappa^{-n}\norm{v}_{C^1} \left(
  \sup_{x\in Y} |\psi(x)|/ e^{\epsilon r^{(n)}(x)} \right)
  e^{n\alpha(\epsilon)}.
  \end{align*}
If $\epsilon$ is small enough, $\kappa^{-1}
e^{\alpha(\epsilon)}<1$. This concludes the proof.
\end{proof}

We now turn to the proof of \eqref{eq:ContractC1Iteres}. As a
preliminary estimate, let us first consider the case
$\psi_i=e^{\epsilon r^{(N)}}$ for all $i$, in the following
lemma.

\begin{lem}
\label{lem:ContractLepsilon} There exist $N_0>0$, $\theta_0<1$,
$C>0$, $\epsilon_0>0$ and a function $\alpha:(0,\epsilon_0)\to
\R_+$ tending to $0$ when $\epsilon\to 0$, satisfying the
following property. For any $N\geq N_0$ and $\epsilon<
\epsilon_0$, for any $C^1$ function $v:Y\to \C$,
  \begin{equation}
  \norm{D(\Lp_\epsilon^{N} v)}_{C^0} \leq
  \theta_0^{N}\norm{Dv}_{C^0} + Ce^{N\alpha(\epsilon)} \norm{v}_{L^2}.
  \end{equation}
\end{lem}
\begin{proof}
We have $\Lp_\epsilon^N v= \sum_{h\in \HH_N} J^{(N)}(hx)
e^{\epsilon r^{(N)}(hx)} v(hx)$. By Lemma
\ref{lem:ControleJacobien}, $J^{(N)}(hx)\leq C J^{(N)}(hy)$,
and $\norm{D(J^{(N)}\circ h)(x)}\leq C J^{(N)}(hx)$. Moreover,
since $h$ contracts the distances by at least $\kappa^N$,
$|v(hx)|\leq |v(hy)|+ C \kappa^{-N} \norm{Dv}$. Hence,
  \begin{equation*}
  J^{(N)}(hx) e^{\epsilon r^{(N)}(hx)}|v(hx)| \leq C J^{(N)}(hy)
  e^{\epsilon r^{(N)}(hy)} |v(hy)| + C \kappa^{-N} J^{(N)}(hy)
  e^{\epsilon r^{(N)}(hy)} \norm{Dv}_{C^0}.
  \end{equation*}
Integrating this equation over $y$ and summing over the inverse
branches, we conclude
  \begin{equation}
  \Lp_\epsilon^N |v|(x)\leq C \int e^{\epsilon r^{(N)}} |v| +
  C \kappa^{-N} \norm{Dv}_{C^0} \int e^{\epsilon r^{(N)}}.
  \end{equation}
But $\int e^{\epsilon r^{(N)}}=\int \Lp_\epsilon^N 1 \leq
e^{N\alpha(\epsilon)}$ by Lemma \ref{lem:ContractDebile}. In
the same way,
  \begin{equation}
  \int e^{\epsilon r^{(N)}} |v| \leq \norm{v}_{L^2} \left(\int
  e^{2\epsilon r^{(N)}} \right)^{1/2}\leq
  \norm{v}_{L^2}e^{N\alpha(2\epsilon)/2}.
  \end{equation}
We obtain (for some different function $\alpha(\epsilon)$)
  \begin{equation}
  \Lp_\epsilon^N |v|(x) \leq C e^{N\alpha(\epsilon)} \norm{v}_{L^2}+ C
  \kappa^{-N}e^{N\alpha(\epsilon)} \norm{Dv}_{C^0}.
  \end{equation}
Let us now bound $D(\Lp_\epsilon^N v)$. We can differentiate
$J^{(N)}(hx)$. As $\norm{D(J^{(N)} \circ h)(x)}\leq C J^{(N)}\circ
h$, we obtain a term which is bounded by $C\Lp_\epsilon^N |v|$. If
we differentiate $v\circ h(x)$, the resulting term is bounded by
  \begin{equation}
  \kappa^{-N} \sum J^{(N)}(hx) e^{\epsilon r^{(N)}(hx)}
  \norm{Dv}_{C^0}
  \leq C \kappa^{-N} \norm{Dv}_{C^0} \int e^{\epsilon
  r^{(N)}},
  \end{equation}
bounded by $C \kappa^{-N}
e^{N\alpha(\epsilon)}\norm{Dv}_{C^0}$. We have proved that
  \begin{equation}
  \norm{D(\Lp_\epsilon^N v)}_{C^0} \leq C
  \kappa^{-N}e^{N \alpha(\epsilon)}\norm{Dv}_{C^0} + C e^{N \alpha(\epsilon)}\norm{v}_{L^2}.
  \end{equation}
Taking $\epsilon_0$ small enough so that
$\kappa^{-1}e^{\alpha(\epsilon_0)} <1$, and $N_0$ large enough,
this implies the lemma.
\end{proof}

The following lemma essentially proves \eqref{eq:ContractC1Iteres}.
\begin{lem}
\label{lem:ControleC1Iteres} There exist $N_0>0$, $\theta_0<1$,
$C>0$, $\epsilon_0>0$ and a function $\alpha:(0,\epsilon_0)\to
\R_+$ tending to $0$ when $\epsilon\to 0$ such that, for any
$N\geq N_0$, for any $A\geq 1$, the following holds. Let
$\epsilon< \epsilon_0$, let $\psi_1,\dots,\psi_n \in
\CC^{A,\epsilon}_{N}$, let $v:Y\to\C$ be a $C^1$ function. Let
$v^0=v$ and $v^i=\Lp^{N}(\psi_i v^{i-1})$. Then
  \begin{equation}
  \norm{v^n}_{C^1} \leq CA\left(\prod_{i=1}^n \norm{\psi_i}_{\CC^{A,\epsilon}_{N}}\right)
  \left(\theta_0^{Nn}\norm{v}_{C^1}+e^{Nn\alpha(\epsilon)}\norm{v}_{L^2}\right).
  \end{equation}
\end{lem}
\begin{proof}
Note first that two points $x$ and $y$ of $Y$ can be joined by a
path of uniformly bounded length, since $\diam(Y)<\infty$. If $v$ is
a $C^1$ function, this implies $|v(x)|\leq C\norm{Dv}_{C^0}
+|v(y)|$. Integrating with respect to $y$,
  \begin{equation}
  \label{BorneC0}
  \norm{v}_{C^0} \leq C \norm{Dv}_{C^0} + \int |v|.
  \end{equation}

Let us first prove a preliminary inequality. For any $C^1$ function
$w$ and any integer $i$,
  \begin{equation}
  \norm{D(\Lp_\epsilon^{Ni} w)}_{C^0} \leq \theta_0^{Ni}
  \norm{Dw}_{C^0}+C e^{Ni\alpha(\epsilon)} \norm{w}_{L^2},
  \end{equation}
by Lemma \ref{lem:ContractLepsilon} (applied to the time $Ni$).
Applying \eqref{BorneC0} to $\Lp_\epsilon^{Ni} w$, we obtain
  \begin{equation}
  \norm{ \Lp_\epsilon^{Ni} w}_{C^0} \leq C \theta_0^{Ni} \norm{Dw}_{C^0} +
  C e^{Ni\alpha(\epsilon)}\norm{w}_{L^2}.
  \end{equation}
Let now $w$ be a Lipschitz function. It is a uniform limit of $C^1$
functions $w_n$, with $\norm{ Dw_n}_{C^0} \leq C \Lip(w)$. Taking
limits in the previous equation for $w_n$, we get
  \begin{equation}
  \norm{ \Lp_\epsilon^{Ni} w}_{C^0} \leq C \theta_0^{Ni} \Lip(w) +
  Ce^{Ni\alpha(\epsilon)} \norm{w}_{L^2}.
  \end{equation}
Let finally $v$ be a $C^1$ function. The function $|v|$ is
Lipschitz, and its Lipschitz coefficient is bounded by
$\norm{Dv}_{C^0}$. We conclude
  \begin{equation}
  \label{BorneC0Lip}
  \norm{ \Lp_\epsilon^{Ni}|v|}_{C^0} \leq C\theta_0^{Ni} \norm{Dv}_{C^0} +
  Ce^{Ni\alpha(\epsilon)}\norm{v}_{L^2}.
  \end{equation}

We can now prove the lemma itself. We will write $\gamma_i=
\norm{\psi_i}_{\CC^{A,\epsilon}_N}$. In particular,
$|\psi_i(x)|\leq \gamma_i e^{\epsilon r^{(N)}(x)}$. Hence,
$|v^i|\leq \gamma_i \dots \gamma_1 \Lp_\epsilon^{Ni}|v^0|$. As
$v^i(x)=\sum_{h\in \HH_N} J^{(N)}(hx) \psi_i(hx)v^{i-1}(hx)$,
we have
  \begin{align*}
  \norm{Dv^i(x)} \leq \gamma_i \Biggl(& \sum \norm{D(J^{(N)}\circ h)(x)}e^{\epsilon
  r^{(N)}(hx)}|v^{i-1}(hx)|  \\&+ \sum J^{(N)}(hx) A
  e^{\epsilon r^{(N)}(hx)}|v^{i-1}(hx)| \\& + \sum J^{(N)}(hx)
  e^{\epsilon r^{(N)}(hx)} \norm{Dh(x)} \norm{Dv^{i-1}(hx)}\Biggr).
  \end{align*}
We will bound these three terms. For the first one,
$\norm{D(J^{(N)}\circ h)(x)} \leq C J^{(N)}(hx)$. This term is
therefore bounded by $C \gamma_i\dots \gamma_1
\norm{\Lp_\epsilon^{Ni} |v^0|}_{C^0}$, which can be estimated with
\eqref{BorneC0Lip}. For the second term, we have a similar bound,
with an additional factor $A$.

For the third term, we bound $\norm{Dh(x)}$ by $\kappa^{-N}$,
and $\sum J^{(N)}(hx) e^{\epsilon r^{(N)}(hx)} =\Lp_\epsilon^N
1 (x) \leq e^{N \alpha(\epsilon)}$ by Lemma
\ref{lem:ContractDebile}. Taking $\epsilon$ small enough, we
can ensure that $\kappa^{-1}e^{\alpha(\epsilon)} \leq \theta_0$
(increasing $\theta_0$ if necessary).

We have proved that
  \begin{equation}
  \norm{Dv^i}_{C^0} \leq (1+A) \gamma_i\dots \gamma_1 (C\theta_0^{Ni} \norm{Dv}_{C^0} +
  Ce^{Ni\alpha(\epsilon)}\norm{v}_{L^2}) + \gamma_i \theta_0^N
  \norm{Dv^{i-1}}_{C^0}.
  \end{equation}
Iterating this equation inductively over $i$ yields
  \begin{align*}
  \norm{Dv^n}_{C^0}&\leq \left(\prod_{i=1}^n \gamma_i\right)
  \left( (1+A) \sum_{i=1}^n \theta_0^{N (n-i)} (C\theta_0^{Ni} \norm{Dv}_{C^0} +
  Ce^{Ni\alpha(\epsilon)}\norm{v}_{L^2}) + \theta_0^{Nn}
  \norm{Dv}_{C^0}
  \right)
  \\ &
  \leq \left(\prod_{i=1}^n \gamma_i\right) \left( C(1+A)n
  \theta_0^{Nn} \norm{Dv}_{C^0} + C e^{Nn\alpha(\epsilon)}
  \norm{v}_{L^2} + \theta_0^{Nn} \norm{Dv}_{C^0}\right)
  \\ &
  \leq C \left(\prod_{i=1}^n \gamma_i\right) \left( (1+A)
  \theta_0^{Nn/2} \norm{Dv}_{C^0} + C e^{Nn\alpha(\epsilon)}
  \norm{v}_{L^2} \right).
  \end{align*}
This gives the estimate of the lemma for $\norm{Dv^n}_{C^0}$. Thanks
to \eqref{BorneC0}, this also implies the desired bound for
$\norm{v^n}_{C^0}$.
\end{proof}

The following technical lemma will be needed later on.
\begin{lem}
\label{lem:SnPhiYOK} There exists a constant $\Cc>0$ such that, for
any $n\in \N$, for any $x\in Y$,
  \begin{equation*}
  \sum_{h\in\HH_n}J^{(n)}(hx) \norm{D(S_n^Y \phi_Y \circ h)(x)}^4  \leq
  \Cc^4.
  \end{equation*}
\end{lem}
\begin{proof}
If $h=h_n\circ \dots \circ h_1$, then $S_n^Y \phi_Y(x)=\sum_{i=1}^n
(\phi_Y\circ h_i)(h_{i-1}\dots h_1 x)$. Thus,
  \begin{equation}
  \norm{D (S_n^Y\phi_Y\circ h)(x)}^4 \leq C\left(\sum_{i=1}^n
  r(h_i\dots h_1 x)\kappa^{-i+1}\right)^4.
  \end{equation}
We will use the convexity inequality $(\sum a_i x_i)^4 \leq
(\sum a_i)^3 \sum a_i x_k^4$, which comes from the convexity of
$x\mapsto x^4$ when $\sum a_i=1$ (the general case can be
reduced to that specific case). We take $a_i=\kappa^{-i+1}$ and
$x_i=r(h_i\dots h_1 x)$, and obtain
  \begin{equation}
  \norm{D (S_n^Y\phi_Y\circ h)(x)}^4 \leq C \sum
  \kappa^{-i}r(h_i\dots  h_1 x)^4.
  \end{equation}
Let $F_n(x)=\sum_{h_1,\dots,h_n\in \HH} \left(\sum_{i=1}^n
\kappa^{-i}r(h_i\dots h_1 x)^4 \right) J^{(n)}(h_n\dots h_1
x)$, The sum that we want to estimate is bounded by $C F_n(x)$.
As $J^{(n)}(hx)\leq C J^{(n)}(hy)$ by Lemma
\ref{lem:ControleJacobien}, we have $F_n(x)\leq C F_n(y)$.
Hence, $F_n(x)\leq C \int F_n$. Finally, a change of variables
yields,
  \begin{equation}
  \int F_n = \sum_{i=1}^n \kappa^{-i}\int_Y r(T_Y^{n-i}x)^4\dd\mu_Y(x)
  =\sum_{i=1}^n \kappa^{-i} \int_Y r^4 \leq \frac{\int_Y r^4}{\kappa-1}
  \qedhere.
  \end{equation}
\end{proof}

\subsection{Contraction for Dolgopyat's norms}
\label{app:ContractDolgo}

To prove the contraction for Dolgopyat's norms, we will
essentially follow Dolgopyat's arguments as they are presented
in \cite[Section 7]{AGY:teich}, with additional technical
complications due to the facts that the involved functions are
unbounded, and that we want estimates which are uniform in $M$
in Theorem \ref{thm:MainContraction}.

We will need the following lemma, proved in \cite[Lemma
7.5]{AGY:teich}.
\begin{lem}
\label{lem:ExisteBump} There exist constants $\Ca>1$ and
$\Cb>0$ such that, for any ball $B(x,\Ca r)$ which is compactly
included in $Y$, there exists a $C^1$ function $\rho:Y\to
[0,1]$, vanishing outside $B(x,\Ca r)$, equal to $1$ on
$B(x,r)$ and with $\norm{\rho}_{C^1} \leq \Cb/r$.
\end{lem}

Later on, we will use oscillatory integral arguments. To do that, it
will be important that the phases of $e^{ik S^Y_N \phi_Y \circ h}$
vary at various speeds when one uses different inverse branches $h$.
This is ensured by the following lemma.
\begin{lem}
\label{lem:ExisteBranches} There exist $\Cd>0$ and an integer
$N_0>0$ such that, for any $N\geq N_0$, there exist inverse
branches $h_1,h_2\in \HH_N$ and a continuous unitary vector
field $y(x)$ on $Y$ such that, for any $x\in Y$,
  \begin{equation}
  |D(S^Y_N\phi_Y \circ h_1)(x)\cdot y(x) - D(S^Y_N\phi_Y \circ h_2)(x)\cdot y(x)|\geq
  \Cd.
  \end{equation}
\end{lem}
\begin{proof}
\emph{First step}. Let us show that there exist $C'$ and $N'$ such
that, for any $N\geq N'$, there exist inverse branches $h_1,h_2 \in
\HH_N$, a point $x\in Y$ and a unit tangent vector $y$ at $x$ such
that
  \begin{equation}
  \label{eq:ConvergePas0Champ}
  |D(S^Y_N\phi_Y \circ h_1)(x)\cdot y - D(S^Y_N\phi_Y \circ h_2)(x)\cdot y|>
  C'.
  \end{equation}
We argue by contradiction, so assume it is not the case.

Let us fix an inverse branch $h\in \HH$, and consider the
sequence of inverse branches $h^n$. Then $D(S_n^Y\phi_Y\circ
h^n)(x)\cdot y=\sum_{k=1}^n D(\phi_Y\circ h)(h^{k-1}x)
Dh^{k-1}(x)\cdot y$. As $\norm{D(\phi_Y\circ h)}$ is bounded
and $\norm{Dh^{k-1}(x)}\leq \kappa^{-k+1}$, this series
converges normally, to a continuous $1$-form $\omega(x)\cdot
y$. Let $x_0$ be any point in $Y$, the series
$\sum_{k=1}^\infty (\phi_Y\circ h^k -\phi_Y\circ h^k(x_0))$
even converges in $C^1$, and its sum $\psi$ is a $C^1$ function
with $D\psi=\omega$.

Let now $h'\in \HH$ be another inverse branch. Let us consider
$h_n=h^{n-1}\circ h' \in \HH_n$. Since we assume that
\eqref{eq:ConvergePas0Champ} does not hold, $D(S_n^Y \phi_Y
\circ h_n) - D(S_n^Y \phi_Y \circ h^n)$ converges pointwise to
$0$. But $D(S^Y_n \phi_Y\circ h_n)=D(\phi_Y\circ h')+
\sum_{k=1}^{n-1} D(\phi_Y\circ h)Dh^{k-1} Dh'$. Letting $n$
tend to infinity, we get
  \begin{equation}
  D\psi(x)\cdot y = D(\phi_Y\circ h')(x)\cdot y + D\psi(h'x)
  Dh'(x)\cdot y.
  \end{equation}
Hence, $D( (\phi_Y +\psi -\psi\circ T_Y)\circ h')=0$.
Therefore, the function $\phi_Y+\psi-\psi\circ T_Y$ is constant
on each set $h'(Y)$, $h'\in \HH$. This contradicts the fact
that $\phi_Y$ is not cohomologous to a locally constant
function, and concludes the proof of the first step.

\emph{Second step}. Let us fix an arbitrary branch $h\in \HH$.
Then $D(S^Y_p \phi_Y \circ h^p)=\sum_{k=0}^{p-1} D(\phi_Y\circ
h) Dh^k$ is uniformly bounded independently of $p$, by a
constant $c_0$. Fix $N\geq N'$ (given by the first step) such
that $c_0 \kappa^{-N}\leq C'/4$. Let $h_1$ and $h_2$ be the
inverse branches given by the first step, at time $N$, and let
$x_0$ and $y_0$ be a point in $Y$ and a tangent vector at this
point, satisfying the conclusions of the first step. We extend
$y_0$ to a continuous vector field on a neighborhood $U$ of
$x_0$, still satisfying \eqref{eq:ConvergePas0Champ}.

Since $\mu_Y$ has full support in $Y$, $\mu_Y(U)>0$. Hence, $U$
intersects $\bigcap_{k>0} \bigcup_{h\in \HH_k}h(Y)$, since
$\mu_Y$ is supported on this last set. Let $x_1$ be a point in
the intersection, and let $\ell_k\in \HH_k$ be the inverse
branch of $T_Y^k$ such that $x_1\in \ell_k(Y)$. Since the
diameter of $\ell_k(Y)$ tends to $0$ when $k\to\infty$,
$\ell_k(Y)$ is included in $U$ for large enough $k$. In
particular, there exist $k>0$ and an inverse branch $\ell \in
\HH_k$ such that $\ell(Y)\subset U$.

Let $y_1(x)=D\ell(x)^{-1}\cdot y_0(\ell x)$. For any $p\in \N$,
and $j\in \{1,2\}$, we have
  \begin{multline*}
  |D(S^Y_{p+N+k} \phi_Y \circ h^p\circ h_j \circ \ell)(x)\cdot y_1(x)
  -D(S^Y_{N+k}\phi_Y\circ h_j\circ \ell)(x)\cdot y_1(x)|
  \\
  =|D(S^Y_p \phi_Y \circ h^p)(h_j\ell x)Dh_j(\ell x)\cdot y_0(\ell x)|
  \leq c_0 \norm{Dh_j(\ell x)} \leq c_0\kappa^{-N}\leq C'/4.
  \end{multline*}
Moreover,
  \begin{multline*}
  |D(S^Y_{N+k}\phi_Y\circ h_1\circ \ell)(x)\cdot y_1(x)-
  D(S^Y_{N+k}\phi_Y\circ h_2\circ \ell)(x)\cdot y_1(x) |
  \\
  = |D (S_N^Y \phi_Y\circ h_1)(x)\cdot y_0(x) - D (S_N^Y \phi_Y \circ h_2)(x)\cdot
  y_0(x)|
  \geq C'.
  \end{multline*}
Adding these estimates, we obtain
  \begin{equation*}
  | D(S^Y_{p+N+k} \phi_Y \circ h^p\circ h_1 \circ \ell)(x)\cdot y_1(x)
  - D(S^Y_{p+N+k} \phi_Y \circ h^p\circ h_2 \circ \ell)(x)\cdot y_1(x)
  | \geq C'/2.
  \end{equation*}
We conclude the proof by taking $y(x)=y_1(x)/\norm{y_1(x)}$.
\end{proof}

We recall that we defined a constant $\Cc$ in Lemma
\ref{lem:SnPhiYOK}, and a constant $\Cd$ in Lemma
\ref{lem:ExisteBranches}.

\emph{We fix once and for all a constant $C_0 \geq
\max(4\Cc,10)$. We also fix an integer $N$ which is larger than
the integers $N_0$ given by Lemmas \ref{lem:ControleC1Iteres}
and \ref{lem:ExisteBranches}, and such that $\kappa^{-N} \leq
1/1000$ and $\Cd \geq 20 \kappa^{-N}C_0$.}

From this point on, the $D_k$ norms and the cones $\E_k$ will
always be defined with respect to the constant $C_0$. The
following lemma essentially proves \eqref{eq:PerdPasTrop}.

\begin{lem}
\label{lem:DolgoPerdPas} There exists a function
$\alpha:(0,\epsilon_0) \to \R_+$ which tends to $0$ when
$\epsilon$ tends to $0$ such that, for any
$\epsilon<\epsilon_0$, $M>0$ and $A>0$, there exists $K>0$ such
that, for any $|\ell|\geq |k|\geq K$, for any $C^1$ function
$v:Y\to\C$ and any function $\psi \in \CC^{A,\epsilon}_{MN}$,
  \begin{equation}
  \norm{\Lp_k^{MN} (\psi v)}_{D_\ell} \leq \norm{\psi}_{\CC^{A,\epsilon}_{MN}}
  e^{MN \alpha(\epsilon)}
  \norm{v}_{D_{2^M \ell}}.
  \end{equation}
\end{lem}
\begin{proof}
Let $u$ be such that $(u,v)\in \E_{2^M \ell}(C_0)$. Let
   \begin{equation}
  \tilde u=\norm{\psi}_{\CC^{A,\epsilon}_{MN}} \left(\sum_{h\in
  \HH_{MN}} J^{(MN)}(hx) u(hx)^2\right)^{1/2},
   \end{equation}
we will show that there exists $\alpha(\epsilon)$ (independent
of $M$) such that $(e^{MN\alpha(\epsilon)} \tilde u, \Lp_k^{MN}
(\psi v))\in \E_{\ell}(C_0)$.

We have
  \begin{equation}
  |\Lp^{MN}_k (\psi v)|\leq \sum_{h\in \HH_{MN}} J^{(MN)}(hx)
  \psi(hx)u(hx).
  \end{equation}
We bound $\psi(hx)$ by $\norm{\psi}_{\CC^{A,\epsilon}_{MN}}
e^{\epsilon r^{(MN)}(hx)}$, and use Cauchy-Schwarz inequality.
We conclude
  \begin{align*}
  |\Lp^{MN}_k (\psi v)| & \leq \norm{\psi}_{\CC^{A,\epsilon}_{MN}} \left(
  \sum J^{(MN)}(hx) e^{2\epsilon r^{(MN)}(hx)} \right)^{1/2} \cdot
  \left(\sum J^{(MN)}(hx) u(hx)^2 \right)^{1/2}
  \\ &
  = \Lp_{2\epsilon}^{MN}1(x)^{1/2} \cdot \tilde
  u(x).
  \end{align*}
The coefficient $\Lp_{2\epsilon}^{MN}1(x)^{1/2}$ is bounded by
a coefficient of the form $e^{MN\alpha(\epsilon)}$ by Lemma
\ref{lem:ContractDebile}.

Let us now estimate the derivative of
  \begin{equation}
  \Lp_k^{MN}(\psi v)(x)=\sum_{h\in \HH_{MN}} J^{(MN)}(hx) e^{-ik
  S^Y_{MN}\phi_Y(hx)} \psi(hx) v(hx).
  \end{equation}
If we differentiate $J^{(MN)}(hx)$, its derivative is bounded
by $C J^{(MN)}(hx)$ by Lemma \ref{lem:ControleJacobien}, and
the resulting term is therefore bounded by par $C
e^{MN\alpha(\epsilon)} \tilde u(x)$ as above. If we
differentiate $e^{-ik S^Y_{MN} \phi_Y(hx)}$, we use
Cauchy-Schwarz inequality and Lemma \ref{lem:SnPhiYOK} to
obtain a bound
  \begin{align*}
  & |k|\norm{\psi}_{\CC^{A,\epsilon}_{MN}}
  \left( \sum J^{(MN)}(hx) \norm{D (S^Y_{MN}\phi_Y\circ h)(x)}^4
  \right)^{1/4} \cdot \left( \sum J^{(MN)}(hx) e^{4\epsilon
  r^{(MN)}(hx)}\right)^{1/4}
  \\& \hspace{3cm} \cdot \left(\sum J^{(MN)}(hx) u(hx)^2
  \right)^{1/2}
  \\ & \leq \Cc |k| e^{MN\alpha(\epsilon)} \tilde u(x).
  \end{align*}
The derivative of $\psi\circ h$ is bounded by $A e^{\epsilon
r^{(MN)}(hx)} \norm{\psi}_{\CC^{A,\epsilon}_{MN}}$, and the
resulting term is therefore bounded by $A
e^{MN\alpha(\epsilon)} \tilde u(x)$. Finally, if we
differentiate $v(hx)$, we use the inequality $\norm{Dv(hx)}
\leq C_0 \kappa^{-MN} 2^M |\ell| u(hx)$, so that the resulting
term is bounded by $C_0 \kappa^{-MN} 2^M |\ell|
e^{MN\alpha(\epsilon)} \tilde u(x)$. Finally,
  \begin{equation}
  \norm{D(\Lp_k^{MN}(\psi v))(x)}\leq (C+A+ \Cc|k| + C_0 \kappa^{-MN} 2^M
  |\ell|) e^{MN\alpha(\epsilon)} \tilde u(x).
  \end{equation}
The choice of $N$ and $C_0$ implies that this term is bounded
by $C_0 |\ell| e^{MN\alpha(\epsilon)} \tilde u(x)$ if $K$ is
large enough.

Let us finally bound the derivative of $\tilde u$, or rather of
$\tilde u^2(x)=\norm{\psi}_{\CC^{A,\epsilon}_{MN}}^2 \sum
J^{(MN)}(hx) u(hx)^2$. If we differentiate the jacobian, the
resulting term is bounded by $C \tilde u^2$. If we differentiate
$u^2$, this is bounded by
  \begin{multline*}
  2\norm{\psi}_{\CC^{A,\epsilon}_{MN}}^2 \sum J^{(MN)}(hx)\kappa^{-MN} u(hx) \norm{Du(hx)}
  \\
  \leq 2\norm{\psi}_{\CC^{A,\epsilon}_{MN}}^2
  \kappa^{-MN} \cdot 2^M |\ell| C_0 \sum J^{(MN)}(hx) u(hx)^2
  = 2|\ell| 2^M\kappa^{-MN} C_0 \tilde u^2.
  \end{multline*}
Hence,
  \begin{equation}
  2\tilde u(x) \norm{D\tilde u(x)}=\norm{D\tilde u(x)^2} \leq
  2(C/2+ 2^M \kappa^{-MN} C_0 |\ell|) \tilde u(x)^2.
  \end{equation}
Dividing by $2\tilde u(x)$ and using $\kappa^{-N} \leq 1/1000$,
we obtain the desired bound $\norm{D\tilde u(x)}\leq C_0|\ell|
\tilde u(x)$ if $|\ell|$ is large enough.

We have proved that $(e^{MN\alpha(\epsilon)} \tilde u,
\Lp_k^{MN} (\psi v))\in \E_{\ell}(C_0)$. Hence,
  \begin{equation}
  \norm{\Lp_k^{MN} (\psi v)}_{D_\ell} \leq e^{MN\alpha(\epsilon)}
  \norm{\tilde u}_{L^4}
  \leq e^{MN\alpha(\epsilon)}
  \norm{\psi}_{\CC^{A,\epsilon}_{MN}} \norm{u}_{L^4}.
  \end{equation}
Taking the infimum over the quantities $\norm{u}_{L^4}$ for
$(u,v)\in \E_{2^M \ell}(C_0)$, we obtain the lemma.
\end{proof}

From this point on, we concentrate on the proof of
\eqref{eq:GagneL4}. For $v\in C^1(Y)$ and $\psi \in
\CC^{A,4\epsilon}_{MN}$, we will estimate $\Lp^{MN}_k(\psi v)$
by starting from $\psi v$ and applying $M$ times the operator
$\Lp_k^N$, which has good contraction properties thanks to the
phase compensation phenomenon given by Lemma
\ref{lem:ExisteBranches}. A technical issue in this argument is
the fact that the functions $\psi v, \Lp_k^N(\psi
v),\dots,\Lp^{(M-1)N}_k(\psi v)$ are not $C^1$ on $Y$, since
the function $\psi$ is quite wild at the beginning (it is only
bounded by $e^{4\epsilon r^{(MN)}(x)}$, so smoothness is only
regained after application of $\Lp_k^{MN}$). To deal with this
issue, we will introduce intermediate degrees of smoothness,
keeping track of the smoothness that has not yet been regained,
as follows.

If $Z$ is a subset of $Y$, $n\in\N$ and $\epsilon\geq 0$, we
will say that $(u,v)\in \E_k(C_0,Z,n,\epsilon)$ if the
functions $u$ and $v$ are $C^1$ on $Z$ and $|v|\leq e^{\epsilon
r^{(n)}}u$, $\norm{Du}\leq C_0 |k| u$ and $\norm{Dv} \leq C_0
|k| e^{\epsilon r^{(n)}} u$ on $Z$. In particular,
$\E_k=\E_k(C_0,Y,0,\epsilon)$ for any $\epsilon\geq 0$. We will
also write $\norm{v}_{D_k(Z,n,\epsilon)}$ for the infimum of
$\norm{u}_{L^4}$ over the functions $u$ such that $(u,v)\in
\E_k(C_0,Z,n,\epsilon)$.

\begin{lem}
\label{lem:ItereCones} There exists a function
$\alpha:(0,\epsilon_0)\to \R_+$ which tends to $0$ when
$\epsilon \to 0$ such that, for any $A>0$, $n>0$,
$\epsilon<\epsilon_0$, and for any $Z\subset Y$, there exists
$K>0$ such that, for any $|\ell|\geq |k|\geq K$, for any pair
of functions $(u,v)\in \E_{9\ell}(C_0, T_Y^{-N}Z,
nN,\epsilon)$, for any $C^1$ function $\chi: T_Y^{-N}Z \to
[3/4,1]$ with $\norm{D \chi}_{C^0} \leq |k|$ such that
$|\Lp_k^N v(x)|\leq \Lp^N(e^{\epsilon r^{(Nn)}}\chi u)(x)$,
holds
  \begin{equation}
  (e^{N\alpha(\epsilon)}\Lp^{N}(\chi^2 u^2)^{1/2}, \Lp_k^N v) \in
  \E_\ell(C_0, Z, (n-1)N, \epsilon).
  \end{equation}
\end{lem}
Note that the lemma also applies for $(u,v)\in \E_{\ell}(C_0,
T_Y^{-N}Z, nN,\epsilon)$ or $\E_{3\ell}(C_0, T_Y^{-N}Z,
nN,\epsilon)$, since these cones are contained in
$\E_{6\ell}(C_0, T_Y^{-N}Z, nN,\epsilon)$.
\begin{proof}[Proof of Lemma \ref{lem:ItereCones}]
The proof is similar to the proof of Lemma
\ref{lem:DolgoPerdPas}. One should only check that the
additional terms coming from the function $\chi$ are harmless
in the estimates. This is ensured by the choice of $N$ and
$C_0$.
\end{proof}

By Lemma \ref{lem:ExisteBranches}, we can fix two inverse branches
$h_1$ and $h_2$ of $T_Y^N$ as well as a vector field $y_0(x)$
satisfying the conclusion of the Lemma. Smoothing it, we obtain a
$C^1$ vector field $y$ such that $1\leq \norm{y} \leq 2$ and, for
any $x\in Y$,
  \begin{equation*}
  | D(S^Y_N\phi_Y\circ h_1)(x)\cdot y(x) - D(S^Y_N\phi_Y \circ
  h_2)(x)\cdot y(x)|\geq \Cd/2.
  \end{equation*}
Since $\norm{Dh_j(x)}\leq \kappa^{-N}$ and $\Cd\geq 20
\kappa^{-N}C_0$, this implies that
  \begin{equation*}
  | D(S^Y_N\phi_Y\circ h_1)(x)\cdot y(x) - D(S^Y_N\phi_Y \circ
  h_2)(x)\cdot y(x)| \geq 5 C_0 \max(\norm{D h_1(x)\cdot y(x)},
  \norm{Dh_2(x)\cdot y(x)}).
  \end{equation*}

Informally, this equation ensures that the difference between
the arguments of $e^{-ik S^Y_N \phi_Y(h_1 x)}$ and $e^{-ik
S^Y_N \phi_Y(h_{2} x)}$ varies quickly when $x$ moves slightly
in the direction of $y(x)$. Using this, it is possible to prove
the following lemma (see \cite[Lemma 7.13]{AGY:teich} for a
detailed proof):

\begin{lem}
\label{petite_boule_compense} There exist $\delta>0$ and
$\zeta>0$ satisfying the following property. Let $|k|\geq 10$
and $x_0\in Y$ be such that the ball $B=B(x_0,
(\zeta+\delta)/|k|)$ is compactly contained in $Y$. Consider
$(u,v)\in \E_{3k}(C_0,h_1B \cup h_2B,0,0)$. Then there exist
$x_1$ with $d(x_0,x_1)\leq \zeta/|k|$, and $j\in \{1,2\}$, such
that, for any $x\in B(x_1,\delta/|k|)$,
  \begin{multline*}
  |e^{-ik S^Y_N \phi_Y(h_j x)} J^{(N)}(h_jx)v(h_jx)+ e^{-ik S^Y_N
  \phi_Y(h_{2-j} x)} J^{(N)}(h_{2-j} x)v(h_{2-j} x)| \\
  \leq \frac{3}{4} J^{(N)}(h_j x)u(h_j x)+ J^{(N)}(h_{2-j}x)
  u(h_{2-j}x).
  \end{multline*}
\end{lem}

If $H$ is a set of inverse branches of $T_Y^n$, we will write
$H(Y)=\bigcup_{h\in H} h(Y)$.

\begin{lem}
\label{lem:DolgoClassique} There exist $\theta_1<1$ and a
function $\alpha:(0,\epsilon_0)\to \R_+$ tending to $0$ when
$\epsilon\to 0 $ satisfying the following property. Let $n>0$,
let $H$ be a finite subset of $\HH_{nN}$. Denote by $H^{(n-1)N}
\subset\HH_{(n-1)N}$ the set of inverse branches $T_Y^N\circ h$
for $h\in H$. Then, for any $H$, there exists $K(H)$ such that,
for any $|k|\geq K(H)$, for any function $v$, for any
$\epsilon<\epsilon_0$,
  \begin{equation}
  \norm{\Lp_k^N v}_{D_k(H^{(n-1)N}(Y),\epsilon, (n-1)N)} \leq \theta_1^N e^{N\alpha(\epsilon)}
  \norm{v}_{D_{3k}(H(Y),\epsilon,nN)}.
  \end{equation}
\end{lem}
\begin{proof}
Increasing $H$ if necessary, we can assume that, for any $h\in
H^{(n-1)N}$, the branches $h_1 \circ h$ and $h_2 \circ h$ belong to
$H$. Let $(u,v)\in \E_{3k}(C_0, H(Y),\epsilon,nN)$.

Let $h\in H^{(n-1)N}$, we will work on $h(Y)$, and use the weak
Federer property for the constant $C=\Ca(\zeta/\delta+1)$
(where $\Ca$ is given by Lemma \ref{lem:ExisteBump}).
Definition \ref{def:WeakFedererProperty} provides us with
constants $D>0$ and $\eta_0(h(Y),C)$. Since the weak Federer
property is uniform over the inverse branches of $T_Y$, we can
even choose $D$ depending only on $C$, and not on $h$.

We apply the definition of the weak Federer property to
$\eta=\delta/(\Ca|k|)$. If $|k|$ is large enough, we indeed
have $\eta < \eta_0(h(Y),C)$ for any $h\in H^{(n-1)N}$ (here,
the finiteness of $H$ is crucial). We obtain disjoint balls
$B(x_1, \Ca(\zeta/\delta+1)\eta),\dots, B(x_k,
\Ca(\zeta/\delta+1) \eta)$ compactly contained in $h(Y)$, and
sets $A_1,\dots,A_k$ contained in $B(x_i, D \eta)$, whose union
covers $h(Y)$, and such that, for any $x'_i\in B(x_i,
(\Ca(\zeta/\delta+1)-1) \eta)$, holds $\mu_Y(B(x'_i, \eta))
\geq \mu_Y(A_i)/D$.

On each ball $B=B(x_i, \Ca(\zeta/\delta+1) \eta)=B(x_i,
(\zeta+\delta)/|k|)$, we apply Lemma
\ref{petite_boule_compense} to the pair of functions
$(u(x)e^{\epsilon r^{(nN)}(x)}, v(x))$ (which belongs to
$\E_{3k}(C_0, T_Y^{-N}B,0,0)$). The conclusion of this lemma
gives a ball $B'_i=B(x'_i, \delta/|k|)$ as well as an index
$j\in \{1,2\}$. We will write $\type(B'_i)=j$. Let
$B''_i=B(x'_i, \delta/(\Ca k))=B(x'_i, \eta)$. By Lemma
\ref{lem:ExisteBump}, there exists a function $\rho_i$ equal to
$1$ on $B''_i$, vanishing outside of $B'_i$, whose $C^1$ norm
is bounded by $C|k|$.

Let us then define a function $\rho$ on $T_Y^{-N}(hY)$ by
$\rho=(\sum_{\type(B'_i)=j} \rho_i)\circ T_Y^N$ on $h_j(hY)$
(for $j=1,2$) and $\rho=0$ elsewhere. Finally, let $\chi=1-c
\rho$ where $c$ is small enough. Then $\norm{\chi}_{C^1}\leq
|k|$ if $c$ is small enough, and $|\Lp_k^N v|\leq \Lp^N(\chi
ue^{\epsilon r^{(nN)}})$ by construction (using Lemma
\ref{petite_boule_compense}). Hence, Lemma \ref{lem:ItereCones}
implies that $(e^{N\alpha(\epsilon)}\Lp^N(\chi^2 u^2)^{1/2},
\Lp^N_k v)\in \E_k(C_0, h(Y), (n-1)N,\epsilon)$.

We glue together the different functions $\chi$ obtained by
varying $h$, to obtain a function (that we still denote by
$\chi$) on $H(Y)$. We sill have
$(e^{N\alpha(\epsilon)}\Lp^N(\chi^2 u^2)^{1/2}, \Lp^N_k v) \in
\E_k(C_0, H^{(n-1)N}(Y), (n-1)N,\epsilon)$. If we can prove
that $\norm{\Lp^N(\chi^2 u^2)^{1/2}}_{L^4} \leq \beta
\norm{u}_{L^4}$ where $\beta<1$ is a constant which is
independent of everything else, then the proof will be
finished.

Let $\tilde u= \Lp^N(\chi^2 u^2)^{1/2}$. We have
  \begin{equation*}
  \tilde u(x)^4 = \left(\sum_{h\in \HH_n} J^{(N)}(hx) \chi(hx)^2
  u(hx)^2\right)^2
  \leq \left(\sum_{h\in \HH_N} J^{(N)}(hx) \chi(hx)^4\right)\cdot
  \left( \sum_{h\in \HH_N} J^{(N)}(hx) u(hx)^4 \right).
  \end{equation*}
Let $Y_1=\bigcup B''_i$, and let $Y_2$ be its complement. On $Y_1$,
the factor $\sum_{h\in \HH_N} J^{(N)}(hx) \chi(hx)^4$ is bounded by
a uniform constant $\beta_0<1$, hence $\tilde u(x)^4 \leq \beta_0
\Lp^N (u^4)(x)$. On $Y_2$, we only have $\tilde u(x)^4 \leq
\Lp^N(u^4)(x)$.

Let $w=\Lp^N(u^4)$. Since $\norm{Du}\leq 3 C_0|k| u$, there
exists a constant $C$ such that $\norm{Dw}\leq C|k|w$.
Integrating this inequality along a path between two points
yields $w(x) \leq e^{C|k| d(x,y)} w(y)$ for any $x,y$. In
particular, since $A_i\subset B(x_i, CD\delta/(C_2 |k|))$,
there exists $C$ such that, for any $x\in A_i$ and $y\in
B''_i$, we have $w(x)\leq C w(y)$. Integrating this inequality,
  \begin{equation*}
  \frac{\int_{A_i}w}{\mu_Y(A_i)} \leq C \frac{\int_{B''_i}
  w}{\mu_Y(B''_i)}.
  \end{equation*}
But $\mu_Y(A_i)\leq D \mu_Y(B''_i)$ by definition of the sets
$A_i$, hence $\int_{A_i} w \leq C \int_{B''_i} w$. The balls
$B''_i$ are pairwise disjoint, so we conclude $\int_{Y_2} w
\leq C' \int_{Y_1} w$ for some constant $C'$.

Let $E$ be large enough so that $(E+1)\beta_0 + C' \leq E$. Then
  \begin{equation*}
  (E+1)\int \tilde u^4 \leq (E+1)\int_{Y_1} \beta_0 w +
  (E+1)\int_{Y_2} w
  \leq (E+1)\beta_0 \int_{Y_1} w + E\int_{Y_2} w + C' \int_{Y_1}w
  \leq E \int w.
  \end{equation*}
Hence, $\norm{\tilde u}_{L^4}^4 \leq \frac{E}{E+1} \int
w=\frac{E}{E+1} \int u^4$. This is the desired inequality.
\end{proof}

\begin{lem}
\label{lem:PresqueFini} There exist $\theta_2<1$ and a function
$\alpha : (0,\epsilon_0)\to \R_+$ which tends to $0$ when
$\epsilon\to 0$ satisfying the following property. For any
$M>0$, $\epsilon<\epsilon_0$ and $A>0$, there exists $K>0$ such
that, for any $C^1$ function $v:Y\to \C$ and for any $\psi\in
\CC^{A,\epsilon}_{MN}$, for any $|k|\geq K$,
  \begin{equation}
  \norm{\Lp_k^{MN}(\psi v)}_{D_k} \leq e^{MN\alpha(\epsilon)}
  \theta_2^{MN}\norm{\psi}_{\CC^{A,\epsilon}_{MN}} \norm{v}_{D_{2^M k}}.
  \end{equation}
\end{lem}
\begin{proof}
We will give the proof for odd $M$ (the proof for even $M$ is
analogous and even simpler).

We will decompose  $\HH_{MN}$ as the union of a finite set
$H_1$ (to which we will apply Lemma \ref{lem:DolgoClassique})
and a set $H_2$ which will yield a small enough contribution.
Let $H \subset \HH$ have finite complement. We will take for
$H_1$ the set of inverse branches in $\HH_{MN}$ which are the
composition of branches not belonging to $H$, and for $H_2$ its
complement.

Let $w=1_{H_1(Y)} \psi v$ and $w'=1_{H_2(Y)}\psi v$. We will
first estimate $\norm{\Lp^{MN}_k w'}_{D_k}$. Let $u$ be such
that $(u,v)\in \E_{2^M k}$. Let $\tilde
u=\norm{\psi}_{\CC^{A,\epsilon}_{MN}} \left(\sum_{h\in H_2}
J^{(MN)}(hx) u(hx)^2\right)^{1/2}$, the computation made in the
proof of Lemma \ref{lem:DolgoPerdPas} shows that
$(e^{MN\alpha(\epsilon)}\tilde u, \Lp_k^{MN}w')\in \E_k(C_0)$.

We have
  \begin{equation}
  \tilde u^2 \leq \norm{\psi}_{\CC^{A,\epsilon}_{MN}}^2
  \sum_{a+b= MN-1} \Lp^a \Lp_{0, H} \Lp^b
  u^2,
  \end{equation}
where $\Lp_{0,H}$ is similar to the operator $\Lp$, but the sum is
only done over branches belonging to $H$ (this operator has already
been defined before Lemma \ref{lem:ContractDebile}). This lemma
shows that, if $H$ is chosen small enough, then
$\norm{\Lp_{0,H}}_{L^2 \to L^2}$ can be made arbitrarily small.
Hence, if $H$ is small enough (in terms of $M$ and $\epsilon$), we
have
  \begin{equation}
  \label{eq:ContractGraceTroncature}
  \norm{\Lp^{MN}_k w'}_{D_k} \leq (\theta_1^{MN/3}-\theta_1^{MN/2})
  \norm{\psi}_{\CC^{A,\epsilon}_{MN}}\norm{v}_{D_{2^M k}}.
  \end{equation}

Let us fix such an $H$. Since $M$ is odd, it can be written as
$M=2m+1$. The set $H_1$ is finite and fixed. In particular, there
exists a constant $B$ such that, for any $x\in H_1(Y)$, $\norm{D
\psi(x)}\leq B \norm{\psi}_{\CC^{A,\epsilon}_{MN}}$. If $|k|$ is
large enough (in terms of $B$), this yields
  \begin{equation}
  \norm{w}_{D_{3^M k}(H_1(Y), MN,\epsilon)}\leq
  \norm{\psi}_{\CC^{A,\epsilon}_{MN}} \norm{v}_{D_{2^M k}}.
  \end{equation}
Iterating $m$ times Lemma \ref{lem:ItereCones} (with $\chi=1$), we
obtain
  \begin{equation}
  \norm{\Lp_k^{mN} w}_{D_{3k}(H_1^{(m+1)N}(Y), (m+1)N, \epsilon)} \leq e^{mN
  \alpha(\epsilon)} \norm{\psi}_{\CC^{A,\epsilon}_{MN}} \norm{
  v}_{D_{2^M k}}.
  \end{equation}
We then apply inductively Lemma \ref{lem:DolgoClassique}. If $|k|$
is large enough, we obtain for $i>m$
  \begin{equation}
  \norm{\Lp^{iN}_k w}_{D_k(H_1^{(M-i)N}(Y), (M-i)N,\epsilon)} \leq
  e^{i N\alpha(\epsilon)} \norm{\psi}_{\CC^{A,\epsilon}_{MN}} \theta_1^{(i-m)N}
  \norm{v}_{D_{2^M k}}.
  \end{equation}
For $i=M=2m+1$, we conclude
  \begin{equation}
  \label{eq:ContractGraceDolgo}
  \norm{\Lp^{MN}_k w}_{D_k} \leq e^{MN\alpha(\epsilon)}
  \norm{\psi}_{\CC^{A,\epsilon}_{MN}} \theta_1^{MN/2} \norm{v}_{D_{2^M
  k}}.
  \end{equation}

Adding up the inequalities \eqref{eq:ContractGraceTroncature} and
\eqref{eq:ContractGraceDolgo}, we get the conclusion of the lemma.
\end{proof}

\begin{proof}[Proof of Theorem \ref{thm:MainContraction}]
We choose $\theta\in (2^{-1/(1010N)},1)$ such that
$\theta^{100}$ is larger than the constants $\theta_0$ given by
Lemmas \ref{lem:contractSKJFDL} and \ref{lem:ControleC1Iteres},
and than $\theta_2$ given by Lemma \ref{lem:PresqueFini}. If
$\epsilon>0$ is small enough, Lemma \ref{lem:ControleC1Iteres}
(applied to $MN$) shows \eqref{eq:ContractC1Iteres}. Moreover,
\eqref{eq:ContracteTriviale} is implied by Lemma
\ref{lem:contractSKJFDL}. Finally, \eqref{eq:PerdPasTrop} is a
consequence of Lemma \ref{lem:DolgoPerdPas}, and
\eqref{eq:GagneL4} follows from Lemma \ref{lem:PresqueFini}.
\end{proof}

\bibliography{biblio}
\bibliographystyle{alpha}

\end{document}